\documentclass[11pt]{amsart} 
\calclayout

\usepackage{setspace} 
\usepackage{graphicx}
\usepackage{microtype}
\usepackage{bbm}

\usepackage{pgfplots}
\pgfplotsset{compat=1.15}
\usepackage{tikz}
\usepackage{mathrsfs}
\usetikzlibrary{arrows}

\usepackage[english]{babel}
\usepackage{csquotes}
\usepackage{amssymb}
\usepackage{amsmath}
\usepackage{amsthm}
\usepackage{float}

\usepackage{stmaryrd}
\usepackage{slashed}
\usepackage{tikz-cd}
\usetikzlibrary{positioning}
\usetikzlibrary{calc,intersections}
\usepackage{hyperref}
\usepackage{todonotes}
\usepackage[all,cmtip]{xy}

\usepackage[justification=centering]{caption}

\newtheorem{thm}{Theorem}[section]
\newtheorem{lem}[thm]{Lemma}
\newtheorem{prop}[thm]{Proposition}

\theoremstyle{remark}
\newtheorem{rema}[thm]{Remark}
\newtheorem{exa}[thm]{Example}
\newtheorem{nota}[thm]{Notation}
\newtheorem{defi}[thm]{Definition}

\DeclareFontFamily{U} {MnSymbolA}{}
\DeclareFontShape{U}{MnSymbolA}{m}{n}{
  <-6> MnSymbolA5
  <6-7> MnSymbolA6
  <7-8> MnSymbolA7
  <8-9> MnSymbolA8
  <9-10> MnSymbolA9
  <10-12> MnSymbolA10
  <12-> MnSymbolA12}{}
\DeclareFontShape{U}{MnSymbolA}{b}{n}{
  <-6> MnSymbolA-Bold5
  <6-7> MnSymbolA-Bold6
  <7-8> MnSymbolA-Bold7
  <8-9> MnSymbolA-Bold8
  <9-10> MnSymbolA-Bold9
  <10-12> MnSymbolA-Bold10
  <12-> MnSymbolA-Bold12}{}

\DeclareSymbolFont{MnSyA} {U} {MnSymbolA}{m}{n}
\DeclareMathSymbol{\lcirclearrowright}{\mathrel}{MnSyA}{252}
\DeclareMathSymbol{\rcirclearrowleft}{\mathrel}{MnSyA}{250}

\DeclareMathOperator{\pr}{pr}
\DeclareMathOperator{\id}{Id}

\DeclareMathOperator{\dec}{Dec}
\DeclareMathOperator{\Hom}{Hom}
\DeclareMathOperator{\sym}{Sym}
\DeclareMathOperator{\lfun}{Lie}

\newcommand{\cA}{\mathcal{A}}

\newcommand{\cC}{\mathcal{C}}
\newcommand{\cD}{\mathcal{D}}

\newcommand{\cF}{\mathcal{F}}
\newcommand{\cG}{\mathcal{G}}
\newcommand{\cH}{\mathcal{H}}

\newcommand{\cK}{\mathcal{K}}
\newcommand{\cL}{\mathcal{L}}

\newcommand{\cV}{\mathcal{V}}
\newcommand{\cW}{\mathcal{W}}

\newcommand{\g}{\mathfrak{g}}
\newcommand{\h}{\mathfrak{h}}
\newcommand{\fk}{\mathfrak{k}}

\newcommand{\fX}{\mathfrak{X}}

\renewcommand{\gg}{\mathfrak{g}}            				       
\newcommand{\G}{\mathcal{G}}									   
\newcommand{\hh}{\mathfrak{h}}            						   
\newcommand{\XX}{\mathfrak{X}}                                     
\newcommand{\Rr}{\mathbb{R}}                                       
\newcommand{\F}{\mathcal{F}}									   
\newcommand{\A}{\mathcal{A}}									   
\renewcommand{\H}{\mathcal{H}}									   

\newcommand{\Lie}{\mathcal{L}}                                     
\providecommand{\abs}[1]{\lvert#1\rvert}						   
\newcommand{\dd}[1]{\frac{d}{d #1}\Big{\vert}_{#1 =0}}			   

\newcommand{\vJoin}{\cdot} 

\let\oldtocsection=\tocsection

\let\oldtocsubsection=\tocsubsection

\let\oldtocsubsubsection=\tocsubsubsection

\renewcommand{\tocsection}[2]{\bf\hspace{0em}\oldtocsection{#1}{#2}}
\renewcommand{\tocsubsection}[2]{\hspace{1em}\oldtocsubsection{#1}{#2}}
\renewcommand{\tocsubsubsection}[2]{\hspace{2em}\oldtocsubsubsection{#1}{#2}}

\title{The van Est homomorphism for strict Lie 2-groups}

\author{Camilo Angulo}

\address{School of Mathematics, Jilin University, No. 2699 Qianjin Street, Changchun City, Jilin Province, China}

\email{cangulo@jlu.edu.cn}

\author{Miquel Cueca}

\address{Departement of Mathematics, KU Leuven. Celestijnenlaan 200B,  Leuven (Heverlee), B-3001, Belgium}

\email{miquel.cuecaten@kuleuven.be}


\date{} 

\begin{document}

\begin{abstract} 
We construct a van Est map for strict Lie $2$-groups from the Bott-Shulman-Stasheff double complex of the strict Lie $2$-group to the Weil algebra of its associated strict Lie $2$-algebra.  We show that, under appropriate connectedness assumptions, this map induces isomorphisms in cohomology. As an application,  we differentiate the Segal $2$-form on the loop group.
\end{abstract}

\maketitle

\tableofcontents

\section{Introduction}

Since their introduction in \cite{cou:cou}, Courant algebroids have become a fundamental object in mathematics and physics. In \cite{roy:sym}, through the lens of supergeometry, Courant algebroids are reinterpreted as Lie $2$-algebroids endowed with a \emph{symplectic} structure. Using this point of view, it is proposed in  \cite{sev:some} that Courant algebroids are the infinitesimal counterpart of \emph{symplectic Lie $2$-groupoids}. Today one canh find many attempts to formalize such an integration in the literature, see e.g. \cite{Cueca_Zhu:2021, Lesdiablerets, sev:cou, raj:from, cc:cou, sev:int}; nevertheless, a concrete differentiation procedure is conspicuously lacking. 

To better understand this problem, in this article, we consider the toy model of a strict Lie 2-group and achieve two goals:
\begin{itemize}
    \item Firstly, given a strict Lie 2-group $G$, we outline a differentiation procedure depending only on the nerve of $\overline{W}G$, which is a priori a simplicial manifold and not necessarily a simplicial group, and produces a degree 2 $NQ$-manifold dual to its strict Lie 2-algebra, see Theorem~\ref{Theo:SimplicialDiffn}. Thus, this has the potential to be generalized to differentiate higher Lie groupoids in general, see also \cite{hen:int, cc:dif, sev:dif}.
    \item Secondly, we relate the simplicial de Rham complex, also known as the \emph{Bott-Shulman-Stasheff double-complex},  of the strict Lie $2$-group with the Weil algebra of its Lie $2$-algebra via a van Est map, see Theorem~\ref{Theo:CommutesWithDiffs}. This is an extension of the classical van Est map for Lie groups, see e.g. \cite{abad:VE, Meinrenken_Li-Bland:2015, mein:VEint, ve:group}.
\end{itemize}

As a direct consequence of our main achievements, in $\S\ref{sec-Applications}$ we differentiate the $2$-shifted symplectic Lie $2$-group $(\mathbb{G}_\bullet, \omega_S)$ constructed in \cite[$\S 3$]{Cueca_Zhu:2021} and show its relation to the quadratic Lie algebra $(\g, [\cdot,\cdot], \langle\cdot,\cdot\rangle)$; thus, producing a concrete differentiation that gives back a Courant algebroid over a point. 
 
There are a number of van Est map formulas scattered throughout the literature, see e.g. \cite{Angulo:2022, abad:VE,  ale:ve, alan:ex}. As stated in \cite{Meinrenken_Li-Bland:2015}, these formulas can be given a more ``conceptual explanation'' by deducing them as a consequence of the Perturbation Lemma of homological algebra \cite{Gugenheim:1972}; furthermore, this viewpoint has the advantage that ``the basic properties of the van Est map follow rather easily''. 
Let us briefly summarize the construction of the classical van Est map and spot the main challenges of our work. 

For a Lie group $G$ it is known that its \emph{d\'{e}calage} $\dec_\bullet(G)$ gives a simplicial model of the universal principal $G$-bundle $\partial_0:E_\bullet G\longrightarrow B_\bullet G$. 
The tangent bundle to the fibers of this bundle defines a simplicial object in the category of Lie algebroids 
$\F_\bullet=\ker T\partial_0\subset T(E_\bullet G)$. 
Dualizing, one obtains a double complex $$C^{p, q}=\Gamma(\bigwedge\nolimits^q\F^*_p)$$ whose columns are Lie algebroid cochain 
complexes with their Chevalley-Eilenberg differentials $\delta$ and whose rows coincide with groupoid cochain 
complexes with values in different representations with their natural simplicial differentials 
$\partial$. Pullback along the projection map yields a map of complexes
$$\partial_0^*:C^\infty(G_\bullet)\longrightarrow C^{\bullet,\bullet}$$
and evaluation at units gives a map
$\iota_0^*:C^{\bullet,\bullet}\to\bigwedge^\bullet\gg^*$ dual to the inclusion $\gg\hookrightarrow TG=\F_0$; however, since this inclusion is in general not a Lie algebroid morphism, the map $\iota_0^*$ is not a cochain map, in general. 

The simplicial manifold $E_\bullet G$ admits a canonical simplicial deformation 
retraction to a point. This determines a homotopy operator $h$ for the simplicial 
differential $\partial$ on the double complex $C^{p,q}$. It follows from 
the Perturbation Lemma of homological algebra that the composition 
     $$\iota_0^*\circ (1+\delta\circ h)^{-1}:{\textnormal{Tot}}(C^{\bullet,\bullet})\longrightarrow\bigwedge\nolimits^\bullet\gg^*$$ 
is a cochain map. The composition $\iota_0^*\circ (1+\delta\circ h)^{-1}\circ\partial_0^*:C^\infty(G_\bullet)\to\bigwedge\nolimits^\bullet\gg^* $ is then taken as the definition of the van Est map.
We point out that this process depends only on the simplicial structure of the nerve of the Lie group and, interestingly, that the Lie algebra structure can be recovered as the dual structure to the first page of the spectral sequence of the filtration by columns of $C^{p,q}$ that first computes the horizontal cohomology. 

We aim at reproducing the above argument for strict Lie $2$-groups to obtain the novel van Est map. In doing so, the main difficulty lies in that for a strict Lie 2-group $G_\bullet$, at first sight, there is no \emph{natural} simplicial Lie 2-algebroid, as opposed to the simplicial Lie algebroid $\F_\bullet$ above. In this sense, the main contribution of this article is the auxiliary simplicial Lie 2-algebroid constructed in Proposition~\ref{prop:A}, which relies on the existence of the {\it mysterious} simplicial splitting of Proposition~\ref{prop:eps}. This simplicial Lie $2$-algebroid can be understood as a $2$-step differentiation using Diagram \eqref{dia:dec-w}.

With this additional structure at hand, the main difficulty is to get an explicit formula for the van Est map thus defined. The outline of the construction of the van Est map and the deduction of its properties follow along the same lines as they do in the group case. The main tool of this part was to construct a homotopy, a task that we perform in $\S\ref{sec:hom}$.  

In general, in a Lie 2-groupoid (see Definition~\ref{def:lie-n-gpd}), the map $(\partial_2,\partial_0):X_2\to X_1{}_{\partial_0}\times_{\partial_1}X_1$  is not only a surjective submersion, but, in fact, a bigon-principal bundle (see \cite{z:tgpd-2} for details). It is perhaps worth mentioning here that the simplicial splitting of Proposition \ref{prop:eps} and the homotopy in $\S\ref{sec:hom}$ admit an immediate extension as long as $X_2$ is trivial as a principal bundle, in which case, we call the Lie 2-groupoid \emph{strict}. The results of this work extend thus to strict Lie 2-groupoids, although we opted out of adding this to the exposition, as the formulas involving the differential of the Weil algebra of the associated algebroid would require choosing a connection that makes the bookkeeping of Section~\ref{sec:VE} a lot more cumbersome. The future challenge is understanding how the splitting of Proposition~\ref{prop:eps} can be replaced in the non-strict case.

 In Section~\ref{sec:pre} we lay down preliminaries and conventions. Section \ref{sec:dif} explains how to differentiate strict Lie $2$-groups anew. Section~\ref{sec:VE} is the heart of the article, in Subection~\ref{sec:hom} we introduce the main tools necessary to apply the Perturbation Lemma, while in Subsections \ref{sec:func} and \ref{sec:proof} we use it to deduce the formula for the van Est map at the level of cochains and at the level of forms, respectively. In Subsection~\ref{subsec-vanEstIso}, we explain how our main tool can be used to get a van Est-type theorem relating the cohomologies of a Lie 2-group and its Lie 2-algebra.  In Section~\ref{sec-Applications}, we use the van Est map to differentiate a 3-shifted symplectic structure and the famous Segal form. Lastly, Appendix \ref{sec:ap} deals with the technical computations that assert that the homotopies we use throughout are so.

\vspace{3mm}

\noindent{\bf Acknowledgements.}
We thank Chenchang Zhu for her encouragement and support, and for bringing our attention to the d\'{e}calage construction. We also thank Alejandro Cabrera and Kateryna Hlyniana for fruitful discussions and suggestions on previous drafts of this paper and all the members of the Higher Structure Seminar at G\"{o}ttingen for inspiring conversations. Lastly, we thank the anonymous Referees for helping in the improvement of the article. C.A. would like to thank the hospitality of Georg-August Universit\"{a}t at G\"{o}ttingen where some of this research was carried out.

\section{Preliminaries}\label{sec:pre}

In this section, we recall basic constructions about strict Lie $n$-groups and their infinitesimal counterparts. The content of this section is standard material, for more detailed expositions see e.g. \cite{goe:book, wol:des, z:tgpd-2}.

\subsection{Simplicial objects in different categories}

Recall that  a \emph{simplicial object $X_\bullet$ on a category $\cC$} is a contravariant
functor from $\Delta$, the category of finite ordinals 
\[
[0]=\{0\}, \qquad [1]=\{0, 1\},\quad \dotsc,\quad
[l]=\{0, 1,\dotsc, l\},\quad\dotsc,
\]
with order-preserving maps, to $\cC$. More precisely, $X_\bullet$ consists of $X_l$ objects in $\cC$, \emph{face} maps $d^l_k: X_l \to X_{l-1}$ and \emph{degeneracy} maps $s^l_k: X_l \to X_{l+1}$, for $k=0, \dots, l$ satisfying the following simplicial identities
\begin{equation*}\label{eq:face-degen}
    \begin{array}{lll}
        d^{l-1}_i d^{l}_j =& d^{l-1}_{j-1} d^l_i &\text{if}\; i<j,  \\
       s^{l}_i s^{l-1}_j =& s^{l}_{j+1} s^{l-1}_i & \text{if}\; i\leq j,
    \end{array}\qquad  d^l_i s^{l-1}_j =\left\{\begin{array}{ll}
    s^{l-2}_{j-1} d^{l-1}_i  & \text{if}\; i<j, \\
    \id  & \text{if}\; i=j, j+1,\\
    s^{l-2}_j d^{l-1}_{i-1} & \text{if}\; i> j+1.
 \end{array}\right.
\end{equation*}
Later in our article, when the context is clear, we drop the upper indices and write instead $d_i$ and $s_i$ for simplicity.  A \emph{morphism $\Phi_\bullet:X_\bullet\to Y_\bullet$} is a family of maps  \(\Phi_n:X_n\to Y_n\) on $\cC$ that intertwine the face and degeneracy maps on \(X_\bullet\) and~\(Y_\bullet\). For more details see e.g. \cite{goe:book, wei:book}.

Among others, we are interested in \emph{simplicial sets} when $\cC$ is the category of sets, \emph{simplicial manifolds} when $\cC$ is the category of smooth manifolds,  \emph{simplicial (Lie) groups} when $\cC$ is the category of (Lie) groups, and \emph{simplicial Lie algebras} when $\cC$ is the category of Lie algebras.

    Let $X_\bullet$ be a simplicial manifold. An important consequence of the simplicial identities is that $C^\infty(X_\bullet)$ inherits the structure of a cochain complex with differential $\partial:C^\infty(X_n)\to C^\infty(X_{n+1})$ given by 
    \begin{equation}\label{eq:simfun}
        \partial f=\sum_{i=0}^n (-1)^id_i^*f.
    \end{equation}
    Another consequence is that  $\widehat{C}^\infty(X_\bullet)=\{ f\in C^\infty(X_\bullet)\ | \ s^*_if=0\ \forall i\}$ is a quasi-isomorphic subcomplex known as the \emph{normalized subcomplex} (see e.g. \cite{behcohsta}).

The following simplicial sets play an important r\^{o}le for us, they are called the \emph{$l$-simplex} $\Delta[l]$ and
the \emph{horn} $\Lambda[l,j]$:
\begin{equation*}\label{eq:simplex-horn}
\begin{split}
(\Delta[l])_k & = \{ f: [k] \to [l] \mid f(i)\leq
f(j),
\forall i \leq j\}, \\
(\Lambda[l,j])_k & = \{ f\in (\Delta[l])_k\mid \{0,\dots,j-1,j+1,\dots,l\}
\nsubseteq \{ f(0),\dots, f(k)\} \}.
\end{split}
\end{equation*}
In fact, the horn $\Lambda[l,j]$ is a simplicial set obtained from the
simplicial $l$-simplex $\Delta[l]$ by taking away its unique
non-degenerate $l$-simplex as well as the $j$-th of its $l+1$
non-degenerate $(l-1)$-simplices.

\begin{defi} \label{def:lie-n-gpd}
 A \emph{Lie $n$-groupoid}, see \cite{get:lie, hen:int, z:tgpd-2},  is a simplicial manifold $X_\bullet$ where the natural projections
\begin{equation*}\label{eq:horn-proj}
 p_{l,j}:X_l=\Hom(\Delta[l],X_\bullet)\to \Hom(\Lambda[l,j], X_\bullet)
\end{equation*}
 are surjective submersions for all $1\le l\ge j\ge 0$ and diffeomorphisms for all $0\le j\le l>n$. If $X_0=pt$, we say that $X_\bullet$ is a \emph{Lie $n$-group}.   
 
  A \emph{strict Lie $n$-group}, see \cite[$\S 6$]{wol:des},  is a simplicial Lie group $K_\bullet$ whose underlying simplicial  manifold is  a Lie $(n-1)$-groupoid. 
\end{defi}

\begin{rema}\label{Rmk:CrossMods1}
    Historically, strict Lie $2$-groups have been much studied \cite{Baez_Lauda:2004, hen:int, mac:cat}. According to Definition \ref{def:lie-n-gpd}, strict Lie 2-groups are groupoid objects in the category of Lie groups. In \cite{Baez_Lauda:2004}, other characterizations of these are extensively studied, chief among which are \emph{crossed module of Lie groups}, i.e. tuples $(G,H,\partial,\alpha)$ where $G$ and $H$ are Lie groups, $\partial:H\to G$ is a Lie group homomorphism, and $\alpha:G\times H\to H$, $\alpha(g,h)=h^g$ is an action by automorphisms satisfying
$$     \partial(h^g)=g\cdot\partial(h)\cdot g^{-1}\quad \text{and}\quad 
h_{2}^{\partial(h_{1})}=h_{1}\cdot h_{2}\cdot h_{1}^{-1}\qquad\text{for } g\in G, h, h_1,h_2\in H.$$    
 A crossed module of Lie groups gives rise to a strict Lie 2-group $H\rtimes G\rightrightarrows G$ whose group of arrows is the semi-direct product with respect to $\alpha$, and whose structural maps as a groupoid are:
 $$d_0(h,g)=g,\quad d_1(h,g)=\partial(h)\cdot g,\quad (h_2, \partial(h_1)\cdot g)\Join (h_1, g)=(h_2h_1,g).$$
 Conversely, the restriction of the target to the kernel of the source map together with conjugation by units in the group of arrows yield a crossed module, see \cite[Prop. 32]{Baez_Lauda:2004} or \cite[\S XII.8]{mac:cat}.   
\end{rema}

\begin{defi}
A \emph{Lie $n$-algebra} $(\cV, \ell_i)$, see e.g. \cite{bc:lie2a}, is a graded vector space $\cV$ concentrated in degrees $-n+1$ to $0$ together with graded antisymmetric linear maps $\ell_i:\cV^{\otimes i}\to \cV[2-i]$ for $i\in\{1,\cdots, n+1\}$ satisfying
\begin{equation}\label{eq:ln}
    \sum_{i+j=r+1}\sum_{\sigma\in Sh(i,j-1)}(-1)^{i(j-1)}K(\sigma)\ell_j\big(\ell_i(v_{\sigma(1)}, \cdots,v_{\sigma(i)}),v_{\sigma(i+1)}, \cdots, v_{\sigma(r)}\big)=0
\end{equation}
where $Sh(i,j-1)$ denotes the $(i,j-1)$-shuffles and $K(\sigma)$ is the graded sign of $\sigma$.

A \emph{strict Lie $n$-algebra} is a simplicial Lie algebra $\fk_\bullet$ whose underlying simplicial vector space is a Lie $(n-1)$-groupoid.     
\end{defi}

\begin{rema}[Lie $n$-algebras as $NQ$-manifolds]\label{rem:CE} 
Let $(\cV, \ell_i)$ be a Lie $n$-algebra. Consider the graded algebra $\sym (\cV[1])^*$ where $\sym$ here denotes the graded symmetric algebra of the positively graded vector space $(\cV[1])^*$. A direct consequence of Eq. \eqref{eq:ln} is that dualizing the structural maps $\ell_i$, one can build a Chevalley-Eilenberg type differential on this algebra that we denote by $\delta_{CE}$.  Using the language of supergeometry, we say that the above is the algebra of functions on the \emph{degree $n$ manifold} $\underline{\cV[1]}$, i.e. $$C(\underline{\cV[1]})=\sym (\cV[1])^*,$$ and $\delta_{CE}$ becomes a degree $1$ vector field on  $\underline{\cV[1]}$ satisfying $[\delta_{CE},\delta_{CE}]=0$. This structure is known as a \emph{degree $n$ $NQ$-manifold} (see \cite{bon:on, mnev:book} for more details).
\end{rema}

\begin{rema}
    One can introduce Lie \emph{$n$-algebroids} in the same way, as a graded vector bundle $\cV\to M$ with an anchor $\rho:\cV_0\to TM$ and brackets $\ell_i$ that are all $C^\infty(M)$-linear except for $\ell_2$, which satisfies the Leibniz rule. In \cite{bon:on} it is shown how to obtain a corresponding degree $n$ $NQ$-manifold.
\end{rema}

\subsection{The functors $\dec$ and $\overline{W}$}\label{sec:dec}

Following \cite{ill:com}, we define the \emph{d\'{e}calage} of a simplicial object $X_\bullet$ on a category $\cC$ as the simplicial object given by
$$\dec_n(X)=X_{n+1}, \quad d_k^{\dec}=d_{k+1}\quad\text{and}\quad s^{\dec}_k=s_{k+1},\quad\text{for } 0\leq k\leq n.$$
In this article, the following properties of the d\'{e}calage construction are relevant for us. 

\begin{rema}\label{prop:dec}
    The d\'{e}calage construction is functorial and satisfies:
    \begin{enumerate}
        \item $d_0:\dec_\bullet(X)\to X_\bullet$ is a simplicial map, indeed is a fibrant replacement\footnote{Here we use the iCFO structure for Lie $n$-groupoids constructed in \cite{cc:icfo}.} for $s_0:X_0\to X_\bullet$, where $X_0$ is the constant simplicial manifold given by the total units of $X_\bullet$.
        \item If $K_\bullet$ is a (strict) Lie $n$-group then so is $\dec_\bullet(K)$.
    \end{enumerate}
\end{rema}

Following \cite{eil:on}, we define the \emph{$\overline{W}$ functor} from simplicial groups to simplicial sets as follows: Let $K_\bullet$ be a simplicial group, then 
$$\overline{W}_0 K=pt, \quad \overline{W}_{n+1} K=K_{n}\times K_{n-1}\times\cdots\times K_1\times K_0$$
and for $\vec{\gamma}=(\gamma^{n};...;\gamma^0)\in K_{n}\times\cdots\times K_0=\overline{W}_{n+1} K$ the face and degeneracy maps are
\begin{align}\label{Eq:BFaceMaps}
\bar{\partial}_i(\vec{\gamma})&=\begin{cases}
(\gamma^{n-1};...;\gamma^0)& \textnormal{if }k=0 \\
(d_{i-1}\gamma^{n};...;d_{1}\gamma^{n-i+2};(d_0\gamma^{n-i+1})\gamma^{n-i};\gamma^{n-i-1};...;\gamma^0)& \textnormal{if }0<i\leq n \\
(d_{n}\gamma^{n};...;d_1\gamma^{1})& \textnormal{if }i=n+1,\end{cases}\\
\bar{\sigma}_i(\vec{\gamma})&=(s_{i-1}\gamma^{n};...;s_0\gamma^{n+1-i};\mathbbm{1}_{n+1-i};\gamma^{n-i};...;\gamma^0).
\end{align}

\begin{nota}
   Combining the above constructions, for a simplicial group $K_\bullet$, we define $$W_\bullet K:=\dec_\bullet(\overline{W}K),$$ whose face and degeneracy maps we denote by $\partial$ and $\sigma$ to avoid confusion.  
\end{nota}
 We state the main properties of these constructions.

\begin{prop}\label{prop:W}
    Let $K_\bullet$ be a simplicial group. The following statements hold:
    \begin{enumerate}
        \item If $K_\bullet$ is a strict Lie $n$-group then $\overline{W}_\bullet K$ is a Lie $n$-group. 
        \item The simplicial group $K_\bullet$ acts on the left on $W_\bullet K$, i.e. we have a simplicial morphism  $$L:K_\bullet\times W_\bullet K\to W_\bullet K,\quad  L_{k^n}(\gamma^n;\dots; \gamma^0)=(k^n\gamma^n;\dots; \gamma^0).$$
        \item $W_\bullet K$ is Morita equivalent to a point. 
    \end{enumerate}
\end{prop}
\begin{proof}
    The first statement follows from \cite[Thm. 6.7]{wol:des}, the second statement is classical and easy to verify, while the last one follows directly from the first item in Remark \ref{prop:dec}. 
\end{proof}

Proposition \ref{prop:W} collects the facts that for a strict Lie $n$-group $K_\bullet$, $\overline{W}_\bullet K$ is a simplicial model for its classifying space, while $W_\bullet K$ is a model for its universal bundle.

\subsection{Strict Lie $2$-groups and their strict Lie $2$-algebras}

The main characters of this article are strict Lie $2$-groups that we denote by $G_\bullet$ with faces given by $d_i$ and degeneracies by $s_j$. As pointed out in \cite{wol:des}, a strict Lie $1$-group is a Lie group $K$ viewed as a constant simplicial group, and $\overline{W}K$ corresponds to its nerve. Similar to this case, a strict Lie $2$-group is a simplicial Lie group $G_\bullet$ whose underlying simplicial manifold is the nerve of the Lie groupoid $G_1\rightrightarrows G_0$. 

\begin{nota}\label{Notation}
    To make clear the distinction between the group operation and the groupoid operation in a strict Lie 2-group $G_\bullet$, 
we adopt the following convention: 
\begin{eqnarray*}
\gamma_1\vJoin\gamma_2, & \gamma_3\Join\gamma_4
\end{eqnarray*}
stand respectively for the group multiplication and the groupoid 
multiplication whenever $(\gamma_1,\gamma_2)\in G_1\times G_1$ and 
$(\gamma_3,\gamma_4)\in G_1\times_{G_0}G_1=G_2$. The notation for the groupoid multiplication intends to reflect the \emph{joining} of consecutive arrows. Also, for $\gamma\in G_1$, $\gamma^{-1}$ and $\overline{\gamma}$ stand respectively for the inverses of the group and the groupoid multiplications.
\end{nota}


Since $G_n$ are Lie groups for each $n$, one can use the usual Lie functor to obtain the following result.  

\begin{prop}\label{prop:lie}
    Let $G_\bullet$ be a strict Lie $2$-group and denote by $\lfun$ the Lie functor from the category of Lie groups to the category of Lie algebras. Then
    \begin{enumerate}
        \item The simplicial manifold $\g_\bullet$ given by $$\g_p=\lfun(G_p),\quad \hat{\partial_i}=\lfun(d_i)\quad\text{and}\quad \hat{\sigma_i}=\lfun(s_i)$$ is a strict Lie $2$-algebra. 
        \item The graded vector space 
    $$\cV_0=\g_0,\quad \cV_{-1}=\h=\ker\hat{\partial}_0 \quad\text{with}\quad \ell_1=\hat{\partial}_1\quad\text{and}\quad \ell_2(\cdot,\cdot)=[\cdot,\cdot]$$
    is a Lie $2$-algebra.
        \item The strict Lie $2$-algebra $\g_\bullet$ acts on $W_\bullet G$.
         For $x\in\gg_0$, we abuse notation and write $x:=\hat{\sigma}_0^p(x)\in\gg_p$, for $y\in\hh$, we write $y_\beta\in\gg_p$ for 
\begin{eqnarray*}
y_\beta=(\underbrace{\hat{\sigma}_0\circ\hat{\partial}_1(y),...,\hat{\sigma}_0\circ\hat{\partial}_1(y)}_\text{$\beta$-1}, y,\underbrace{0,...,0}_\text{p-$\beta$}\big{)}\in\gg_p\subset\gg_1^p.
\end{eqnarray*}
        The infinitesimal generators corresponding to $x,y_\beta\in \g_p$ are denoted respectively by 
\begin{equation*}
    x^{p+1}\in \fX(W_pG)=\fX(\overline{W}_{p+1}G) \qquad\text{and}\qquad y_\beta^{p+1}\in\fX(W_pG)=\fX(\overline{W}_{p+1}G).
\end{equation*}
    \end{enumerate}
\end{prop}
\begin{proof}
    The first item is a clear consequence of the functoriality of $\lfun$. The second follows directly from the main result in \cite{bc:lie2a}. The third is immediate from Proposition \ref{prop:W}.
\end{proof}

\begin{rema}\label{Rmk:CrossMods2}
    Continuing Remark~\ref{Rmk:CrossMods1}, recall that a \emph{crossed modules of Lie algebras} $(\g,\h,\partial,\rho)$ consists of a pair of Lie algebras $\g$, $\h$ together with two Lie algebra homomorphisms $\partial:\hh\to \gg$ and $\rho:\g\to \text{Der}(\h)$ satisfying $$\partial(\rho(x)(y))=[x,\partial(y)]\quad\text{and}\quad \rho(\partial y_1)(y_2)=[y_1,y_2],\qquad\text{for } x\in\gg, y, y_1,y_2\in\hh.$$
    Crossed modules of Lie algebras can also be shown to be equivalent to Lie $2$-algebras with $\ell_3\equiv 0$; indeed, $\ell_1:\cV_{-1}\to\cV_0$ can be made into a Lie algebra homomorphism endowing $\cV_{-1}$ with the Lie bracket $[u,v]:=\ell_2(\ell_1(u),v)$, and conversely, one only needs to forget the bracket on the domain of the crossed module homomorphism, see \cite[Prop. 48]{bc:lie2a}. One can easily differentiate a strict Lie 2-group to a strict Lie 2-algebra either by differentiating the groupoid object in the category of Lie groups to a groupoid object in the category of Lie algebras as in Proposition \ref{prop:lie}, or by looking at the crossed modules \cite{bc:lie2a}. In Section~\ref{sec:dif}, we propose a different method based on the simplicial nature of strict Lie 2-groups.
\end{rema}

\begin{rema}\label{Rmk:2CE}[The Chevalley-Eilenberg complex]
As a consequence of Proposition~\ref{prop:lie}, we can define the Chevalley-Eilenberg complex of a strict Lie $2$-algebra $\g_\bullet$ as the Chevalley-Eilenberg complex of the Lie $2$-algebra $(\cV, \ell_i).$ More concretely, 
$$C^i_{CE}(\h\to\g_0)=\bigoplus_{k+2l=i}\bigwedge\nolimits^ k \g_0^*\otimes \sym^l\h^*$$
and the Chevalley-Eilenberg differential $\delta_{CE}$ is given by the two components
\begin{align*}
    (\delta_{CE}\eta)(x_0,\cdots, x_{k}; y_0,\cdots, y_{l-1})& = \sum_{i,j}(-1)^{i+1} \eta(x_0,\cdots, \widehat{x}_i,\cdots, x_{k}; \ell_2(x_i,y_j), y_0,\cdots,  \widehat{y}_j,\cdots, y_{l-1})\\
    &\quad+\sum_{i<j}(-1)^{i+j} \eta(\ell_2(x_i,x_j), x_0,\cdots, \widehat{x}_i,\cdots, \widehat{x}_j,\cdots, x_{k}; y_0,\cdots, y_{l-1}) \\
    (\delta_{CE}\eta)(x_0,\cdots, x_{k-2}; y_0,\cdots, y_{l})& = (-1)^{l+1}\sum_{i=0}^{l} \eta(x_0,\cdots, x_{k-2}, \ell_1(y_i); y_0,\cdots, \widehat{y}_i,\cdots, y_{l})
\end{align*}
for $\eta\in \bigwedge^k\g_0^*\otimes\sym^l\h^*$, $x_i\in\g$ and $y_j\in\h$. 
\end{rema}   

\begin{exa}\label{Rmk:CrossMods3}
As an example of the construction in Remark~\ref{Rmk:2CE}, recall that the standard Weil algebra of a Lie algebra $\gg$ can be seen as the Chevalley-Eilenberg complex of the {\it tangent Lie 2-algebra}; namely, $\id:\gg\to\gg$ with $\ell_2$ given by the Lie bracket. Using the correspondence of Remark~\ref{Rmk:CrossMods2}, this can also be seen as the crossed module of Lie algebras $(\gg, \g, \id, ad)$ that integrates to the crossed module of Lie groups $(G,G,\id,Ad)$ where $Ad$ denotes the action by conjugation.
\end{exa}

\section{A simplicial differentiation}\label{sec:dif}

In $\S\ref{sec:dec}$ we exhibit a model for the classifying space of a strict Lie $n$-group. Let $G_\bullet$ be a strict Lie $2$-group. In the sequel, we construct the Lie $2$-algebra of $\overline{W}_\bullet G$ using its universal bundle.  Due to the functoriality of $\dec$ and $\overline{W}$, we get the following diagram
\begin{equation}\label{dia:dec-w}
\begin{array}{c}
     \xymatrix{ W_\bullet(\dec G) \ar[d]_{\bar{\partial}_0^{\dec}}\ar[r]^{\quad Wd_0}& W_\bullet G\ar[d]^{\bar{\partial}_0}\\  \overline{W}_\bullet(\dec G)\ar[r]_{\quad\overline{W}d_0}&\overline{W}_\bullet G.}
\end{array}
\end{equation}
An immediate consequence of Proposition \ref{prop:W} is that $W_\bullet(\dec G)$ and $W_\bullet G$ are both Morita equivalent to a point. 

\begin{rema}
    Notice that for $G_\bullet$, one can also construct the classifying space by first taking the bi-simplicial manifold associated with the double groupoid $G_1\rightrightarrows G_0$ and then the Artin-Mazur codiagonal construction as explained in \cite{AM:1966,raj:from}. A simple computation shows that the Lie $2$-group thus obtained is isomorphic to $\overline{W}_\bullet G$ (see \cite{eil:on}).
\end{rema}

Even though $d_0:\dec_\bullet(G)\to G_\bullet$ is a morphisms of simplicial manifolds, $s_0$ is not; thus, $Ws_0$ does not define a section for $Wd_0$. However, one can twist $s_0$ on higher simplicial degrees so as to get an honest section for $Wd_0$. That is the content of the following key Proposition.

\begin{prop}\label{prop:eps}
    The simplicial map $\varepsilon_\bullet:W_\bullet G\to W_\bullet(\dec G)$ given by the recursive formula
    $$\varepsilon_0(g)=s_0(g),\quad \varepsilon_1(\gamma;g)=\big(s_0\gamma;\gamma^{-1}\vJoin s_0d_0(\gamma)\vJoin s_0(g)\big)\quad\text{and}\quad$$ $$\varepsilon_p=\Big(s_0\circ \pr_{p};\big(\bar{\partial}_0^{\dec}\circ\varepsilon_1\circ\partial_2\circ\cdots\circ\partial_p, \pr_{p-1}\big);\bar{\partial}^{\dec}_0\circ\varepsilon_{p-1}\circ\partial_0\Big)$$
   where $\pr_{j}:W_p(G)\to G_j$ denotes the natural projection, is a simplicial splitting of the map $W_\bullet d_0: W_\bullet(\dec G)\to W_\bullet G$.
\end{prop}
\begin{proof}
     First, we show that $\varepsilon_p$ is well-defined. Given $\vec{\gamma}=(\gamma^p;...;\gamma^0)\in W_pG$, then 
     \begin{align*}
        d_0\circ\bar{\partial}_0^{\dec}\circ\varepsilon_1\circ\partial_2\circ ...\circ\partial_{p}(\vec{\gamma}) & = d_0\circ\bar{\partial}_0^{\dec}\circ\varepsilon_1(d_2\circ ...\circ d_{p}(\gamma^p);d_1\circ\dots\circ d_{p-1}(\gamma^{p-1}))\\
        &=d_1\circ\dots\circ d_{p-1}(\gamma^{p-1});
    \end{align*}
    therefore, $\big(\bar{\partial}_0^{\dec}\circ\varepsilon_1\circ\partial_2\circ\cdots\circ\partial_p, \pr_{p-1}\big)\in G_p$.

   Next, we show that $\varepsilon_\bullet$ is compatible with the face maps. Let $(\gamma;g)\in W_1G$, then
    \begin{align*}
        \partial_0^{\dec}(\varepsilon_1(\gamma;g))  =&\  \partial_0^{\dec}\big(s_0(\gamma);\gamma^{-1}\vJoin s_0d_0(\gamma)\vJoin s_0(g)\big)
        =\ (s_0d_1(\gamma)\Join\gamma)\vJoin\gamma^{-1}\vJoin s_0(d_0(\gamma)\vJoin g)\\
        =&\ s_0(d_0(\gamma)\vJoin g)= \varepsilon_0(\partial_0(\gamma;g));\\
        \partial_1^{\dec}(\varepsilon_1(\gamma;g)) =&\  \partial_1^{\dec}\big(s_0(\gamma);\gamma^{-1}\vJoin s_0(d_0(\gamma)\vJoin g)\big)=d_1^{\dec}s_0(\gamma)=s_0{d_1(\gamma)}= \varepsilon_0(\partial_1(\gamma;g)).
    \end{align*}
    Inductively, assume $\partial_k^{\dec}\circ\varepsilon_r=\varepsilon_{r-1}\circ\partial_k$ for all $k$ and all $r\leq p-1$. For $k\geq 2$ it follows from the simplicial identities that 
    \begin{align*}
        \partial_k^{\dec}\circ\varepsilon_p & =\Big{(}d_k^{\dec}\circ s_0\circ \pr_n;d_{k-1}^{\dec}\circ\big(\bar{\partial}_0^{\dec}\circ\varepsilon_1\circ\cdots\circ\partial_n, \pr_{n-1}\big);\partial_{k-2}^{\dec}\circ\bar{\partial}^{\dec}_0\circ\varepsilon_{p-1}\circ\partial_0\Big{)} \\
            & =\Big{(}s_0\circ d_{k}\circ\pr_n;\big(\bar{\partial}_0^{\dec}\circ\varepsilon_1\circ\cdots\circ\partial_n, d_{k-1}\circ \pr_{n-1}\big);\bar{\partial}^{\dec}_0\circ\varepsilon_{p-1}\circ\partial_0\circ\partial_{k}\Big{)} = \varepsilon_{p-1}\circ\partial_k.
    \end{align*}
      For $k=0$, let $\vec{\gamma}=(\gamma^p;\cdots; \gamma^0)\in W_pG$, where each $\gamma^j=(\gamma^j_1,\cdots\gamma^j_j)\in G_j$, then 
    \begin{align*}
        \partial^{\dec}_0\circ\varepsilon_p(\vec{\gamma}) & =\Big{(}d^{\dec}_0(s_0(\gamma^p))\vJoin\big{(}\bar{\partial}^{\dec}_0\circ\varepsilon_1(\gamma^p_1;d_1\gamma^{p-1}_1), \gamma^{p-1}\big{)};\bar{\partial}^{\dec}_0\circ\varepsilon_{p-1}\circ\partial_0(\vec{\gamma})\Big{)} \\
&= \Big{(}\gamma^p\vJoin\big((\gamma^p_1)^{-1}\vJoin s_0d_1(\gamma^{p-1}_1), \gamma^{p-1}\big{)};\bar{\partial}^{\dec}_0\circ\varepsilon_{p-1}\circ\partial_0(\vec{\gamma})\Big{)}\\
        &= \Big{(}s_0\big(d_0(\gamma^p)\vJoin\gamma^{p-1}\big);\big(\bar{\partial}^{\dec}_0\circ\varepsilon_1(\gamma^p_2\vJoin\gamma^{p-1}_1;d_1\gamma^{p-2}_1),\gamma^{p-2}\big);\bar{\partial}^{\dec}_0\circ\varepsilon_{p-2}\circ\partial_0\circ\partial_0(\vec{\gamma})\Big{)}\\
            & =\varepsilon_{p-1}(d_0(\gamma^p)\vJoin\gamma^{p-1};\gamma^{p-2};\cdots;\gamma^0)=\varepsilon_{p-1}\circ\partial_0(\vec{\gamma}).        
    \end{align*}
    Lastly, when $k=1$, on the one hand, we get
    \begin{align*}
        \partial^{\dec}_1\circ\varepsilon_p(\vec{\gamma}) & =\big{(}d_1^{\dec}\circ s_0(\gamma^p);(X, d_0(\gamma^{p-1})\vJoin\gamma^{p-2}); \partial^{\dec}_1\circ\bar{\partial}^{\dec}_0\circ\varepsilon_{p-1}\circ\partial_0(\vec{\gamma})\big{)} \\
            & =\Big{(}s_0\circ d_1(\gamma^p);(X,d_0(\gamma^{p-1})\vJoin\gamma^{p-2}) ;\bar{\partial}^{\dec}_0\circ\varepsilon_{p-2}\circ\partial_0\circ\partial_1(\vec{\gamma})\Big{)}, 
    \end{align*}
where 
\begin{align*}
    X=&d_{0}^{\dec}\bar{\partial}^{\dec}_0\circ\varepsilon_1(\gamma^p_1;d_1\gamma^{p-1}_1)\vJoin\bar{\partial}^{\dec}_0\circ\varepsilon_1(\gamma^p_2\vJoin\gamma^{p-1}_1;d_1\gamma^{p-2}_1)\\
    =& \Big{(}\big((\gamma^p_1)^{-1}\vJoin s_0(d_0\gamma^p_1\vJoin d_1\gamma^{p-1}_1)\big)\Join\gamma^{p-1}_1\Big{)}\vJoin(\gamma^p_2\vJoin\gamma^{p-1}_1)^{-1}\vJoin s_0\big(d_0(\gamma^p_2\vJoin\gamma^{p-1}_1)\vJoin d_1(\gamma^{p-2}_1)\big) .
\end{align*}    
 On the other hand, 
    \begin{align*}
       \varepsilon_p\circ\partial_1(\vec{\gamma}) & =\Big{(}s_0\circ d_1(\gamma^p);\big(X',d_0(\gamma^{p-1})\vJoin\gamma^{p-2}\big);\bar{\partial}^{\dec}_0\circ\varepsilon_{p-2}\circ\partial_0\circ\partial_1(\vec{\gamma})\Big{)}, 
    \end{align*}
    where $X'=(\gamma^p_1\Join\gamma^p_2)^{-1}\vJoin s_0\big(d_0(\gamma^p_2)\vJoin d_1(\gamma^{p-1}_2\vJoin\gamma^{p-2}_1)\big)$; thus, the result follows if $X=X'$, which, due to 
    \begin{align*}
        \big(\gamma^{-1}\vJoin s_0(d_0\gamma\vJoin d_1\gamma')\big)\Join\gamma' & =\big{(}(\gamma^{-1}\vJoin s_0d_0\gamma)\vJoin s_0d_1\gamma')\big{)}\Join\big{(}1\vJoin\gamma'\big{)} \\
            & =\big{(}(\gamma^{-1}\vJoin s_0d_0\gamma)\Join s_0(1)\big{)}\vJoin\big{(}s_0d_1\gamma'\Join\gamma'\big{)}=\gamma^{-1}\vJoin (s_0d_0\gamma)\vJoin\Join\gamma'
    \end{align*}
    holding for all $\gamma,\gamma'\in G_1$, ultimately boils down to prove 
    $$(\gamma^p_1)^{-1}\vJoin s_0d_0(\gamma_p^1)\vJoin(\gamma^p_2)^{-1}=(\gamma^p_1\Join\gamma^p_2)^{-1},$$
    which follows from an easy computation. 
    
    Finally, we show that $\varepsilon_\bullet$ is compatible with the degeneracy maps. For the first degree, we get
     \begin{align*}
        \varepsilon_1(\sigma_0(g))=\varepsilon_1(s_0g;1)=(s_0s_0g;(s_0g)^{-1}\vJoin s_0(d_0s_0g\vJoin 1))=(s_0s_0g;1)=\sigma_0^{\dec}(s_0g)=\sigma_0^{\dec}(\varepsilon_0(g)).
    \end{align*}
    Inductively, assume $\varepsilon_r\circ\sigma_k=\sigma_k\circ\varepsilon_{r-1}$ for all $k$ and all $r\leq p$ then
    \begin{align*}
        \varepsilon_{p+1}(\sigma_0\vec{\gamma}) & =\Big{(}s_0s_{0}(\gamma^p);(\bar{\partial}^{\dec}_0\circ\varepsilon_1\circ\partial_2\circ\cdots\partial_{p+1}\circ\sigma_0, \mathbbm{1}_{p});(\bar{\partial}^{\dec}_0\circ\varepsilon_{p}\circ\partial_0\circ\sigma_0)(\vec{\gamma})\Big{)} \\
            & =\Big{(}s_1s_{0}(\gamma^p);\mathbbm{1}_{p+1};(\bar{\partial}^{\dec}_0\circ\varepsilon_{p})(\vec{\gamma})\Big{)}=\sigma_0^{\dec}\circ\varepsilon_p(\vec{\gamma}).
    \end{align*}
   For $k>0$, it follows from the simplicial identities and the induction hypothesis that 
    \begin{align*}
        \varepsilon_{p+1}\circ\sigma_k & =\Big{(}s_{0}s_{k}\circ\pr_p;\big(\bar{\partial}^{\dec}_0\circ\varepsilon_1\circ\partial_2\circ\cdots\partial_{p}, s_{k-1}\circ \pr_{p-1}\big);\bar{\partial}^{\dec}_0\circ\varepsilon_{p}\circ\partial_0\circ\sigma_k\Big{)} \\
            & =\Big{(}s_{k+1}s_{0}\circ\pr_p;s_{k}\circ\big(\bar{\partial}^{\dec}_0\circ\varepsilon_1\circ\cdots\circ\partial_{p}, \pr_{p-1}\big);\bar{\partial}^{\dec}_0\circ\sigma^{\dec}_{k-1}\circ\varepsilon_{p-1}\circ\partial_0\Big{)} =\sigma^{\dec}_k\circ\varepsilon_{p}.
    \end{align*}
    Notice that the last equality in the above formula for $p=1$ follows from
    \begin{eqnarray*}
        \bar{\partial}^{\dec}_0\circ\varepsilon_{1}\circ\partial_0\circ\sigma_1(\gamma;g) & =s_0(d_0(\gamma)\vJoin g)^{-1}\vJoin s_0d_0s_0(d_0(\gamma)\vJoin g)=1,
    \end{eqnarray*}
    while for higher $p$ the equality is immediate from the definition of $\varepsilon$.
    Lastly, note that a simple inductive argument implies that $\varepsilon$ is indeed a splitting of $Wd_0$, as $Wd_0 \circ \varepsilon_1(\gamma, g)=(\gamma, g)$.
\end{proof}

\begin{rema}\label{Rmk:CrossMods4}
    In view of the correspondence of Remark~\ref{Rmk:CrossMods1}, the reader might be wondering what the map $\varepsilon$ of Proposition~\ref{prop:eps} looks like when expressed in terms of crossed modules. Since letting $H:=\ker d_0$, $G_1\cong H\rtimes G_0$ canonically, then $G_p\cong H^p\times G_0$. Consequently, $W_pG\cong H^{\sum_{\alpha=1}^{p}\alpha}\times G_0^{p+1}$ and $W_p(\dec G)\cong H^{\sum_{\alpha=0}^{p}(\alpha+1)}\times G_0^{p+1}=H^{p+1}\times W_pG$. In this presentation, $\varepsilon_p:W_p G\longrightarrow W_p(\dec G)=H^{p+1}\times W_pG$ assumes the shape $\varepsilon_p=\mu_p\times\id_{W_pG}$, where $\mu_p:W_pG\longrightarrow H^{p+1}$ has coordinates $(1,\mu_p^1,\mu_p^2...,\mu_p^{p-1},\mu_p^p)$. Here, the first entry is the constant map to the identity element of $H$, and $\mu_p^j$ is recursively given by
    \begin{align*}
    \mu_1^1(h,g_1;g_0) & =(h^{-1})^{g_0}, & 
    \mu_p^j & =\mu_1^1\circ\bigcirc_{a=0}^{j-2}\partial_a\circ\bigcirc_{b=j-1}^{p-2}\partial_{b+2}, 
    \end{align*}
    where $\bigcirc$ denotes composition between the corresponding indices. 
    Although this presentation for the map $\varepsilon$ might seem cleaner, it hides the increasing complexity of the face maps when written in terms of crossed modules. We ultimately opted out of taking this angle because, in addition, the actions of Propositions \ref{prop:W} and \ref{prop:lie} are not only conceptually clearer when avoiding crossed modules, but also the formulas for these actions in terms of crossed modules have increasingly cumbersome formulas. The formula for $\mu_p^j$ ressembles Eq.~(80) in \cite{Angulo:2024}; however, as it was the case in the referred paper, there is no clear theoretical explanation for these maps. 
\end{rema}

\subsection{The simplicial strict Lie $2$-algebroid}

To differentiate $G_\bullet$, we introduce the simplicial vector bundles
$$\cF'_\bullet:=\ker T\bar{\partial}_0^{\dec}\subset T(W_\bullet(\dec G))\quad\text{and}\quad \cF_\bullet:=\ker T\bar{\partial}_0\subset T(W_\bullet G)$$
over $W_\bullet(\dec G)$ and $W_\bullet G$, respectively (see e.g. \cite{mat:vb} for a clear exposition on simplicial vector bundles).
 
Since $\dec$ and $\overline{W}$ are functorial, diagram \eqref{dia:dec-w} implies that, defining $\cK:=\ker T(Wd_0)|_{\cF'}$, we get a short exact sequence of simplicial vector bundles
\begin{eqnarray}\label{Eq:dec2}
\begin{array}{c}
     \xymatrix{
    0 \ar[r] & \mathcal{K}_\bullet \ar[d]\ar[r] & \F'_\bullet \ar[d]\ar[r] & \F_\bullet \ar[d]\ar[r] & 0 \\
     & W_\bullet(\dec G) \ar@{=}[r] & W_\bullet(\dec G) \ar[r]_{\quad Wd_0} & W_\bullet G & }
\end{array}
\end{eqnarray}
over different bases. Using the map $\varepsilon$ of Proposition \ref{prop:eps}, we get a new short exact sequence of simplicial vector bundles
\begin{eqnarray}\label{Eq:dec3}
\begin{array}{c}
     \xymatrix{
    0 \ar[r] & \cH_\bullet:=\varepsilon^*_\bullet\mathcal{K}_\bullet \ar[d]\ar[r] & \cA_\bullet:=\varepsilon^*_\bullet\F'_\bullet \ar[d]\ar[r]^{\qquad\rho} & \F_\bullet \ar[d]\ar[r] & 0 \\
     & W_\bullet G  \ar@{=}[r] & W_\bullet G  \ar@{=}[r] & W_\bullet G. & }
\end{array}
\end{eqnarray}

\begin{prop}\label{prop:A}
The following holds:
    \begin{enumerate}
         \item The infinitesimal actions described in Proposition \ref{prop:lie} give trivializations $\cA_p\cong W_pG\times \g_{p+1}$ and $\cH_p\cong W_pG\times\h$.
        \item For a fixed $p$, diagram  \eqref{Eq:dec3} is a strict Lie $2$-algebroid. 
        \item 
        $\cH_\bullet\oplus\cA_\bullet$ is a simplicial strict Lie $2$-algebroid with face and degeneracy maps induced by $\varepsilon_\bullet$. More concretely, the face maps  $\partial^\cA_i:\cA_p\to \cA_{p-1}$ are  given by 
        $$\partial^\cA_0(X,\Vec{\gamma})=(Td_1X\star 0_{\pr_p\varepsilon_p(\vec{\gamma})}, \partial_0\vec{\gamma})\quad \text{and}\quad \partial^\cA_i(X,\Vec{\gamma})=(Td_{i+1}X,\partial_i\vec{\gamma})$$
         where $X\in T_{\pr_{p+1}\varepsilon_p(\vec{\gamma})}G_{p+1}$, and  $\partial^\cH_i:\cH_p\to \cH_{p-1}$ are given by 
        $$\partial^\cH_0(Y,\Vec{\gamma})=(Td_1Y\star 0_{\pr_p\varepsilon_p(\vec{\gamma})}, \partial_0\vec{\gamma})\quad \text{and}\quad \partial^\cH_i(Y,\Vec{\gamma})=(Td_{i+1}Y,\partial_i\vec{\gamma})$$
         where $Y\in \ker T_{\pr_{p+1}\varepsilon_p(\vec{\gamma})}d_0\subset T_{\pr_{p+1}\varepsilon_p(\vec{\gamma})}G_{p+1}$ and $\star$ denotes the tangent of the multiplication $\vJoin$ in $G_p.$
    \end{enumerate}
\end{prop}

\begin{proof}
For a fixed $p$, the action defined in Proposition \ref{prop:W} is given by left translation on the first factor. Therefore, the infinitesimal generators for the action of $\dec_p (G)$ on $W_p(\dec G)$ give a trivialization for $\cF'_p$. We define the right-invariant sections of $\cA_p$ by
$$\xi^{p+1}:=\varepsilon_p^*\widetilde{\xi}\in \Gamma(\cA_p)$$
where $\widetilde{\xi}\in\fX(W_p(\dec G))$ is the vector field corresponding to $\xi\in\g_{p+1}.$  In particular, if $y\in\hh$, $y_1\in\gg_{p+1}$ and, due to the exact sequence \eqref{Eq:dec2}, we get that $\widetilde{y}_1\in\Gamma(\mathcal{K}_{p})$. These sections trivialize $\cK_p$. We define the right-invariant sections of $\cH_p$ accordingly as $$y_1^{p+1}=\varepsilon_p^*\widetilde{y}_1\in\Gamma(\H_p).$$
For the second part, observe that, since $\varepsilon$ is a section for $Wd_0$, we get the identity
\begin{equation}\label{eq:anc}
    f\xi^{p+1}=\varepsilon_p^*((Wd_0)^*f)\widetilde{\xi}
\end{equation}
for each $f\in C^\infty(W_pG)$ and $\xi\in\g_{p+1}$. This allows us to define the brackets on $\Gamma(\cA_p)$ and $\Gamma(\cH_p)$ by the formulas
\begin{eqnarray*}
        \left[\xi^{p+1},\zeta^{p+1}\right]_{\A_p}:=\varepsilon_p^*\left[\widetilde{\xi},\widetilde{\zeta}\right]_{\F'_p}\qquad \text{and}\qquad\left[\xi^{p+1},y_1^{p+1}\right]:=\varepsilon_p^*\left[\widetilde{\xi},\widetilde{y}_1\right]_{\F'_p},
    \end{eqnarray*}
    for $\xi,\zeta\in\g_{p+1}$ and $y\in\h$. It follows directly from \eqref{eq:anc} that the anchor and the $\ell_1=j_p$ brackets are given by 
    $$\rho_{\A_p}(\xi^{p+1})=(\hat{\partial}_0\xi)^{p+1}\qquad\text{and}\qquad j_p(y_1^{p+1})=y_1^{p+1},$$
    where $\hat{\partial}_0$ is the first face map of the simplicial lie algebra $\g_\bullet$, and $(\hat{\partial}_0\xi)^{p+1}\in\fX(W_pG)$ is the infinitesimal generator corresponding to $\hat{\partial}_0\xi\in\g_p$.

    The third part follows from noticing that the Lie 2-algebroid structures of part two can be seen as a right-invariant extension of $\hh\rightarrow\gg_{p+1}, y\mapsto y_1$ for all $p\geq 0$ extended by Leibniz. Indeed, $\hh\rightarrow\gg_{p+1}$ is a strict Lie 2-algebra with cokernel $\gg_{p}$ and 2-bracket 
    $$\ell_2(\xi,y):=\ell_2(\hat{\partial}_{1}\circ...\circ\hat{\partial}_{p+1}(\xi),y),$$
    where the 2-bracket on the right-hand side is the 2-bracket on $\hh\rightarrow\gg_{0}$. Now, when restricting to units, the structural maps restrict to strict Lie $2$-algebra maps. Indeed, differentiating the d\'ecalage fibration $d_0:\dec_\bullet(G)\rightarrow G_\bullet$ yields a map of simplicial Lie algebras $\gg_{\bullet+1}\rightarrow\gg_{\bullet}$ whose kernel is precisely $\hh\rightarrow\gg_{\bullet+1}$. 
\end{proof}

\begin{rema}\label{Rmk:CE2}
    As in Remark \ref{rem:CE}, we can use the language of supergeometry to recast the second point in the above Proposition \ref{prop:A} as saying that, for each $p$, there is a degree $2$ $NQ$-manifold $(\underline{\cH_p[2]\oplus\cA_p[1]}, (\delta_{CE})_p)$. The third point thus says that these $NQ$-manifolds assemble into a simplicial object $(\underline{\cH_\bullet[2]\oplus\cA_\bullet[1]}, \delta_\bullet)$ (see e.g. \cite{raj:qgrou} for more examples of these). The complex we describe in the following paragraph should be understood as the double-complex of functions $(C^\bullet(\underline{\cH_\bullet[2]\oplus\cA_\bullet[1]}), \partial, \delta)$, where $\partial$ is the simplicial differential \eqref{eq:simfun}.
\end{rema}

Dualizing the structure of Proposition \ref{prop:A} yields a double-complex $(C^{\bullet,\bullet}, \partial,  \delta)$ where
\begin{eqnarray}\label{Eq:AuxDoubleCx}
    C^{p,q}:=C_{\textnormal{CE}}^q(\H_p\rightarrow\A_{p})=\bigoplus_{k+2l=q}\Gamma\left(\bigwedge\nolimits^k\A_{p}^*\otimes \sym^l\H_p^*\right),
\end{eqnarray}
 the rows are given by the simplicial diferential\footnote{One can think $\Gamma(\cA^*_p)$ and $\Gamma(\cH^*_p)$ as linear functions on $\cA_p$ and $\cH_p$ and extend $\partial$ linearly to $\bigwedge$ and $\sym$.} $\partial:C^{p,q}\to C^{p+1,q}$ and the columns are given by the Chevalley-Eilenberg differential of the strict Lie $2$-algebroid $\delta:C^{p,q}\to C^{p,q+1}$  $\cH_p\to\cA_p$ as defined in Remark \ref{Rmk:CE2}. Moreover, pullback along the d\'ecalage projection yields a map of complexes 
\begin{eqnarray}\label{Eq:partialZero}
    \xymatrix{\big(\widehat{C}^\infty(\overline{W}_\bullet G), \bar{\partial}\big) \ar[r]^{\bar{\partial}_0^*\qquad} & \big(C^\infty(W_\bullet G)=C^{\bullet,0},\partial\big).}
\end{eqnarray}
The following result says that one can recover the Lie 2-algebra of $G_\bullet$ by approximating the cohomology of $C^{p,q}$. 
\begin{thm}\label{Theo:SimplicialDiffn}
    Let $G_\bullet$ be a strict Lie $2$-group with strict Lie $2$-algebra $\g_\bullet$. Let $(E_r^{p,q},d_r)$ be the spectral sequence of the filtration by columns of the double-complex $(C^{p,q}, \partial,\delta)$ in \eqref{Eq:AuxDoubleCx}. 
    Then, $(E_1^{p,q},d_1)$ is isomorphic to the Chevalley-Eilenberg complex of $\g_\bullet$. 
\end{thm}

\begin{proof}
    It follows from Remark \ref{prop:dec} that $\dec(X_\bullet)$ and $X_0$ are Morita equivalent; indeed, the cohomology of the complex $(C^\infty(\dec_\bullet(X)),\partial)$ is concentrated in degree $0$ and $s_0^*$ gives an explicit homotopy. This holds in particular for $W_\bullet G$, whose homotopy is given by $\bar{\sigma}_0^*$.

    It follows from Lemma \ref{lemma:h0} that $\bar{\sigma}_0^*$ is covered by a homotopy given by dualizing
    \begin{align}\label{Eq:ContrHomotopiesF}
        \xymatrix{h^{\A}_0:\A_p \ar[r] & \A_{p+1},} &\quad  h^{\A}_0(X,\vec{\gamma}):=\big(\hat{\sigma}_0(X\star 0_{s_0\pr_p(\vec{\gamma})}^{-1}), \bar{\sigma}_0(\vec{\gamma})\big),\\ \label{Eq:ContrHomotopiesH}
        \xymatrix{h^{\H}_0:\H_p \ar[r] & \H_{p+1},} &\quad h^{\H}_0(Y, \vec{\gamma}):=\big(\hat{\sigma}_{1}(Y\star 0_{s_0\pr_p(\vec{\gamma})}^{-1}), \bar{\sigma}_0(\vec{\gamma})\big),
    \end{align}
where 
$\star$ stands for the differential of the multiplication $\vJoin$ on $G_{p+1}$ and $\hat{\sigma}_i$ are degeneracy maps of $\g_\bullet$. Thus, $E_1^{p,q}\equiv 0$ for all $p>0$.

 We claim that $E_1^{0,q}\cong C_{\textnormal{CE}}^q(\hh\rightarrow\gg_0)$ as cochain complexes. Consider the following maps:
    \begin{align*}
        \xymatrix{\iota_0^*:C^\infty(G_0) \ar[r] & \Rr,} & (\iota_0^*f):=f(1), \textnormal{ for }f\in C^\infty(G_0), \\
        \xymatrix{\iota_0^*:\Gamma\big(\bigwedge\nolimits^k\A^*_0\otimes\sym^l\H_0^*\big) \ar[r] & \bigwedge^k\gg_0^*\otimes \sym^l\hh^* ,} & (\iota_0^*\omega)(\mathbbm{x};\mathbbm{y}):=\omega_1(\hat{\sigma}_0(x_1),...,\hat{\sigma}_0(x_k);y_1,...,y_l), 
    \end{align*}
    where $\mathbbm{x}=(x_1,...,x_k)$ where $x_a\in\gg_0$, and $\mathbbm{y}=(y_1,...,y_l)$ where $y_b\in\hh$.
    \begin{align*}
        \xymatrix{\pi_0^*:\Rr \ar[r] & C(G_0),}(\pi_0^*\lambda)(g):=\lambda,\textnormal{ for }\lambda\in\Rr, &  \textnormal{ and }g\in G_0 \\
        \xymatrix{\pi_0^*:\bigwedge^k\gg_0^*\otimes \sym^l\hh^* \ar[r] & \Gamma\left(\bigwedge\nolimits^k\cA^*_0\otimes\sym^l\cH_0\right)
        ,} 
    \end{align*}
    \begin{align}\label{Eq:piZero}
         (\pi_0^*\omega)_g(\mathbb{X};\mathbb{Y}):=\omega(\hat{\partial}_1(X_1\star 0_{s_0g}^{-1}),...,\hat{\partial}_1(X_k\star 0_{s_0g}^{-1});Y_1\star 0_{s_0g}^{-1},...,Y_l\star 0_{s_0g}^{-1}),
    \end{align}
     where $\omega\in\bigwedge^k\gg_0^*\otimes \sym^l\hh^*$, $\mathbb{X}=(X_1,...,X_k)$ where $X_a\in(\A_0)_g$, and $\mathbb{Y}=(Y_1,...,Y_l)$ where $Y_b\in(\H_0)_g$. Here and throughout, $0_\gamma$ is the zero in the corresponding fibre over $\gamma\in G_1$, and $\cdot^{-1}$ stands for the inverse in $TG_1$ (consistently with our notation convention, see Notation~\ref{Notation}).  We observe that, $\iota_0^*\circ\pi_0^*=\textnormal{Id}$.
     Conversely, we show that the image of $\pi_0^*$ lies in $E_1^{0,q}$. Let $(X,(\gamma;g))\in\A_1$ for some $(\gamma;g)\in W_1G$, then, for $X=Tm(X_0, X_1)$,
     \begin{align*}
         \partial^\cA_0(X,(\gamma;g)) & =\big(T_{(s_0d_1(\gamma),\gamma)}m(X_0,X_1)\star 0_{\gamma^{-1}\vJoin s_0(d_0\gamma\vJoin g)}, d_0\gamma\vJoin g\big)\in\A_0.
     \end{align*}
     Since, $T_{(s_0d_1(\gamma),\gamma)}m(X_0,X_1)\star 0_{\gamma^{-1}\vJoin s_0(d_0\gamma\vJoin g)}\star 0_{s_0(d_0\gamma\vJoin g)}^{-1}$ equals
     \begin{align*}
         T_{(s_0d_1(\gamma),\gamma)}m(X_0,X_1)\star T_{(s_0d_1(\gamma)^{-1},\gamma^{-1})}m(0_{s_0d_1(\gamma)^{-1}},0_{\gamma^{-1}}))
             = T_{\mathbbm{1}_2}m(X_0\star 0_{s_0(d_1(\gamma)\vJoin g)}^{-1},X_1\star 0_{\gamma}^{-1}) ,
     \end{align*}
     and this computation also works for $X=Tm(Y,0)$, it follows easily that $\partial\pi_0^*\omega\equiv 0$. 
     
     If $\omega\in \Gamma\left(\bigwedge\nolimits^k\cA^*_0\otimes\sym^l\cH_0\right)$ is such that $\partial\omega=0$, then
    \begin{align*}
        (\pi_0^*\iota_0^*\omega)_g(\mathbb{X};\mathbb{Y}) & =(\iota^*_0\omega)(\hat{\partial}_1(X_1\star 0_{s_0g}^{-1}),...,\hat{\partial}_1(X_k\star 0_{s_0g}^{-1});Y_1\star 0_{s_0g}^{-1},...,Y_l\star 0_{s_0g}^{-1}) \notag\\
          & =\omega_1(\hat{\sigma}_0\circ\hat{\partial}_1(X_1\star 0_{s_0g}^{-1}),...,\hat{\sigma}_0\circ\hat{\partial}_1(X_k\star 0_{s_0g}^{-1});Y_1\star 0_{s_0g}^{-1},...,Y_l\star 0_{s_0g}^{-1})=\omega_g(\mathbb{X};\mathbb{Y}),
    \end{align*}
     where the last equation holds as
     $$0=(\partial\omega)_{(1;g)}(\sigma_0(\mathbb{X});\sigma_1(\mathbb{Y}))=\omega_g(\mathbb{X};\mathbb{Y})-\omega_1(\hat{\sigma}_0\circ\hat{\partial}_1(\mathbb{X}\star 0_{s_0g}^{-1});\mathbb{Y}\star 0_{s_0g}^{-1}).$$ Summing up, $\iota_0^*$ restricts to an isomorphism of vector spaces between $E_1^{0,q}$ and $C_{\textnormal{CE}}^q(\hh\rightarrow\gg_0)$. 

       Lastly, we check that the restriction of the vertical differential to $E_1^{0,q}$ coincides with the Chevalley-Eilenberg differential on $\hh\rightarrow\gg_0$. Since both $\iota_0^*$ and $\pi_0^*$ are easily proven to be homomorphisms for the corresponding algebra structures, and, since the differential is generated by the lowest degrees extending by Leibniz rule, it suffices to prove the following:
    \begin{itemize}
        \item For $\lambda\in\Rr$, $\iota_0^*\delta\pi_0^*\lambda\equiv 0$ . Indeed, $\iota_0^*f$ is constant for any $f\in C(G_0)$.
        \item For $\omega\in\gg_0^*$, $\iota_0^*\delta\pi_0^*\omega=\delta_{\textnormal{CE}}\omega\in\bigwedge^2\gg_0^*\oplus\hh^*$. 
        Let $y\in\hh$, recall that, using the notation of Proposition \ref{prop:lie}, 
        $$y^1={y^1_1}|_{s_0(G_0)}\in\Gamma(\H_0).$$ 
        Then, note that
        \begin{align}\label{Eq:iotaTrickH}
            (\iota_0^*\delta\pi_0^*\omega)(y) & =\big{(}(\delta\pi_0^*\omega)(y^1)\big{)}(1);
        \end{align}
        therefore,
        \begin{align*}
            (\iota_0^*\delta\pi_0^*\omega)(y) & =-(\pi_0^*\omega)_1(j_0(y^{1}))=-\omega(\hat{\partial}_1(y))=-\omega(\ell_1(y)).
        \end{align*}
        For any $x\in\gg_0$, we have that 
        $$x^1=\hat{\sigma}_0(x)^1\Big{\vert}_{u(G_0)}\in\Gamma(\A_0),$$
        as in Proposition \ref{prop:lie}. Letting $x_0,x_1\in\gg_0$, note that
        \begin{align}\label{Eq:iotaTrickF}
            (\iota_0^*\delta\pi_0^*\omega)(x_0,x_1) & =\big{(}(\delta\pi_0^*\omega)(x_0^1,x_1^1)\big{)}(1);
        \end{align}
        therefore, as $\rho(x_a^1)=x_a^1\in\Gamma(\F_0)$ for $a\in\lbrace0,1\rbrace$,  
        \begin{equation*}
           (\iota_0^*\delta\pi_0^*\omega)(x_0,x_1)=\Big{(}\Lie_{x_0^1}(\pi_0^*\omega)(x_1^{1})-\Lie_{x_1^1}(\pi_0^*\omega)(x_0^1) 
               -(\pi_0^*\omega)([x_0^{1},x_1^{1}])\Big{)}(1). 
        \end{equation*}
        Now, $(\pi_0^*\omega)_g(x_a^{1})=\omega(x_a)$ is constant; whereas, by definition, $[x_0^1,x_1^1]=\ell_2(x_0,x_1)^1$. So, ultimately, 
        $$(\iota_0^*\delta\pi_0^*\omega)(x_0,x_1)=-\omega(\ell_2(x_0,x_1)).$$
        \item For $\omega\in\hh^*$, $\iota_0^*\delta\pi_0^*\omega=\delta_{\textnormal{CE}}\omega\in \gg_0^*\otimes\hh^*$. Letting $x\in\gg_0$ and $y\in\hh$ and combining the observations of Eq.'s~\eqref{Eq:iotaTrickH} and \eqref{Eq:iotaTrickF},
        \begin{equation*}
           (\iota_0^*\delta\pi_0^*\omega)(x;y)=\Big{(}\Lie_{x^1}(\pi_0^*\omega)(y^{1})
               -(\pi_0^*\omega)([x^1,y^{1}])\Big{)}(1). 
        \end{equation*}
        Again, $(\pi_0^*\omega)_g(y^{1})=\omega(y)$ is constant, and $[x^1,y^{1}]=\ell_2(x,y)^{1}$, thus implying
        \begin{align}
            (\iota_0^*\delta\pi_0^*\omega)(x;y)=-\omega(\ell_2(x,y)).&\qedhere
        \end{align}
    \end{itemize}
\end{proof}

\section{The van Est map}\label{sec:VE}

In this section, we construct a van Est map for a strict Lie $2$-group $G_\bullet$. Our map is similar to that in \cite{abad:VE}; however, our construction is inspired by \cite{Meinrenken_Li-Bland:2015}.

\subsection{The Weil algebra of a Lie $2$-algebroid}\label{ssec-WeilAlg}
Let $(\cV_\bullet\to M,\rho, \ell_i)$ be a Lie $n$-algebroid and let $(\underline{\cV[1]}, \delta_{CE})$ be its corresponding degree $n$  NQ-manifold as per Remark \ref{rem:CE}. Following \cite{teo:mod, raj:qal}, we say that the \emph{Weil algebra} of the Lie $n$-algebroid is the double complex $$\left(\cW^{\bullet,\bullet}(\cV), \delta, d\right):=\left(\Omega^{\bullet,\bullet}(\underline{\cV[1]}), \cL_{\delta_{CE}}, d_{DR}\right)$$
where $\Omega^{\bullet,\bullet}(\underline{\cV[1]})$ denotes the algebra of differential forms on the degree $n$ manifold $\underline{\cV[1]}$, $\cL_{\delta_{CE}}$ is the Lie derivative of differential forms with respect to the degree $1$ vector field $\delta_{CE}$, and $d_{DR}$ denotes the de Rham differential on the degree $n$ manifold $\underline{\cV[1]}$ (see also \cite{mein:weil}).

More concretely, for the strict Lie $2$-algebroid $(\cH_p\xrightarrow{j_p} \cA_p, [\cdot,\cdot])$, we get that $\cW^{q,r}(\cH_p\to\cA_p)$ is
$$\bigoplus_{j=0}^r\bigoplus_{a+b=j}\bigoplus_{k+a+2(l+b)=q}
\Gamma\left(\bigwedge\nolimits^{r-j}T^*(W_pG)\otimes\bigwedge\nolimits^k \cA^*_p\otimes \sym^l\cH^*_p\otimes\sym^a \cA^*_p\otimes\bigwedge\nolimits^b \cH_p^*\right).$$
One can express $\delta$ in terms of symmetric powers of the coadjoint representation up to homotopy of the Lie  $2$-algebroid. Explicit formulas for the differentials are an extension of the ones computed in \cite[$\S 5.4$]{teo:mod}; nevertheless, following an argument close to that in \cite[Example 4.9]{Meinrenken_Li-Bland:2015}, given that the strict Lie $2$-algebroid $(j_p:\cH_p\to \cA_p, [\cdot,\cdot])$ comes with the canonical trivialization of Proposition \ref{prop:A}, these formulas simplify greatly.

\subsection{The map}
Let $G_\bullet$ be a strict Lie $2$-group with strict Lie $2$-algebra $\g_\bullet$. Denote by $(\overline{W}_\bullet G, \bar{\partial}_i, \bar{\sigma}_j)$ the Lie $2$-group of Proposition \ref{prop:W} and by $\h\to\g_0$ the Lie $2$-algebra detailed in Proposition \ref{prop:lie}. As explained in \cite{bott:bss}, the de Rham complexes of $\overline{W}_p G$ fit into a double-complex, that we refer as the \emph{Bott-Shulmann-Stasheff double-complex}, $$(\Omega^\bullet(\overline{W}_\bullet G), \bar{\partial}, d),$$
where $\bar{\partial}$ is the simplicial differential and $d$ is the de Rham differential. Moreover,  \emph{normalized forms} $$\widehat{\Omega}^{p,r}(G)=\{ \omega\in \Omega^r(\overline{W}_p G) \ | \ \bar{\sigma}_k^*\omega=0,\forall k\}$$
form a differential subalgebra with respect to the \textit{cup product} that the simplicial de Rham complex naturally carries \cite{bott:bss,Meinrenken_Li-Bland:2015}. On the other hand, as mentioned above, the Weil algebra of the strict Lie $2$-algebra
$$\left(\cW^{\bullet,\bullet}(\h\to \g_0), \delta, d\right),$$
also defines a double-complex. The van Est map relates these double complexes by means of the infinitesimal action described in Proposition \ref{prop:lie}. The van Est map is defined as a linear combination of permutations of certain operators we move on to define. For $x\in\gg_0$, recall that $x^p\in\fX(\overline{W}_p G)$ denotes the infinitesimal generator of the action of Proposition \ref{prop:lie} and define the operators
\begin{align*}
R_x:\Omega^r(\overline{W}_p G) & \xymatrix{ \ar[r] & \Omega^r(\overline{W}_{p-1}G)} & J_x:\Omega^r(\overline{W}_pG) & \xymatrix{ \ar[r] & \Omega^{r-1}(\overline{W}_{p-1}G)} \\
	\omega  & \xymatrix{\ar@{|->}[r] & \bar{\sigma}_0^*\Lie_{x^p}\omega ,}     &    \omega  & \xymatrix{\ar@{|->}[r] & \bar{\sigma}_0^*\iota_{x^p}\omega .}
\end{align*}
For $y\in\hh$, recall that $y_{\beta}^p\in\fX(\overline{W}_p G)$ denotes the infinitesimal generator of the action of Proposition \ref{prop:lie} and define the operators
\begin{align*}
R_y:\Omega^r(\overline{W}_pG) & \xymatrix{ \ar[r] & \Omega^r(\overline{W}_{p-2}G)} &  J_y:\Omega^r(\overline{W}_pG) & \xymatrix{ \ar[r] & \Omega^{r-1}(\overline{W}_{p-2}G)} \\ 
\omega & \longmapsto \sum_{\beta=0}^{p-2}(-1)^\beta(\bar{\sigma}_0\circ\bar{\sigma}_\beta)^*\Lie_{y_{\beta+1}^p}\omega &  \omega & \longmapsto \sum_{\beta=0}^{p-2}(-1)^\beta(\bar{\sigma}_0\circ\bar{\sigma}_\beta)^*\iota_{y_{\beta+1}^p}\omega.
\end{align*}

For the statement of the main theorem, we use the following convention. One can always write a permutation $\varrho$ in the symmetric group $S_{k+l+a+b}$ as a composition $\varrho=\varrho'\circ\varrho^w\circ\varrho^z\circ\varrho^y\circ\varrho^x$, where $\varrho^x$ permutes only the first $k$ many elements, $\varrho^y$ permutes the next $l$ many elements, $\varrho^z$ permutes the next $a$ many elements, $\varrho^w$ permutes only the last $b$ many elements, and $\varrho'$ is the residual permutation. Then, we define the special sign 
$$\abs{\abs{\varrho}}:=\textnormal{sgn}(\varrho^x)\textnormal{sgn}(\varrho^w)\textnormal{sgn}(\varrho').$$
\begin{thm}\label{Theo:CommutesWithDiffs}
    There is a map of double-complexes 
    $$\Phi:\left(\widehat{\Omega}^{\bullet,\bullet}(G), \bar{\partial}, d\right)\longrightarrow \left(\cW^{\bullet,\bullet}(\h\to \g_0), \delta, d\right),$$
    whose  $\Phi^{klab}\omega\in\bigwedge^k\gg_0^*\otimes \sym^l\hh^*\otimes \sym^a\gg_0^*\otimes\bigwedge^b\hh^*$ component is defined by 
\begin{eqnarray*}
(\Phi^{klab}\omega)
(\mathbbm{x};\mathbbm{y}; \mathbbm{z};\mathbbm{w})=\sum_{\varrho\in S_{k+l+a+b}}\abs{\abs{\varrho}}\varrho\Big{(}R_{x_1}\cdots R_{x_k}R_{y_1}\cdots R_{y_l}J_{z_1}\cdots J_{z_a}J_{w_1}\cdots J_{w_b}\Big{)}\omega
\end{eqnarray*} 
where $\mathbbm{x}=(x_1,...,x_k)$, $\mathbbm{y}=(y_1,...,y_l)$, $\mathbbm{z}=(z_1,...,z_a)$ and  $\mathbbm{w}=(w_1,...,w_b)$ for $x_\alpha,z_m\in\gg_0$ for $1\leq\alpha\leq k$ and $1\leq m\leq a$, and for  $y_\beta,w_n\in\hh$ for $1\leq\beta\leq l$ and $1\leq n\leq b$. We call $\Phi$ the {\bf van Est map}.
\end{thm}

We dedicate the following subsections to setting up the proof of Theorem~\ref{Theo:CommutesWithDiffs}. As outlined in the introduction, the van Est map $\Phi$ is built as a composition of the map \eqref{Eq:partialZero} and the homotopy inverse to the inclusion; schematically,  
$$\hat{\Omega}(G)\longrightarrow \mathcal{C}\xleftarrow{(\ast)}\mathcal{W}(\hh\to\gg_0).$$
In building it so, $\Phi$ is a map of complexes by construction, but we are left to prove that its formula is the one given in the statement of Theorem~\ref{Theo:CommutesWithDiffs}. In Subsection~\ref{sec:hom} below, to set up the necessary ingredients to use the Perturbation lemma, we describe the intermediary complex $\mathcal{C}$ in Eq.~\eqref{Eq:AuxDoubleCxForms}, as well as the homotopy $h$ in Eq.~\eqref{Eq:h} used to produce the homotopy inverse to $(\ast)$ in the schematic diagram above (see Propositions~\ref{Prop:prePertLemma} and \ref{Prop:FormPerturb}). We proceed then by steps, proving that our formula works at the level of cochains first, in Subsection~\ref{sec:func}, and, then at the level of forms in Subsection~\ref{sec:proof}. The resulting formulas for the van Est map at the level of cochains and forms are given in Theorems~\ref{Theo:vanEstAtCochains} and \ref{Theo:vanEstAtForms}, respectively.

\subsection{The homotopy}\label{sec:hom}

The map $\iota_0^*$ in the proof of Theorem~\ref{Theo:SimplicialDiffn} can be seen as dual to the inclusion of $\hh\longrightarrow\gg_0$ into $\H_0\longrightarrow\A_0$; however, without restricting to $E_1^{0,q}$, $\iota_0^*$ is in general 
not a cochain map. This echoes the fact that, as in the group case, this inclusion is in general not a map of Lie $2$-algebroids. Notwithstanding, as in the group case, there exists a homotopy operator $h$ for the simplicial differential $\partial$ on $C^{p,q}$ and the 
following result holds.

\begin{prop}\label{Prop:Perturb}
     The composition 
     $$\iota_0^*\circ (1+\delta\circ h)^{-1}:\textnormal{Tot}^\bullet(C^{p,q})\longrightarrow C_{\textnormal{CE}}^\bullet(\h\to\g_0),$$ 
     where $\delta$ is the Chevalley-Eilenberg differential, is a cochain map homotopy inverse to $\pi_0^*$.
\end{prop}

The proof of this result is given by applying the Perturbation lemma of homological algebra 
\cite{Brown:1965, Gugenheim:1972}, the setup of which is the following: Consider the constant simplicial 
Lie 2-algebra $\hh\longrightarrow\gg_0$ all of whose face and degeneracy maps are identities, and let 
$$D^{p,q}:=C_{\textnormal{CE}}^p(\hh\rightarrow\gg_0)$$ be its associated double complex.

For each $p>0$, there is a projection 
\begin{eqnarray}\label{Eq:Proj}
\begin{array}{c}
    \xymatrix{
    \H_p \ar[d]\ar[r]^{\pi_p} & \hh: \ar[d] &  (Y,\vec{\gamma}) \ar@{|->}[r] & \hat{\partial}_2\circ ...\circ\hat{\partial}_{p+1}(Y\star 0_{s_0(\pr_p\vec{\gamma})}^{-1})\in\hh \\
    \A_p \ar[r]^{\pi_p} & \gg_0: & (X,\vec{\gamma}) \ar@{|->}[r] & \hat{\partial}_1\circ\hat{\partial}_2\circ ...\circ\hat{\partial}_{p+1}(X\star 0_{s_0(\pr_p\vec{\gamma})}^{-1})\in\gg_0
    }
\end{array}
\end{eqnarray}
that extends the projection $\pi_0$ of Eq.~\eqref{Eq:piZero} in the sense that they assemble 
into a morphism of simplicial Lie 2-algebroids.

\begin{lem}\label{Lemma:pi'sGood}
    The maps~\eqref{Eq:Proj} define a strict map of simplicial strict Lie 2-algebroids.
\end{lem}

\begin{proof}
    The proof splits in two: That the maps respect the Lie 2-algebroid structure is equivalent to their duals being morphisms of 
    Chevalley-Eilenberg complexes. As in the proof of Theorem~\ref{Theo:SimplicialDiffn}, this reduces 
    to checking that the differentials commute in the lowest degrees. Because the face maps in the 
    simplicial structure of $D^{p,q}$ are identical, being a simplicial map boils down to proving that $\partial\circ\pi_p^*$ equals zero for even $p$, and equals the identity for odd $p$. The proof is analogous to that for $\pi_0^*$ in Theorem~\ref{Theo:SimplicialDiffn}.
\end{proof}

$D^{p,q}$ can be thus seen as a subcomplex of $C^{p,q}$ as $\pi_p^*$ admits a left inverse $\iota_p^*$, 
for each $p$. These latter are dual to the inclusions of $\hh\longrightarrow\gg_0$ in 
$\H_p\longrightarrow\A_p$ that is, for $p>0$,
\begin{eqnarray}\label{Eq:Incl}
\begin{array}{c}
   \xymatrix{
    \hh \ar[d]\ar[r]^{\iota_p} & \H_p: \ar[d] & \hh\ni y \ar@{|->}[r] & (\hat{\sigma}_1^{p}(y),(\mathbbm{1}_p;...;1)) \in\H_p  \\
    \gg_0 \ar[r]^{\iota_p} & \A_p: & \gg_0\ni x \ar@{|->}[r] & (\hat{\sigma}_0^{p+1}(x),(\mathbbm{1}_p;...;1))\in\A_p
    }
\end{array}
\end{eqnarray}
This map is in general not a map of Lie 2-algebroids and, hence, its dual is in general not a map of 
cochain complexes. However, the canonical simplicial deformation retraction of 
Eq.'s~\eqref{Eq:ContrHomotopiesF} and \eqref{Eq:ContrHomotopiesH} can be extended to a homotopy between 
the identity and $\pi_\bullet^*\circ\iota_\bullet^*$. For $1\leq i\leq p$, let $h^{\A}_i:\A_p \to \A_{p+1}$ be given by
\begin{equation}\label{Eq:HomotopiesF}
    h^{\A}_i(X, \vec{\gamma}):=\left(\hat{\sigma}_0^{i+1}\circ\hat{\partial}_{p+2-i}\circ ...\circ\hat{\partial}_{p+1}(X\star 0_{s_0(\pr_p\vec{\gamma})}^{-1}), \sigma_0^{i+1}\circ\partial_{p+1-i}\circ ...\circ\partial_{p}(\vec{\gamma})\right)
\end{equation}
and let $h^{\H}_i:\H_p \to \H_{p+1}$ be given by
$$h^{\H}_i(Y,\vec{\gamma}):=\left(\hat{\sigma}_{1}(Y\star 0_{s_0(\pr_p\vec{\gamma})}^{-1}),\sigma_0^{i+1}\circ\partial_{p+1-i}\circ ...\circ\partial_{p}(\vec{\gamma})\right).$$
These maps are well-defined because $\varepsilon_\bullet$ is a simplicial map. The homotopy operator is thus given by,
\begin{eqnarray}\label{Eq:h}
    h:=\sum_{i=0}^{p}(-1)^{i}h_i^*:\xymatrix{C^{p+1,q} \ar[r] & C^{p,q}.}
\end{eqnarray}
More precisely, for $\omega\in \Gamma\left(\bigwedge\nolimits^k\cA^*_{p+1}\otimes\sym^l\cH_{p+1}\right)$ with $k+2l=q$ then
\begin{align*}
    (h\omega)_{\vec{\gamma}}(\mathbb{X};\mathbb{Y}) & =\sum_{i=0}^{p}(-1)^{i}\omega_{\sigma_0^{i+1}(\partial_{p+1-i}\circ ...\circ\partial_{p}(\vec{\gamma}))}(h^{\A}_i(X_1, \vec{\gamma}),...,h^{\A}_i(X_{k}, \vec{\gamma});h^{\H}_i(Y_1,\vec{\gamma}),...,h^{\H}_i(Y_l,\vec{\gamma})),
\end{align*}
for $\vec{\gamma}\in W_pG$, $\mathbb{X}=((X_1,\vec{\gamma}),...,(X_{k},\vec{\gamma}))$ and $\mathbb{Y}=((Y_1,\vec{\gamma})...,(Y_l,\vec{\gamma}))$, where $(X_m, \vec{\gamma})\in \A_{p}$ for each $1\leq m\leq k$ and $(Y_n, \vec{\gamma})\in\H_{p}$ for each $1\leq n\leq l$.

\begin{prop}\label{Prop:prePertLemma}
    $h$ preserves $D^{p,q}$, commutes with $\delta$ on $D^{p,q}$ and is a homotopy between the identity and $\pi_\bullet^*\circ\iota_\bullet^*$, that is 
    $$\partial\circ h+h\circ\partial =\textnormal{Id}-\pi_\bullet^*\circ\iota_\bullet^*.$$
\end{prop}

\begin{proof}
    That $h$ preserves $D^{p,q}$ follows from the identities $\pi_{p+1}\circ h_i^\A=\pi_p$ and $\pi_{p+1}\circ h_i^\H=\pi_p$, which hold for all $1\leq i\leq p$ because
    \begin{align*}
        \pi_{p+1}(h_i^\A(X, \vec{\gamma})) & =\hat{\partial}_1\circ ...\circ\hat{\partial}_{p+2}\circ\hat{\sigma}_0^{i+1}\circ\hat{\partial}_{p+2-i}\circ ...\circ\hat{\partial}_{p+1}(X\star 0_{\sigma_0(\gamma^p)}^{-1}) =\pi_{p+1}(X, \vec{\gamma})
    \end{align*}
    and
    \begin{align*}
        \pi_{p+1}(h_i^\H(Y,\vec{\gamma})) & =\hat{\partial}_2\circ ...\circ\hat{\partial}_{p+2}(\hat{\sigma}_{1}(Y\star 0_{\sigma_0(\gamma^p)}^{-1}))=\pi_p(Y,\vec{\gamma}).
    \end{align*}
    That $h$ commutes with $\delta$ on $D^{p,q}$ is equivalent to $h$ restricting to a map of complexes and, 
    thus, to its dual defining a map of Lie 2-algebras. 
    Let $x\in\gg_0$, then $\iota_p(x)=x^{p+1}_{(\mathbbm{1}_p;...;1)}$, where $x^{p+1}\in\Gamma(\A_p)$. Hence, since for all $2\leq p+2-i\leq k\leq p$, $\hat{\partial}_k\circ\hat{\sigma}_0^k=\hat{\sigma}_0^{k-1}$, then 
    \begin{eqnarray*}
        h_i^\A(x^{p+1})=x^{p+2}\in(\A_{p+1})_{(\mathbbm{1}_{p+1};...;1)}.
    \end{eqnarray*}
    Consequently, $\pi_{p+1}(x^{p+2})=x$, and  thus $\pi_{p+1}\circ h_i^\A\circ\iota_p=\textnormal{Id}_{\gg_0}$. Ultimately, this implies that the restriction of $h$ to $D^{p,1}$, 
    \begin{eqnarray*}
        \iota_p^*\circ h\circ\pi_{p+1}^*=\sum_{i=0}^p(-1)^i\iota_p^*\circ (h_i^\A)^*\circ\pi_{p+1}^*=\sum_{i=0}^p(-1)^i\textnormal{Id}_{\gg_0^*},
    \end{eqnarray*}
    which is either zero or the identity depending on the parity of $p$. Analogously, for $y\in\hh$.  Lastly,  that $h$ yields the desired homotopy follows from combining Lemma~\ref{lemma:h0} and Proposition~\ref{Prop:h'}.   
\end{proof}
\begin{proof}[Proof of Proposition~\ref{Prop:Perturb}]
    In sight of Lemma~\ref{Lemma:pi'sGood} and Proposition~\ref{Prop:prePertLemma}, we can invoke the Perturbation 
    Lemma of homological algebra (see e.g. \cite[Lemma B.3 in Appendix B]{Meinrenken_Li-Bland:2015} ) to conclude that 
    $\iota_\bullet^*\circ(1+\delta\circ h)^{-1}:C^{p,q}\longrightarrow D^{p,q}$ is a cochain map for the total 
    differential $\delta +\partial$, and is a homotopy equivalence, with homotopy inverse 
    $(1+h\circ\delta)^{-1}\circ\pi_\bullet^*$. Now, the inclusion of 
    $C_{\textnormal{CE}}^\bullet(\hh\rightarrow\gg_0)=D^{0,\bullet}$ in $D^{p,q}$ is also a homotopy equivalence, 
    with homotopy inverse the projection. Here, the corresponding homotopy operator, as in the proof of 
    Proposition~\ref{Prop:prePertLemma}, is the identity in even degrees, and 0 otherwise. Composing the two 
    homotopy equivalences, and observing that for degree reasons $(1+h\circ\delta)^{-1}\circ\pi_0^*=\pi_0^*$, 
    the result follows.
\end{proof}

As outlined in the introduction, we define the van Est map as the composition of the map from Proposition~\ref{Prop:Perturb} and the map $\partial_0^*$ from Eq.~\eqref{Eq:partialZero}. Thus far, however, our considerations only apply to the level of cochains; to extend them to include differential forms, we replace each occurrence of a Chevalley-Eilenberg complex by a Weil algebra. To make these fit into the setup of the Perturbation lemma, we replace the auxiliary double complex $C^{p,q}$ in \eqref{Eq:AuxDoubleCx} by 
\begin{eqnarray}\label{Eq:AuxDoubleCxForms}
    \cC^{p,q}:=\bigoplus_{a+b=p}\cW^{q,b}(\H_a\rightarrow\A_a).
\end{eqnarray}

Observe that $\cW^{0,r}(\H_p\rightarrow\A_p)=\Omega^r(W_p\G)$ and that pullback along the d\'ecalage projection yields a map of complexes
\begin{eqnarray}\label{Eq:partialZeroSuite}
    \xymatrix{\Big(\textnormal{Tot}^{\bullet}\ \widehat{\Omega}^{p,r}(G), \bar{\partial}+d\Big) \ar[r]^{\bar{\partial}_0^*\qquad} & \Big(\textnormal{Tot}^{\bullet}\ \Omega^r(W_p\G)=\cC^{\bullet, 0}, \partial+d\Big).}
\end{eqnarray}
The retracting subcomplex $\cD^{p,q}$ is defined accordingly as the contraction of the triple complex of Weil algebras of the constant simplicial Lie $2$-algebra $\hh\longrightarrow\gg_0$, where the contraction is taken as in~\eqref{Eq:AuxDoubleCxForms},
\begin{eqnarray}\label{Eq:AuxDoubleSubCxForms}
    \cD^{p,q}:=\bigoplus_{a+b=p}\cW^{q,b}(\hh\rightarrow\gg_0)_a.
\end{eqnarray}
Note that the subindex $a$ is mute here, as the Weil algebra does not change as $a$ increases.

The maps ~\eqref{Eq:Proj} and \eqref{Eq:Incl} remain unchanged but the homotopy operator of Eq.~\eqref{Eq:h} needs a replacement too. To define
$$h:\cW^{q,r}(\cH_{p+1}\to\cA_{p+1})\longrightarrow \cW^{q,r}(\cH_{p}\to\cA_{p}), $$
we use the pullback along the degeneracy map $\bar{\sigma}_0$ to extend $h_0^*$ to account for differential forms. Then, for 
$\omega\in \Gamma\left(\bigwedge\nolimits^{s}T^*(W_{p+1}G)\otimes\bigwedge\nolimits^k \cA^*_{p+1}\otimes \sym^l\cH^*_{p+1}\otimes\sym^a \cA^*_{p+1}\otimes\bigwedge\nolimits^b \cH_{p+1}^*\right)$, with $a+b+s=r$ and $k+a+2(l+b)=q$,
we define  
\begin{align}\label{Eq:hFinal}
    (h\omega)_{\vec{\gamma}}(\mathbb{X};\mathbb{Y};\mathbb{Z};\mathbb{W}\vert V)=\sum_{i=0}^{p}&(-  1)^{i}\omega_{\bar{h}_i(\vec{\gamma})}(h^{\A}_i(X_1, \vec{\gamma}),...,h^{\A}_i(X_{k},\vec{\gamma});h^{\H}_i(Y_1,\vec{\gamma}),...,h^{\H}_i(Y_l,\vec{\gamma}); \\
         & \qquad h^{\A}_i(Z_1,\vec{\gamma}),...,h^{\A}_i(Z_{k},\vec{\gamma});h^{\H}_i(W_1,\vec{\gamma}),...,h^{\H}_i(W_l,\vec{\gamma})\vert T\bar{h}_i(V)), \nonumber
\end{align}
where $\bar{h}_i=\bar{\sigma}_0^{i+1}\circ\bar{\partial}_{p+2-i}\circ ...\circ\bar{\partial}_{p+1}$, $\vec{\gamma}\in W_{p} G$, $V\in\bigwedge^sT_{\vec{\gamma}}(W_{p}\G)$,   $\mathbb{X}$ and $\mathbb{Y}$ are as in \eqref{Eq:h} and $\mathbb{Z}=((Z_1, \vec{\gamma})...,(Z_a, \vec{\gamma}))$ and $\mathbb{W}=((W_1, \vec{\gamma})...,(W_b, \vec{\gamma}))$ for $(Z_m,\vec{\gamma})\in\A_{p}$ with $1\leq m\leq a$ and $(W_n, \vec{\gamma})\in\H_{p}$ with $1\leq n\leq b$. 

Mimicking the proofs of Proposition \ref{Prop:Perturb} and \ref{Prop:prePertLemma}, we get. 

\begin{prop}\label{Prop:FormPerturb}
\
\begin{enumerate}
    \item The homotopy $h$ preserves $\cD^{p,q}$, commutes with $\delta$ on $\cD^{p,q}$ and is a homotopy between the identity and $\pi_\bullet^*\circ\iota_\bullet^*$.
    \item Then, the composition 
     $\iota_\bullet^*\circ (1+\delta\circ h)^{-1}:\cC^{\bullet,\bullet}\longrightarrow \cW^{\bullet,\bullet}(\hh\rightarrow\gg_0),$ 
     where $\delta$ is the Chevalley-Eilenberg differential, is a morphism of bidifferential spaces and a homotopy inverse to $\pi_0^*$ with respect to $\partial+d$, where $d$ is the de Rham differential of the Weil algebra.
\end{enumerate}
\end{prop}

\subsection{The van Est map at the level of functions}\label{sec:func}

In this section, we give a proof of Theorem \ref{Theo:CommutesWithDiffs} at the level of cochains. Our proof makes use of the simplicial differentiation explained in $\S\ref{sec:dif}$ and the construction outlined in Subsection $\S\ref{sec:hom}$

\begin{nota}\label{not:vE}
    Let $k$ and $l$ be integers such that $k+2l=p$, $x_1,...,x_k\in\gg_0$ and $y_1,...,y_l\in\hh$. Let's assume the following notational conventions:
    \begin{itemize}
        \item We write $\mathbbm{x}:=(x_1,...,x_k)$ and $\mathbbm{y}=(y_1,...,y_l)$;
        \item For sequences of integers $1\leq a_1<...<a_r\leq k$ and $1\leq b_1<...<b_s\leq l$,
        \begin{align*}
            \mathbbm{x}(a_1,...,a_r) & :=(x_1,...,x_{{a_1}-1}, x_{{a_1}+1},...,x_{{a_r}-1}, x_{{a_r}+1},...,x_k), \\
            \mathbbm{y}(b_1,...,b_s) & :=(y_1,...,y_{{b_1}-1}, y_{{b_1}+1},...,y_{{b_s}-1}, y_{{b_s}+1},...,y_l);
        \end{align*}
        \item For $1\leq\beta\leq r\leq p$, we write $\mathbbm{x}^r$ and $\mathbbm{y}^r_\beta$, respectively, for the $k$-tuple and the $l$-tuple of right-invariant sections  
        \begin{align*}
            \mathbbm{x}^r & :=(x_1^r,...,x_k^r), & \mathbbm{y}^r_\beta & :=((y_1)^r_\beta,...,(y_l)^r_\beta),
        \end{align*}
        where $x_a^r\in\Gamma(\A_{r-1})$ and $(y_b)^r_\beta\in\Gamma(\H_{r-1})$ as in the proof of Proposition \ref{prop:A}. 
    \end{itemize}
\end{nota}

\begin{thm}\label{Theo:vanEstAtCochains}
    The composition of the map $\iota_0^*\circ (1+\delta\circ h)^{-1}$ from Proposition~\ref{Prop:Perturb} and the map $\bar{\partial}_0^*$ from Eq.~\eqref{Eq:partialZero} yields the van Est map $\Phi$ at the level of cochains.
\end{thm}

\begin{proof}
    First, observe that for a normalized cochain $f\in \widehat{C}^\infty(\overline{W}_pG)$, $h\bar{\partial}_0^*f=f$. Since $\delta\circ h$ has bidegree $(-1,1)$, the map from Proposition~\ref{Prop:Perturb} can be computed using the following ‘zig-zag’ form: $\iota_0^*\circ(\delta\circ h_0^*)^p\bar{\partial}_0^*f$,
    where the $h_i$'s for $i>0$ vanish because $f$ is a normalized cochain.

    As in the proof of Theorem~\ref{Theo:SimplicialDiffn}, $\iota_0^*$ can be computed using  
    \begin{align}\label{Eq:0thStep}
        (\iota_0^*\circ(\delta\circ h_0^*)^p\bar{\partial}_0^* f)(\mathbbm{x};\mathbbm{y}) & =\bar{\sigma}_0^*\Big{(}(\delta\circ(h_0^*\circ\delta)^{p-1})f(\mathbbm{x}^1;\mathbbm{y}^1)\Big{)}.
    \end{align}
    When computing the Chevalley-Eilenberg differential, there are four types of terms: 
    \begin{equation}\label{eq:I-IV}\def\arraystretch{1.3}
        \begin{array}{rcl}
            \textnormal{I.} & \Lie_{x_a^1}(h_0^*\circ\delta)^{p-1}f(\mathbbm{x}^1(a);\mathbbm{y}^1) & \text{ for } 1\leq a\leq k\\
            \textnormal{II.} & (h_0^*\circ\delta)^{p-1}f(\mathbbm{x}^1(a);[x_{a}^1,y_b^1],\mathbbm{y}^1(b)) &  \text{ for } 1\leq a\leq k \text{ and } 1\leq b\leq l\\
            \textnormal{III.}& (h_0^*\circ\delta)^{p-1}f([x_{a_1}^1,x_{a_2}^1],\mathbbm{x}^1(a_1,a_2);\mathbbm{y}^1) & \text{ for } 1\leq a_1<a_2\leq k\\
           \textnormal{IV.} &(h_0^*\circ\delta)^{p-1}f(\mathbbm{x}^1,y_b^1;\mathbbm{y}^1(b))&\text{ for } 1\leq b\leq l
        \end{array}
    \end{equation}

    We refer to the terms of type II and III as {\it bracket terms}; these will ultimately vanish due to $f$ being normalized, so we argue separately for these. Note that composing with the degeneracy of Eq.~\eqref{Eq:0thStep}, the terms of type I and IV become respectively into:
    \begin{equation*}
    \def\arraystretch{1.3}
    \begin{array}{rcl}
        1. & R_{x_a}(h_0^*\circ\delta)^{p-1}f(\mathbbm{x}^1(a);\mathbbm{y}^1) & \text{ for } 1\leq a\leq k\\
        2. & \bar{\sigma}_0^*(h_0^*\circ\delta)^{p-1}f(\mathbbm{x}^1,y_b^1;\mathbbm{y}^1(b)) & \text{ for } 1\leq b\leq l.
    \end{array}
    \end{equation*}
    We collect the sum of these terms and record them using the following table 
    \begin{eqnarray}\label{Eq:Graph1}
    \begin{array}{|c|c|}
    \hline
       (k \ 0 \ l )  & \begin{array}{c}
            \underset{k-1,0,l}{x}  \\
            \hline
            \underset{k,1,l-1}{0_y}
       \end{array} \\
         \hline
    \end{array}
     \end{eqnarray}
     Here, the under-labels are triples that represent the number of arguments as follows: 
    The under-label $k',l',l''$ indicates that the argument takes values in $\Gamma(\bigwedge^{k'+l'}\A_\bullet\otimes \sym^{l''}\H_\bullet)$; it further indicates that the first $k'$ many arguments come from sections of $\A_\bullet$ and the subsequent $l'$ many come from sections of $\H_\bullet$ regarded via the inclusion $j_\bullet$ as sections of $\A_\bullet$. As for the larger labels, $x$ represents the sum $\sum_{a=1}^k(-1)^{a-1}R_{x_a}$ of terms of type 1, and the larger label $0_y$ represents the sum $(-1)^k\sum_{b=1}^l\bar{\sigma}_0^*$ of terms of type 2. 
    
    One can compute $h_0^*$ using right-invariant sections and evaluating at degeneracies. In doing so, the terms of type 1 and 2 respectively give back: 
    \begin{equation*}
        \begin{array}{rcl}
        1'. & R_{x_a}\bar{\sigma}_0^*\Big{(}(\delta\circ(h_0^*\circ\delta)^{p-2})f(\mathbbm{x}^2(a);\mathbbm{y}^2_1)\Big{)} & \text{ for } 1\leq a\leq k\\
        2'. & (\bar{\sigma}_0^*)^2\Big{(}(\delta\circ(h_0^*\circ\delta)^{p-2})f(\mathbbm{x}^2,(y_b)^2_2;\mathbbm{y}^2_1(b))\Big{)} & \text{ for } 1\leq b\leq l.
    \end{array}
    \end{equation*}
    Computing the subsequent Chevalley-Eilenberg differential, and neglecting the bracket terms, 
    yields:
    \begin{equation*}\label{eq:1.i}
     \def\arraystretch{1.3}
        \begin{array}{rcl}
        1.1 & R_{x_{a_1}}R_{x_{a_2}}(h_0^*\circ\delta)^{p-2}f(\mathbbm{x}^2(a_1,a_2);\mathbbm{y}^2_1) & \text{ for } 1\leq a_1,a_2\leq k\\
        1.2 & R_{x_a}\bar{\sigma}_0^*(h_0^*\circ\delta)^{p-2}f(\mathbbm{x}^2(a),(y_b)^2_1;\mathbbm{y}^2_1(b)) & \text{ for } 1\leq a\leq k\quad \text{and}\quad 1\leq b\leq l.
    \end{array}
    \end{equation*}
    from the terms of type 1', and 
\begin{equation*}\label{eq:2.i}
     \def\arraystretch{1.3}
        \begin{array}{rcl}
        2.1 & R_{x_a}\bar{\sigma}_1^*(h_0^*\circ\delta)^{p-2}f(\mathbbm{x}^2(a),(y_b)^2_2;\mathbbm{y}^2_1(b)) & \text{ for } 1\leq a\leq k\quad \text{and}\quad 1\leq b\leq l\\
        2.2 & R_{y_b}(h_0^*\circ\delta)^{p-2}f(\mathbbm{x}^2;\mathbbm{y}^2_1(b)) & \text{ for } 1\leq b\leq l\\
        2.3 & (\bar{\sigma}_0^*)^2(h_0^*\circ\delta)^{p-2}f(\mathbbm{x}^2,(y_{b_1})^2_2,(y_{b_2})^2_1;\mathbbm{y}^2_1(b_1,b_2))& \text{ for } 1\leq b_1,b_2\leq l
    \end{array}
    \end{equation*}
Note that the shape terms in $2.1$ assume follows from the simplicial identities and the fact that $x_a^1$ is $\bar{\sigma}_1$-related to $x_a^2$. Indeed,
\begin{align*}
        (\bar{\sigma}_0\circ\bar{\sigma}_0)^*\Lie_{x_a^2}(h_0^*\circ\delta)^{p-2}f(\mathbbm{x}(a)^2,(y_b)^2_2;\mathbbm{y}^2_1(b)) & =(\bar{\sigma}_1\circ\bar{\sigma}_0)^*\Lie_{x_a^2}(h_0^*\circ\delta)^{p-2}f(\mathbbm{x}(a)^2,(y_b)^2_2;\mathbbm{y}^2_1(b)) \\
            & =\bar{\sigma}_0^*\circ\bar{\sigma}_1^*\Lie_{x_a^2}(h_0^*\circ\delta)^{p-2}f(\mathbbm{x}(a)^2,(y_b)^2_2;\mathbbm{y}^2_1(b)) \\
            & =\bar{\sigma}_0^*\Lie_{x_a^1}\bar{\sigma}_1^*(h_0^*\circ\delta)^{p-2}f(\mathbbm{x}(a)^2,(y_b)^2_2;\mathbbm{y}^2_1(b)).
        \end{align*}
        On the other hand, the shape terms on $2.2$ assume follows because the anchor of $(y_b)^2_2\in\Gamma(\A_1)$ being $(y_b)^2_1\in\XX(\overline{W}_2G)$, thus implying 
        \begin{align*}
        (\bar{\sigma}_0^*)^2\Lie_{\rho((y_b)^2_2)}(h_0^*\circ\delta)^{p-2}f(\mathbbm{x}^2;\mathbbm{y}^2_1(b)) & =(\bar{\sigma}_0\circ\bar{\sigma}_0)^*\Lie_{(y_b)^2_1}(h_0^*\circ\delta)^{p-2}f(\mathbbm{x}^2;\mathbbm{y}^2_1(b)).
        \end{align*} 
Recording these terms in the table~\eqref{Eq:Graph1} induces the following branching:
    \begin{eqnarray}\label{Eq:Graph2}
    \begin{array}{|c|c|}
    \hline
      (k\ 0 \ l)   & \begin{array}{c|c}
         \underset{k-1,0,l}{x}  & \begin{array}{c}
              \underset{k-2,0,l}{xx}  \\
              \hline 
             \underset{k-1,1,l-1}{x0_y} 
         \end{array} \\
         \hline 
         \underset{k,1,l-1}{0_y}  & \begin{array}{c}
               \underset{k-1,1,l-1}{x1_y} \\
               \hline
             \underset{k,0,l-1}{y^1}\\
             \hline
             \underset{k,2,l-2}{0_y0_y}
         \end{array}
      \end{array} \\
      \hline
    \end{array}
    \end{eqnarray}
    Here, the label $xx$ represents 
    \begin{equation}\label{eq:xx}
        \sum_{a_1=1}^k(-1)^{a_1-1}R_{x_{a_1}}\left(\sum_{a_2=1}^{a_1-1}(-1)^{a_2-1}R_{x_{a_2}}+\sum_{a_2=a_1+1}^{k}(-1)^{a_2}R_{x_{a_2}}\right)
    \end{equation}
     summing over terms of type 1.1; the labels $x0_y$ and $x1_y$ stand respectively for 
     \begin{equation}\label{eq:x0ax1}
         \sum_{a=1}^k(-1)^{(a-1)+(k-1)}R_{x_a}\sum_{b=1}^l\bar{\sigma}_0^*\quad\text{and}\quad \sum_{b=1}^l\sum_{a=1}^k(-1)^{(a-1)+k}R_{x_a}\bar{\sigma}_1^*
     \end{equation} summing over terms of type 1.2 and type 2.1 respectively; and, lastly, the labels $y^1$ and $0_y0_y$ stands respectively for   
     \begin{equation}\label{eq:00}
        (-1)^k\sum_{b=1}^lR_{y_b}\quad\text{and}\quad(-1)^{k+(k+1)}\sum_{b_1=1}^l\sum_{b_2=1}^{l-1}(\bar{\sigma}_0^*)^2    
     \end{equation}
     summing over terms of type 2.2 and of type 2.3, respectively. 
    
    In general, $\alpha:=((k-k')+2(l-l'')-l')$ many steps are required to reach a vertex under-labelled $k',l',l''$, with $0\leq k'\leq k$ and $0\leq l'+l''\leq l$. This process corresponds to computing a type of term in the $\alpha$-th iteration of $h_0^*\circ\delta$. 
    Let the larger label on $k',l',l''$ be $\overrightarrow{P}(d_1)_y...(d_{l'})_y$, where $\overrightarrow{P}:=\varrho(x...xy^{\beta_1}...y^{\beta_{l-(l'+l'')}})$ for some $\varrho\in S_{\alpha-l+l''}$, and $k-k'\geq d_1\geq ...\geq d_{l'}\geq 0$. This label represents a sum over $1\leq b_1,...,b_{l'}\leq l$ and an alternating sum over $1\leq a_1,...,a_{k-k'}\leq k$ of the operators 
    $$\varrho\big{(}R_{x_{a_1}}...R_{x_{a_{k-k'}}}R^{(\beta_1)}_{y_{b_1}}...R^{(\beta_{l-(l'+l'')})}_{y_{b_{l-(l'+l'')}}}\big{)}(\bar{\sigma}_{d_1}^*\circ...\circ\bar{\sigma}_{d_{l'}}^*),$$ 
    where $R^{(\beta)}_{y}:=(\bar{\sigma}_0\circ\bar{\sigma}_{\beta-1})^*\Lie_{y^\bullet_\beta}$ (which evidences there are some restrictions on the values $\beta_i$ can assume depending on the permutation), evaluated at
    \begin{align*}
        (h_0^*\circ\delta)^{p-\alpha}f\big{(}\mathbbm{v}^\alpha_{d_1...d_{l'}}\big{)},
    \end{align*}
    where, setting $b'_i:=b_{l-(l'+l'')+i}$, 
       $$\mathbbm{v}^\alpha_{d_1...d_{l'}}=\left(\mathbbm{x}^\alpha(a_1,...,a_{k-k'}),(y_{b'_1})^\alpha_{d_1+l'},(y_{b'_2})^\alpha_{d_2+l'-1},...,(y_{b'_{l'-1}})^\alpha_{d_{l'-1}+2},(y_{b'_{l'}})^\alpha_{d_{l'}+1};\mathbbm{y}^\alpha_1(b_1,...,b_{l-l'})\right).$$

    Inductively, one computes the homotopy as before, right-extending sections and evaluating at degeneracies, thus yielding the same total sum, though evaluated at
    \begin{align}\label{Eq:oneStepHomotopy}
        \bar{\sigma}_0^*\Big{(}(\delta\circ(h_0^*\circ\delta)^{p-\alpha-1})f(\mathbbm{v}^{\alpha+1}_{(d_1+1)...(d_{l'}+1)} )\Big{)}.
    \end{align}
    Disregarding the bracket terms again, there are three types of terms coming out of computing the Chevalley-Eilenberg differential:

\begin{itemize}
    \item     Type 1: If $k'>0$, the first $k'$ many terms can be written as an alternating sum over $1\leq a_{k-k'+1}\leq k'$ of operators $\Lie_{x^{\alpha+2}_{a_{k-k'+1}}}$ evaluated at $(h_0^*\circ\delta)^{p-\alpha-1}f\big{(}\mathbbm{v}^{\alpha+1, I(a_{k-k'+1})}_{(d_1+1)...(d_{l'}+1)}\big{)}$, where
    $$\mathbbm{v}^{\alpha,I(a_{k-k'+1})}_{d_1...d_{l'}}=\left(\mathbbm{x}^\alpha(a_1,...,a_{k-k'},a_{k-k'+1}),(y_{b'_1})^\alpha_{d_1+l'},...,(y_{b'_{l'}})^\alpha_{d_{l'}+1};\mathbbm{y}^\alpha_1(b_1,...,b_{l-l'})\right).$$
    Indeed, since $0\leq d_i$ for each $1\leq i\leq l'$, the simplicial identities imply that
    \begin{align*}
        \bar{\sigma}_{d_1}^*\circ...\circ\bar{\sigma}_{d_{l'}}^*\circ\bar{\sigma}_0^*  
          =\bar{\sigma}_0^*\circ\bar{\sigma}_{d_{1}+1}^*\circ...\circ\bar{\sigma}_{d_{l'}+1}^*;
    \end{align*}
    therefore, since $x^{\alpha+2-l'}_{a_{k-k'+1}}\in\XX(\overline{W}_{\alpha+2-l'}G)$ is $(\bar{\sigma}_{d_{l'}+1}\circ...\circ\bar{\sigma}_{d_{1}+1})$-related to $x^{\alpha+2}_{a_{k-k'+1}}$, we get that 
    \begin{align*}
        (\bar{\sigma}_{d_1}^*\circ...\circ\bar{\sigma}_{d_{l'}}^*\circ\bar{\sigma}_0^*)\Lie_{x^{\alpha+2}_{a_{k-k'+1}}} & = R_{x_{a_{k-k'+1}}}(\bar{\sigma}_{d_{1}+1}^*\circ...\circ\bar{\sigma}_{d_{l'}+1}^*).
    \end{align*}

    \item Type 2: Each of the next $l'$ many terms are of the form $(-1)^{k'+i-1}\Lie_{(y_{b'_i})^{\alpha+2}_{d_{i}+l'+2-i}}$ evaluated at $(h_0^*\circ\delta)^{p-\alpha-1}f\big{(}\mathbbm{v}^{\alpha+1, II(i)}_{(d_1+1)...(d_{l'}+1)}\big{)}$, where $\mathbbm{v}^{\alpha,II(i)}_{d_1...d_{l'}}$ is defined as
    $$\left(\mathbbm{x}^\alpha(a_1,...,a_{k-k'}),(y_{b'_1})^\alpha_{d_1+l'},...,(y_{b'_{i-1}})^{\alpha+1}_{d_{i-1}+l'+2-i},(y_{b'_{i+1}})^{\alpha+1}_{d_{i+1}+l'-i},...,(y_{b'_{l'}})^\alpha_{d_{l'}+1};\mathbbm{y}^\alpha_1(b_1,...,b_{l-l'})\right).$$
    Indeed, since $d_i\leq d_{i'}$ for all $i'\leq i$ and $0\leq d_{i'}<d_{i'+1}$ for each $1\leq i'\leq l'$, the simplicial identities imply that
    \begin{align*}
        \bar{\sigma}_{d_1}^*\circ...\circ\bar{\sigma}_{d_{l'}}^*\circ\bar{\sigma}_0^* & =(\bar{\sigma}_0\circ\bar{\sigma}_{d_i})^*\circ\bar{\sigma}_{d_{1}+2}^*\circ...\circ\bar{\sigma}_{d_{i-1}+2}^*\circ\bar{\sigma}_{d_{i+1}+1}^*\circ...\circ\bar{\sigma}_{d_{l'}+1}^*.
    \end{align*}
    For each $1\leq i'\leq i-1$, $d_i+1\leq d_{i'}+1<d_{i'}+2$, thereby implying $(y_{b'_i})^{\alpha+2-l'}_{d_{i}+1}\in\XX(\overline{W}_{\alpha+2-l'}G)$ is $(\bar{\sigma}_{d_{i-1}+2}\circ...\circ\bar{\sigma}_{d_{1}+2})$-related to $(y_{b'_i})^{\alpha+3-i}_{d_{i}+1}$; in order, as $d_{i'+1}\leq d_{i'}$ for all ranging values of $i'$, the latter is $(\sigma_{d_{l'}+1}\circ...\circ\sigma_{d_{i+1}+1})$-related to $(y_{b'_i})^{\alpha+2}_{d_{i}+l'+2-i}$. Thus, we ultimately conclude that 
    \begin{align*}
       (\bar{\sigma}_{d_1}^*\circ...\circ\bar{\sigma}_{d_{l'}}^*\circ\bar{\sigma}_0^*)\Lie_{(y_{b'_i})^{\alpha+2}_{d_{i}+l'+2-i}} & = R^{(d_i+1)}_{y_{b'_i}}(\bar{\sigma}_{d_{1}+2}^*\circ...\circ\bar{\sigma}_{d_{i-1}+2}^*\circ\bar{\sigma}_{d_{i+1}+1}^*\circ...\circ\bar{\sigma}_{d_{l'}+1}^*).
    \end{align*}

   \item  Type 3: If $l''>0$, there are $l''$ many terms coming from the terms of type IV above. These get written as 
    \begin{align*}
        (-1)^{k'+l'}\sum_{b_{l'+1}=1}^{l''}\bar{\sigma}_0^*(h_0^*\circ\delta)^{p-\alpha-1}f\big{(}\mathbbm{v}^{\alpha+1, III(b_{l-l'+1})}_{(d_1+1)...(d_{l'}+1)}\big{)}, 
    \end{align*}
    where 
    $$\mathbbm{v}^{\alpha,III(b_{l-l'+1})}_{d_1...d_{l'}}=\left(\mathbbm{x}^\alpha(a_1,...,a_{k-k'}),(y_{b'_1})^\alpha_{d_1+l'},...,(y_{b'_{l'}})^\alpha_{d_{l'}+1},(y_{b'_{l'+1}})^{\alpha}_{1};\mathbbm{y}^\alpha_1(b_1,...,b_{l-l'},b_{l-l'+1})\right).$$
    \end{itemize}
    Ultimately, these observations can be synthesised continuing table~\eqref{Eq:Graph2} using the following branching rules for a vertex under-labelled $k',l',l''$: If $k'=l''=0$, there are $l'$ branches; if $k'=0$ and $l''>0$, or $k'>0$ and $l''=0$, there are $l'+1$ branches; and if $k'>0$ and $l''>0$, there are $l'+2$ branches. The larger labels are assigned as follows: 
    \begin{itemize}
        \item  If $k'\neq 0$, the first branch is under-labelled $k'-1,l',l''$ and labelled $$\overrightarrow{P}x(d_1+1)_y...(d_{l'}+1)_y.$$
        \item  For $1\leq i\leq l'$, the $(i+1)$-st branch is under-labelled $k',l'-1,l''$ and labelled $$\overrightarrow{P}y^{d_i+1}(d_1+2)_y...(d_{i-1}+2)_y(d_{i+1}+1)_y...(d_{l'}+1)_y.$$
        \item Lastly, if $l''\neq 0$, the $(l'+2)$-nd branch is under-labelled $k',l'+1,l''-1$ and labelled $$\overrightarrow{P}(d_1)_y...(d_{l'})_y0_y.$$ 
    \end{itemize}
    With these rules at hand, it is fairly straightforward to check that the formula for the van Est map gets recovered. Indeed, after $p$ many steps, all under-labels equal $0,0,0$; $k'=0$ implies there are $k$ many $R_x$ operators, $l'=0$ implies the operator does not have any residual degeneracies and, together with $l''=0$ implies there are $l$ many $R_y$ operators. Because of how the tables are built, it is easy to inspect that every possible permutation of symbols occurs and with the right sign. 
    
    To finish the proof, observe that, using the same type of reasoning behind the tables, the bracket terms can be seen to vanish because, in their branches, there will always be a residual degeneracy.
\end{proof}

\begin{rema}
    In~\cite{Angulo:2022}, the author defines a van Est map between the naturally associated double complexes which one gets by regarding strict Lie 2-groups and strict Lie 2-algebras as groupoids. Given that the complexes therein are weakly isomorphic to the ones presented here, one can try and find an alternative proof for the previous result presenting the van Est map $\Phi$ as a composition of maps of complexes. Indeed, though we were able to prove this is the case in the first degrees, we failed to find a general formula for the weak isomorphisms. We did not further pursue this approach as nothing could be said neither about the algebraic structure nor about the generalization to forms; however, we learned from it the plausibility of Theorem \ref{Theo:CommutesWithDiffs} (cf.~\cite[Theorem 14 and Corollary 15]{Angulo:2022}).
\end{rema}

\subsection{The van Est map at the level of forms}\label{sec:proof}
We are ready to prove the main Theorem \ref{Theo:CommutesWithDiffs}.   Observe that, because the operator $h$ from Proposition~\ref{Prop:FormPerturb} is defined to preserve all degrees and because of the degree of $\delta$, the van Est map for $r$-forms can be computed using a `zig-zag' analogous to that in the proof of Theorem~\ref{Theo:vanEstAtCochains} though on the double complex of $r$-forms on the graded manifolds $\underline{\H_\bullet[2]\oplus\A_\bullet[1]}$ rather than on the Chevalley-Eilenberg double complex, which can be thought of as their functions. The technique thus parallels that of the proof of Theorem~\ref{Theo:vanEstAtCochains}: we keep track of the sums that appear in the differentials by laying out some labeled tables. We use Notation \ref{not:vE}.

\begin{thm}\label{Theo:vanEstAtForms}
    The composition of the map $\iota_\bullet^*\circ (1+\delta\circ h)^{-1}$ from Proposition~\ref{Prop:FormPerturb} and the map $\bar{\partial}_0^*$ from Eq.~\eqref{Eq:partialZeroSuite} yields the van Est map $\Phi$ of Theorem ~\ref{Theo:CommutesWithDiffs}.
\end{thm}

\begin{proof}
     Let $\omega\in\widehat{\Omega}^{p,r}(G)$ then $$\iota_\bullet^*\circ (1+\delta\circ h)^{-1}\bar{\partial}_0^*\omega=\iota_r^*\circ(\delta\circ h_0^*)^p\bar{\partial}_0^*\omega$$ because the $h_i$'s for $i>0$ vanish due to the fact that $\omega$ is normalized. To evaluate the above element we denote by
     \begin{equation*}
        \mathbbm{x}:=(x_1,...,x_k), \quad \mathbbm{y}=(y_1,...,y_l), \quad \mathbbm{z}:=(z_1,...,z_a)\quad \text{and}\quad \mathbbm{w}=(w_1,...,w_b),
     \end{equation*}
     with $a+b=r$ and  $k+2l=p-a-2b$ for $x_i,z_i\in \g_0$, $y_i,w_i\in\h$. We compute $\iota_r^*$ by 
     restricting right-invariant sections to $\bar{\sigma}_0$, thus yielding
    \begin{align}\label{Eq:0thStepForm}
        (\iota_r^*\circ(\delta\circ h_0^*)^p\bar{\partial}_0^*\omega)(\mathbbm{x};\mathbbm{y};\mathbbm{z};\mathbbm{w}) & =\bar{\sigma}_0^*\Big{(}(\delta\circ(h_0^*\circ\delta)^{p-1})\omega(\mathbbm{x}^1;\mathbbm{y}^1;\mathbbm{z}^1;\mathbbm{w}^1)\Big{)}.
    \end{align}
In principle, when computing the Chevalley-Eilenberg differential, there are twelve different types of terms. However, all our vector bundles are trivial and come endowed with the natural connection given by right-invariant sections (see Proposition \ref{prop:A}).  Thus, as in \cite[Example 4.9]{Meinrenken_Li-Bland:2015}, all curvature terms vanish and the differential simplifies to nine types of terms. The first four are, mutatis mutandis, given in \eqref{eq:I-IV}, and the new extra terms are

 \begin{equation*}\def\arraystretch{1.3}
        \begin{array}{rcl}
            \textnormal{V.}& (h_0^*\circ\delta)^{p-1}\omega(\mathbbm{x}^1(\alpha);\mathbbm{y}^1;[x_{\alpha}^1,z_m^1],\mathbbm{z}^1(m);\mathbbm{w}^1) & \text{ for } 1\leq\alpha\leq k\text{ and } 1\leq m\leq a\\
           \textnormal{VI.} &(h_0^*\circ\delta)^{p-1}\omega(\mathbbm{x}^1(\alpha);\mathbbm{y}^1;\mathbbm{z}^1;[x_{\alpha}^1,w_n^1],\mathbbm{w}^1(n))&\text{ for } 1\leq\alpha\leq k \text{ and } 1\leq n\leq b\\
             \textnormal{VII.} & (h_0^*\circ\delta)^{p-1}\omega(\mathbbm{x}^1;\mathbbm{y}^1;\mathbbm{z}^1(m);\mathbbm{w}^1\vert\rho(z_m^1)) 
             &\text{ for }1\leq m\leq a\\
              \textnormal{VIII.} & (h_0^*\circ\delta)^{p-1}\omega(\mathbbm{x}^1;\mathbbm{y}^1;\mathbbm{z}^1,w_n^1;\mathbbm{w}^1(n))&\text{ for } 1\leq n\leq b\\
               \textnormal{IX.} & (h_0^*\circ\delta)^{p-1}\omega(\mathbbm{x}^1;\mathbbm{y}^1(\beta);\mathbbm{z}^1(m);[y_{\beta}^1,z_m^1],\mathbbm{w}^1)&\text{ for }1\leq\beta\leq l \text{ and } 1\leq m\leq a.
        \end{array}
    \end{equation*}
 We refer again to the terms of type II, III, V, VI, and IX as {\it bracket terms}; and, as in the case of cochains, these ultimately vanish due to $\omega$ being normalized. Focusing on the remaining terms, note that composing with the degeneracy of Eq.~\eqref{Eq:0thStepForm}, the terms of type I, IV, VII, and VIII turn respectively into:
\begin{equation*}
    \def\arraystretch{1.3}
    \begin{array}{rcl}
        1. & R_{x_\alpha}(h_0^*\circ\delta)^{p-1}\omega(\mathbbm{x}^1(\alpha);\mathbbm{y}^1;\mathbbm{z}^1;\mathbbm{w}^1) & \text{ for } 1\leq\alpha\leq k\\
        2. & \bar{\sigma}_0^*(h_0^*\circ\delta)^{p-1}\omega(\mathbbm{x}^1,y_\beta^1;\mathbbm{y}^1(\beta);\mathbbm{z}^1;\mathbbm{w}^1) & \text{ for }1\leq\beta\leq l\\
        3. & J_{z_m}(h_0^*\circ\delta)^{p-1}\omega(\mathbbm{x}^1;\mathbbm{y}^1;\mathbbm{z}^1(m);\mathbbm{w}^1) & \text{ for }1\leq m\leq a\\
        4. & \bar{\sigma}_0^*(h_0^*\circ\delta)^{p-1}\omega(\mathbbm{x}^1;\mathbbm{y}^1;\mathbbm{z}^1,w_n^1;\mathbbm{w}^1(n)) & \text{ for } 1\leq n\leq b.
    \end{array}
    \end{equation*}
We collect the sum of these terms and record them using the following table
 \begin{eqnarray}\label{Eq:Graph1Form}
    \begin{array}{|c|c|}
    \hline
       (k \ 0 \ l \ a \ 0 \ b )  & \begin{array}{c}
            \underset{k-1,0,l,a,0, b}{x}  \\
            \hline
            \underset{k,1,l-1, a,0,b}{0_y}\\
            \hline 
            \underset{k,0,l,a-1,0, b}{z}\\
             \hline 
            \underset{k,0,l,a,1, b-1}{0_w}
       \end{array} \\
         \hline
    \end{array}
     \end{eqnarray}
 The under-label $k',l',l'',a',b',b''$ indicates that the argument takes values in $$\Gamma\left(\bigwedge\nolimits^{r-(a'+b'+b'')}T(W_{\bullet}G)\otimes\bigwedge\nolimits^{k'+l'}\A_\bullet\otimes \sym^{l''}\H_\bullet\otimes \sym^{a'+b'}\A_\bullet\otimes\bigwedge\nolimits^{b''}\H_\bullet\right);$$ it further indicates that, in the antisymmetric part of $\A_\bullet$, the first $k'$ many arguments come from sections of $\A_\bullet$ and the subsequent $l'$ many come from sections of $\H_\bullet$ regarded via the inclusion $j_\bullet$ as sections of $\A_\bullet$, and, similarly, in the symmetric part of $\A_\bullet$, the first $a'$ many arguments come from sections of $\A_\bullet$ and the subsequent $b'$ many come from sections of $\H_\bullet$. As for the larger labels, $x$ represents the sum $\sum_{\alpha=1}^k(-1)^{\alpha-1}R_{x_\alpha}$ of terms of type 1; the larger label $0_y$ represents the sum $(-1)^k\sum_{\beta=1}^l\bar{\sigma}_0^*$ of terms of type 2; the larger label $z$ represents the sum $(-1)^k\sum_{m=1}^aJ_{z_m}$ of terms of type 3; and, lastly, the larger label $0_y$ represents the sum $\sum_{n=1}^b(-1)^{k+n-1}\bar{\sigma}_0^*$ of terms of type 4. 

We compute $h_0^*$ using right-invariant sections and evaluating at degeneracies, thus turning terms of type 1 through 4 respectively into:

\begin{equation*}
    \def\arraystretch{1.3}
    \begin{array}{rcl}
        1'. & R_{x_\alpha}\bar{\sigma}_0^*\Big{(}(\delta\circ(h_0^*\circ\delta)^{p-2})\omega(\mathbbm{x}^2(\alpha);\mathbbm{y}^2_1;\mathbbm{z}^2;\mathbbm{w}^2_1)\Big{)} & \text{ for } 1\leq\alpha\leq k\\
        2'. & (\bar{\sigma}_0^*)^2\Big{(}(\delta\circ(h_0^*\circ\delta)^{p-2})\omega(\mathbbm{x}^2,(y_\beta)^2_2;\mathbbm{y}^2_1(\beta);\mathbbm{z}^2;\mathbbm{w}^2_1)\Big{)} & \text{ for }1\leq\beta\leq l\\
        3'. & J_{z_m}\Big{(}(\delta\circ(h_0^*\circ\delta)^{p-2})\omega(\mathbbm{x}^2;\mathbbm{y}^2_1;\mathbbm{z}^2(m);\mathbbm{w}^2_1)\Big{)} & \text{ for }1\leq m\leq a\\
        4'. & (\bar{\sigma}_0^*)^2\Big{(}(\delta\circ(h_0^*\circ\delta)^{p-2})\omega(\mathbbm{x}^2;\mathbbm{y}^2_1;\mathbbm{z}^2,(w_n)^2_2;\mathbbm{w}^2_1(n))\Big{)} & \text{ for } 1\leq n\leq b.
    \end{array}
    \end{equation*}
Computing the subsequent Chevalley-Eilenberg type differential and neglecting the bracket terms yields four new terms for each type. For type 1', we get the two terms corresponding to \eqref{eq:1.i} plus the new terms
\begin{equation*}
    \def\arraystretch{1.3}
    \begin{array}{rcl}
        1.3 & R_{x_{\alpha}}J_{z_m}(h_0^*\circ\delta)^{p-2}\omega(\mathbbm{x}^2(\alpha);\mathbbm{y}^2_1;\mathbbm{z}^2(m);\mathbbm{w}^2_1) & \text{ for }1\leq\alpha\leq k\text{ and } 1\leq m\leq a\\
        1.4 & R_{x_\alpha}\bar{\sigma}_0^*(h_0^*\circ\delta)^{p-2}\omega(\mathbbm{x}^2(\alpha);\mathbbm{y}^2_1;\mathbbm{z}^2,(w_n)^2_1;\mathbbm{w}^2_1(n)) & \text{ for } 1\leq\alpha\leq k\text{ and } 1\leq n\leq b,
    \end{array}
    \end{equation*}
    while for type 2', we get the terms corresponding to \eqref{eq:2.i} plus the extra terms
    \begin{equation*}
    \def\arraystretch{1.3}
    \begin{array}{rcl}
        2.4 & J_{z_m}\bar{\sigma}_1^*(h_0^*\circ\delta)^{p-2}\omega(\mathbbm{x}^2,(y_{\beta})^2_2;\mathbbm{y}^2_1(\beta);\mathbbm{z}^2(m);\mathbbm{w}^2_1) & \text{ for } 1\leq\beta\leq l \text{ and } 1\leq m\leq a \\
        2.5& (\bar{\sigma}_0^*)^2(h_0^*\circ\delta)^{p-2}\omega(\mathbbm{x}^2,(y_\beta)^2_2;\mathbbm{y}^2_1(\beta);\mathbbm{z}^2,(w_n)^2_1;\mathbbm{w}^2_1(n)) & \text{ for } 1\leq\beta\leq l \text{ and } 1\leq n\leq b.
    \end{array}
    \end{equation*}
     The types 3' and 4' are new and give respectively the following terms:
    \begin{equation*}
    \def\arraystretch{1.3}
    \begin{array}{rcl}
        3.1 & J_{z_m}R_{x_{\alpha}}(h_0^*\circ\delta)^{p-2}(h_0^*\circ\delta)^{p-2}\omega(\mathbbm{x}^2(\alpha);\mathbbm{y}^2_1;\mathbbm{z}^2(m);\mathbbm{w}^2_1) & \text{ for } 1\leq\alpha\leq k \text{ and } 1\leq m\leq a\\
        3.2 & J_{z_m}\bar{\sigma}_0^*(h_0^*\circ\delta)^{p-2}\omega(\mathbbm{x}^2,(y_\beta)^2_1;\mathbbm{y}^2_1(\beta);\mathbbm{z}^2(m);\mathbbm{w}^2_1) & \text{ for }1\leq\beta\leq l \text{ and } 1\leq m\leq a\\
        3.3 & J_{z_{m_1}}J_{z_{m_2}}(h_0^*\circ\delta)^{p-2}\omega(\mathbbm{x}^2;\mathbbm{y}^2_1;\mathbbm{z}^2(m_1,m_2);\mathbbm{w}^2_1) & \text{ for }1\leq m_1,m_2\leq a\\
        3.4 & J_{z_m}\bar{\sigma}_0^*(h_0^*\circ\delta)^{p-2}\omega(\mathbbm{x}^2;\mathbbm{y}^2_1;\mathbbm{z}^2(m),(w_n)^2_1;\mathbbm{w}^2_1(n)) & \text{ for } 1\leq m\leq a \text{ and } 1\leq n\leq b,\\
    \end{array}
    \end{equation*}
    \begin{equation*}
    \def\arraystretch{1.3}
    \begin{array}{rcl}
        4.1 & R_{x_\alpha}\bar{\sigma}_1^*(h_0^*\circ\delta)^{p-2}\omega(\mathbbm{x}^2(\alpha);\mathbbm{y}^2_1;\mathbbm{z}^2,(w_n)^2_2;\mathbbm{w}^2_1(n)) & \text{ for } 1\leq\alpha\leq k \text{ and } 1\leq n\leq b\\
        4.2 & (\bar{\sigma}_0^*)^2(h_0^*\circ\delta)^{p-2}\omega(\mathbbm{x}^2,(y_\beta)^2_1;\mathbbm{y}^2_1(\beta);\mathbbm{z}^2,(w_n)^2_2;\mathbbm{w}^2_1(n)) & \text{ for }1\leq\beta\leq l\text{ and }1\leq n\leq b\\
        4.3 & J_{z_m}\bar{\sigma}_1^*(h_0^*\circ\delta)^{p-2}\omega(\mathbbm{x}^2;\mathbbm{y}^2_1;\mathbbm{z}^2(m),(w_n)^2_2;\mathbbm{w}^2_1(n)) & \text{ for }1\leq m\leq a  \text{ and } 1\leq n\leq b\\
        4.4 & J_{w_n}(h_0^*\circ\delta)^{p-2}\omega(\mathbbm{x}^2;\mathbbm{y}^2_1;\mathbbm{z}^2;\mathbbm{w}^2_1(n)) & \text{ for } 1\leq n\leq b\\
        4.5& (\bar{\sigma}_0^*)^2(h_0^*\circ\delta)^{p-2}\omega(\mathbbm{x}^2;\mathbbm{y}^2_1;\mathbbm{z}^2,(w_{n_1})^2_2,(w_{n_2})^2_1;\mathbbm{w}^2_1(n_1,n_2)) & \text{ for } 1\leq n_1,n_2\leq l.
    \end{array}
    \end{equation*}
Recording these terms in the table~\eqref{Eq:Graph1Form} induces the branching 
 \begin{eqnarray*}\label{Eq:Graph2Form}
    \begin{array}{|c|c|}
    \hline
      (k\ 0 \ l\ a\ 0 \ b)   & \begin{array}{c|c}
         \underset{k-1,0,l,a, 0, b}{x}  & \begin{array}{c}
              \underset{k-2,0,l, a, 0, b}{xx}  \\
              \hline 
             \underset{k-1,1,l-1, a, 0, b}{x0_y} \\
             \hline 
             \underset{k-1,0,l, a-1, 0, b}{xz}\\
              \hline 
             \underset{k-1,0,l, a, 1, b-1}{x0_w}
         \end{array} \\
         \hline 
         \underset{k,1,l-1, a, 0, b}{0_y}  & \begin{array}{c}
               \underset{k-1,1,l-1, a, 0, b}{x1_y} \\
               \hline
             \underset{k,0,l-1, a,0,b}{y^1}\\
             \hline
             \underset{k,2,l-2, a,0,b}{0_y0_y}\\
             \hline
             \underset{k,1,l-1, a-1,0,b}{z1_y}\\
             \hline
             \underset{k,1,l-1, a,1,b-1}{0_y0_w}
         \end{array}
      \end{array} \\
      \hline
    \end{array}\quad     \begin{array}{|c|c|}
    \hline
      (k\ 0 \ l\ a\ 0 \ b)   & \begin{array}{c|c}
         \underset{k,0,l,a-1, 0, b}{z}  & \begin{array}{c}
              \underset{k-1,0,l, a-1, 0, b}{zx}  \\
              \hline 
             \underset{k,1,l-1, a-1, 0, b}{z0_y} \\
             \hline 
             \underset{k,0,l, a-2, 0, b}{zz}\\
              \hline 
             \underset{k,0,l, a-1, 1, b-1}{z0_w}
         \end{array} \\
         \hline 
         \underset{k,0,l, a, 1, b-1}{0_w}  & \begin{array}{c}
               \underset{k-1,0,l, a, 1, b-1}{x1_w} \\
               \hline
             \underset{k,1,l-1, a,1,b-1}{0_w0_y}\\
             \hline
             \underset{k,0,l, a-1,1,b-1}{z1_w}\\
             \hline
             \underset{k,1,l, a,0,b-1}{w^1}\\
             \hline
             \underset{k,0,l, a,2,b-2}{0_w0_w}
         \end{array}
      \end{array} \\
      \hline
    \end{array}
    \end{eqnarray*}
As in the proof of Theorem~\ref{Theo:vanEstAtCochains}, the labels $xx$, $x0_y$, $x1_y$, $y^1$ and $0_y0_y$ represents the operators described in \eqref{eq:xx}--\eqref{eq:00} summing over terms of types 1.1, 1.2, 2.1, 2.2 and 2.3. In a similar manner, the labels $xz$ and $zx$ stand respectively for 
$$
\sum_{\alpha=1}^k(-1)^{\alpha+k}R_{x_{\alpha}}\sum_{m=1}^{a}J_{z_m}\quad \text{ and }\quad (-1)^k\sum_{m=1}^{a}J_{z_m}\sum_{\alpha=1}^k(-1)^{\alpha-1}R_{x_{\alpha}}
$$ 
summing over terms of types 1.3 and 3.1 respectively. The labels $z1_y$ and $z0_y$ represent 
$$
\sum_{m=1}^a(-1)^dJ_{z_m}\sum_{\beta=1}^l \bar{\sigma}_d^*
$$ 
summing over terms of types 2.4 for $d=1$ and 3.2 for $d=0$. The label $zz$ represents $$
\sum_{m_1=1}^aJ_{z_{m_1}}\sum_{m_2=1,m_2\neq m_1}^{a}J_{z_{m_2}}
$$ 
summing over terms of type 3.3 and the label $w^1$ stands for the operator $\sum_{n=1}^b(-1)^{n-1}J_{w_n}$ summing over terms of type 4.4. The labels $x0_w$ and $x1_w$ represent  
$$
\sum_{\alpha=1}^k(-1)^{(\alpha-1)+(k-1)}R_{x_\alpha}\sum_{n=1}^b(-1)^{n-1}\bar{\sigma}_d^*
$$ 
summing over terms of types 1.4 for $d=0$ and 4.1 for $d=1$; and, analogously, the labels $z0_w$ and $z1_w$ represent 
$$
\sum_{m=1}^aJ_{z_m}\sum_{n=1}^b(-1)^{n-1}\bar{\sigma}_d^*
$$ 
summing over terms of types 3.4 for $d=0$ and 4.3 for $d=1$. The labels $0_y0_w$ and $0_w0_y$ stand respectively for 
$$
(-1)^{k+(k+1)}\sum_{\beta=1}^l\sum_{n=1}^{b}(-1)^{n-1}(\bar{\sigma}_0^*)^2\quad\text{and}\quad \sum_{n=1}^{b}\sum_{\beta=1}^l(-1)^{n-1}(\bar{\sigma}_0^*)^2
$$ summing over terms of types 2.5. and 4.2. Lastly,  the label $0_w0_w$ represents the composition 
$$
\sum_{n_1=1}^l(-1)^{n_1-1}\left(\sum_{n_2=1}^{n_1-1}(-1)^{n_2-1}(\bar{\sigma}_0^*)^2+\sum_{n_2=n_1+1}^{b}(-1)^{n_2}(\bar{\sigma}_0^*)^2\right)
$$ 
summing over terms of type 4.5. 

 In general, $\eta:=((k-k')+2(l-l'')-l'+(a-a')+2(b-b'')-b')$ many steps are required to reach a vertex under-labelled $k',l',l'',a',b',b''$, with $0\leq k'\leq k$, $0\leq l'+l''\leq l$, $0\leq a'\leq a$ and $0\leq b'+b''\leq b$. This process corresponds to computing a type of term in the $\eta$-th iteration of $h_0^*\circ\delta$. 

Let the larger label on $k',l',l'',a',b',b''$ be $\overrightarrow{P}(d_1)_{v_1}...(d_{l'+b'})_{v_{l'+b'}}$, where 
$$
\overrightarrow{P}:=\varrho(x...xy^{\mu_1}...y^{\mu_{l-(l'+l'')}}z...zw^{\nu_1}...w^{\nu_{b-(b'+b'')}})\quad\text{for } \quad \varrho\in S_{\eta-l+l''-b+b''},
$$ 
with $v_i\in\lbrace y,w\rbrace$ for all $1\leq i\leq l'+b'$, the cardinality of the set $\mathcal{I}_y=\lbrace i\vert v_i=y\rbrace$ is $l'$, and $(k-k')+(a+a')\geq d_1\geq ...\geq d_{l'+b'}\geq 0$. This label represents a sum over $1\leq\beta_1,...,\beta_{l'}\leq l$, a sum over $1\leq m_1,...,m_{a-a'}\leq l$ and a alternating sums over $1\leq\alpha_1,...,\alpha_{k-k'}\leq k$ and $1\leq n_1,...,n_{b'}\leq k$ of the operators 
$$\varrho\left(R_{x_{\alpha_1}}...R_{x_{\alpha_{k-k'}}}R^{(\mu_1)}_{y_{\beta_1}}...R^{(\mu_{l-(l'+l'')})}_{y_{\beta_{l-(l'+l'')}}}J_{z_{m_1}}...J_{z_{m_{a-a'}}}J^{(\nu_1)}_{w_{n_1}}...J^{(\nu_{b-(b'+b'')})}_{w_{n_{b-(b'+b'')}}}\right)(\bar{\sigma}_{d_1}^*\circ...\circ\bar{\sigma}_{d_{l'+b'}}^*),$$ 
where, as before, $R^{(\mu)}_{y}:=(\sigma_0\circ\sigma_{\mu-1})^*\Lie_{y^\bullet_\mu}$ and $J^{(\nu)}_{w}:=(\sigma_0\circ\sigma_{\nu-1})^*\iota_{w^\bullet_\nu}$ (which, again, indicates there are permutation depending restrictions on the values $\mu_j$ and $\nu_j$). 

To express at what these operators are evaluated, consider the ordered sets $\mathcal{I}_y$ and $\mathcal{I}_w=\lbrace i\vert v_i=w\rbrace$, and define $\mathcal{I}_y+\mathcal{I}_w$ as the union of the two sets ordered by declaring that every element in $\mathcal{I}_y$ be smaller than every element in $\mathcal{I}_w$. There exists a unique permutation $\varsigma\in S_{l'+b'}$ enumerating this ordering, then, defining $\beta'_i:=\beta_{l-(l'+l'')+i}$, $n'_i:=n_{b-(b'+b'')+i}$ and $\chi(i):=d_{\varsigma(i)}+l'+b'-1+\varsigma(i)$, the operators are evaluated at $(h_0^*\circ\delta)^{p-\eta}\omega(\mathbbm{v}^\eta)$, where 
    \begin{align*}
        \mathbbm{v}^{\eta}&=\left(  \mathbbm{x}^\eta(\alpha_1,...,\alpha_{k-k'}),(y_{\beta'_1})^\eta_{\chi(1)},(y_{\beta'_2})^\eta_{\chi(2)},..., (y_{\beta'_{l'-1}})^\eta_{\chi(l'-1)},(y_{\beta'_{l'}})^\eta_{\chi(l')};\mathbbm{y}^\eta_1(\beta_1,...,\beta_{l-l'}); \right.\\
            & \left. \mathbbm{z}^\eta(m_1,...,m_{a-a'}),(w_{n'_1})^\eta_{\chi(l'+1)},(w_{n'_2})^\eta_{\chi(l'+2)},...,(w_{n'_{b'-1}})^\eta_{\chi(l'+b'-1)},(w_{n'_{b'}})^\eta_{\chi(l'+b')};\mathbbm{w}^\eta_1(n_1,...,n_{b-b'})\right).
    \end{align*}

Inductively, the homotopy is computed to be the same total sum ranging over  $(h_0^*\circ\delta)^{p-\eta-1}\omega(\mathbbm{v}^{\eta+1})$, where, in the same fashion as Eq.~\eqref{Eq:oneStepHomotopy}, each of the subindices $\chi(i)$ also gets replaced by $\chi(i)+1$. Disregarding the bracket terms one last time, six types of terms are coming out of computing the Chevalley-Eilenberg differential. Using the same reasoning of the proof of Theorem~\ref{Theo:vanEstAtCochains}, it follows from the simplicial identities that these six types of terms correspond to the following branching rules for a vertex under-labeled $k',l',l'',a',b',b''$: For each non-zero $k'$, $l''$, $a'$ and $b''$ there is a single branch, and, additionally, there are $l'+b'$ branches. The larger labels are assigned as follows: 
\begin{itemize}
    \item If $k'\neq 0$, the first branch is under-labelled $k'-1,l',l'',a',b',b''$ and labelled $$\overrightarrow{P}x(d_1+1)_{v_1}...(d_{l'+b'}+1)_{v_{l'+b'}}.$$
    \item For $1\leq i\leq l'$, the $(i+1)$-st branch is under-labelled $k',l'-1,l'',a',b',b''$ and labelled $$\overrightarrow{P}y^{d_{\varsigma(i)}+1}(d_1+2)_{v_1}...(d_{\varsigma(i)-1}+2)_{v_{\varsigma(i)-1}}(d_{\varsigma(i)+1}+1)_{v_{\varsigma(i)+1}}...(d_{l'+b'}+1)_{v_{l'+b'}}.$$ 
    \item If $l''\neq 0$, the $(l'+2)$-nd branch is under-labelled $k',l'+1,l''-1,a',b',b''$ and labelled $$\overrightarrow{P}(d_1+1)_{v_1}...(d_{l'+b'}+1)_{v_{l'+b'}}0_y.$$
    \item If $a'\neq 0$, the $(l'+3)$-rd branch is under-labelled $k',l',l'',a'-1,b',b''$ and labelled $$\overrightarrow{P}z(d_1+1)_{v_1}...(d_{l'+b'}+1)_{v_{l'+b'}}.$$
    \item For $1\leq i\leq b'$, the $(l'+i+3)$-rd branch is under-labelled $k',l',l'',a',b'-1,b''$ and labelled $$\overrightarrow{P}w^{d_{\varsigma(l'+i)}+1}(d_1+2)_{v_1}...(d_{\varsigma(l'+i)-1}+2)_{v_{\varsigma(l'+i)-1}}(d_{\varsigma(l'+i)+1}+1)_{v_{\varsigma(l'+i)+1}}...(d_{l'+b'}+1)_{v_{l'+b'}}.$$
    \item Lastly, if $b''\neq 0$, the $(l'+b'+4)$-th branch is under-labelled $k',l',l'',a,b'+1,b''-1$ and labelled $$\overrightarrow{P}(d_1+1)_{v_1}...(d_{l'+b'}+1)_{v_{l'+b'}}0_w.$$  
\end{itemize}

    The combinatorics of these labeled tables give back the formula for the van Est map $\Phi$, as desired. After $p$ many steps, all under-labels equal $0,0,0,0,0,0$; $k'=0$ implies there are $k$ many $R_x$ operators, $a'=0$ implies there are $a$ many $J_z$ operators, $l'=b'=0$ imply the operator does not have any residual degeneracies and, together with $l''=0$ and $b''=0$, they imply there are $l$ many $R_y$ operators and $b$ many $J_w$ operators. Every possible permutation of symbols appears and they carry a sign that coincides with the one in the formula. 
\end{proof}

\subsection{The van Est isomorphism}\label{subsec-vanEstIso}

As the van Est map can be viewed as a differentiation procedure, in some situations, it is possible to obtain an integration procedure in the opposite direction. To obtain a cochain map in the other direction, we need a homotopy operator with respect to the vertical differential $\delta$. 

In the setup, we already have the map $\bar{\partial}_0^*$ dual to the d\'ecalage projection. Observe that $\bar{\sigma}_0$ defines a natural section thereof, and thus its dual defines a map
\begin{eqnarray}\label{Eq:splitDec}
    \bar{\sigma}_0^*:\cW^{0,r}(\H_p\rightarrow\A_p)\longrightarrow\Omega^r(\overline{W}_p\G),
\end{eqnarray}
left inverse to $\bar{\partial}_0^*$. The degeneracy $\bar{\sigma}_0^*$ is a cochain map with respect to $\delta$ and $d$; however, since $\bar{\sigma}_0$ is not a simplicial map, $\bar{\sigma}_0^*$ is not a cochain map with respect to $\partial$.

Under the assumption that $G_0$ are $k$-connected and $G_1$ is source $(k-1)$-connected, one concludes that the fibers of each of the principal bundles in the d\'ecalage fibration are at least $(k-1)$-connected. Furthermore, using homological perturbation theory, given a homotopy operator, one can construct the inverse in degrees less than or equal to $k$ at the level of cochains.  The following result explains how this connectedness hypothesis is related to the cohomology of the Lie 2-algebroid $\H_p\rightarrow\A_p$.

\begin{lem}\label{lemma:vanCoh}
    If $G_p$ is $(k-1)$-connected then  $$H^0(\H_p\rightarrow\A_p)=C^\infty(\overline{W}_p\G)\quad \text{and}\quad H^n(\H_p\rightarrow\A_p)=0$$
    for all $1\leq n< k$.
\end{lem}
\begin{proof}
    The Chevalley-Eilenberg complex of $\H_p\rightarrow\A_p$ can be regarded as the total complex of the double complex $$E_0^{m,n}=\Gamma\left(\bigwedge\nolimits^{m-n}\A_p^*\otimes\sym^n\H_p^*\right)$$ whose vertical differentials are dual to the 1-bracket $j_p$ and whose horizontal differentials are dual to the 2-bracket. Let $(E_r^{m,n},d_r)$ be the spectral sequence of the filtration by columns of $E_0^{m,n}$. We prove in the sequel that $(E_1^{m,n},d_1)$ is isomorphic to the Chevalley-Eilenberg complex of the foliation $\F_p$; thus, the result follows from \cite[Theorem 2]{Crainic:2003}.

    To prove that the cohomology of  
    \begin{eqnarray}\label{Eq:vertCoh}
         \cdots\xrightarrow{}\Gamma\left(\bigwedge\nolimits^{m-n}\A_p^*\otimes\sym^n\H^*_p\right)\xrightarrow{j_p^*}\Gamma\left(\bigwedge\nolimits^{m-(n+1)}\A_p^*\otimes\sym^{n+1}\H^*_p\right)\xrightarrow{}\cdots
    \end{eqnarray}
    is concentrated in degree zero amounts to a `skeletal replacement'. Indeed, recall from Proposition~\ref{prop:A} that right-multiplication gives trivializations 
    $$\A_p\cong W_pG\times\gg_{p+1}\quad\text{and}\quad\H_p\cong W_p\G\times\hh.$$ Therefore,  since $j_p$ is injective, the only non-vanishing row is $n=0$, where $C^\infty(W_pG)\otimes\bigwedge^{m}\gg_p^*\cong\Gamma(\bigwedge^{m}\cF_p)$, as claimed.
\end{proof}

Lemma~\ref{lemma:vanCoh} implies that the double complex~\eqref{Eq:AuxDoubleCx} has vanishing $\delta$-cohomology in bidegree $(p,q)$ for all $q\leq k-1$ and also in bidegree $(0,k)$. This suffices to prove the following result:

\begin{thm}\label{Theo:IndIsosCochains}
Let $G_\bullet$ be a strict Lie $2$-group and assume $G_0$ is $k$-connected and $G_1$ is source $(k-1)$-connected. The van Est map $\Phi$ restricted to cochains induces isomorphisms in cohomology for all $p\leq k$ and it is injective in cohomology for $p=k+1$.
\end{thm}
\begin{proof}
    Let $(E_r^{p,q},d_r)$ be the spectral sequence of the filtration by rows of the double complex \eqref{Eq:AuxDoubleCx}. Then, Lemma~\ref{lemma:vanCoh} implies $E_1^{p,q}\equiv 0$ for $0<q<k-1$ and for $(p,q)=(0,k)$ and, otherwise, $E_1^{p,0}\cong C^\infty(\overline{W}_{p} G)$. Thus, the map~\eqref{Eq:splitDec} induces an isomorphism between $E_2^{p,0}$ and $H^p(\overline{W}_\bullet G)$ with inverse~\eqref{Eq:partialZero} for $p\leq k$ and it is injective in cohomology for $p=k+1$. Composing with the map induced by $\pi_0^*$ in cohomology, we get the desired isomorphism.
\end{proof}

\begin{rema}\label{Rmk:invFormula}
    An explicit formula for the inverse to the van Est map is sometimes available at the level of cochains, see e.g. \cite{mein:VEint}. Letting $\tau_{\leq k}C^{p,q}$ be the truncated complex, if there is a homotopy operator $\eta$ for it, with $\delta\circ\eta+\eta\circ\delta=\textnormal{Id}-\sigma_0^*\circ \delta_0^*$; 
    by the Perturbation Lemma, the composition  $\sigma_0^*\circ(1+\partial\circ\eta)^{-1}$ is a cochain map for the total differential, and its composition with $\pi_0^*$ gives the desired inverse to $\Phi$.
\end{rema}
\begin{rema}\label{Rmk:OnHyp}
    The proof of Theorem~\ref{Theo:IndIsosCochains} suggests that, presumably, the hypotheses can be relaxed. 
    The Morita invariance of cohomology in \cite{behcohsta} implies that the cohomologies of both Lie 2-groups and Lie 2-algebras are isomorphic to the cohomologies of their skeletal replacements. Therefore, whenever the orbit space of $G_\bullet$ is smooth, the conclusion of Theorem~\ref{Theo:IndIsosCochains} still holds if one asks that the isotropies and the orbit space of $G_\bullet$ are, respectively, $(k-1)$-connected and $k$-connected. This can thus be interpreted as a connectedness hypothesis on the stack presented by $G_\bullet$.
\end{rema}

\begin{rema}
    Theorem \ref{Theo:IndIsosCochains} can be extended to forms. This follows immediately from the fact that the proof of Theorem~\ref{Theo:vanEstAtForms} at level $r$-forms only relies on the double complex of $r$-forms (see Subsection~\ref{ssec-WeilAlg}). That there are isomorphisms at the level of $r$-forms follows similarly by noticing that the $\delta$-cohomology of the complexes of $r$-forms can be computed to vanish. This is analogous to the method used in \cite[Lemma 5.3]{abad:VE}.
\end{rema}

\section{Applications}\label{sec-Applications}

In this section, we use the van Est map to differentiate shifted symplectic structures on strict Lie $2$-groups in the sense of \cite{Lesdiablerets}.

\subsection{The semi-direct product with the coadjoint representation}\label{sec:app2}

Let $G$ be a Lie group and consider the strict Lie $2$-group $\cG$ associated with the crossed module $(G,\gg^*,1,Ad^*)$ as per Remark~\ref{Rmk:CrossMods1}. 
By Proposition \ref{prop:lie} and Remark~\ref{Rmk:CrossMods2}, the corresponding infinitesimal crossed module $(\gg,\gg^*,0,ad^*)$ gives rise to a strict Lie $2$-algebra whose associated degree $2$ NQ-manifold is usually denoted by $\g[1]\oplus\g^*[2]=T^*[3]\g[1]$  when equipped with the degree $1$ vector field $\cL_{\delta_{CE}}.$ 

$(T^*[3]\g[1], \cL_{\delta_{CE}})$ is a degree $3$ NQ-manifold that carries a canonical symplectic form
$$\omega^{\inf}=-d\theta^{\inf}\in\Omega^{2,3}(T^*[3]\g[1]),$$
where the tautological 1-form $\theta^{\inf}\in\Omega^{1,3}(T^*[3]\g[1])$ is defined by the property $\theta^{\inf}(\widetilde{X})=X$ for any $X\in\fX(\g[1])$ with $\widetilde{X}$ denoting the cotangent lift. When picking coordinates $|\alpha^i|=1$ and $|\beta_i=\frac{\partial}{\partial\alpha^i}|=2$, then $\theta^{\inf}=\beta_id\alpha^i$.

In what follows, we construct an integration of the symplectic structure $\omega^{\inf}$ in the Lie 2-group $\mathbb{K}_\bullet=\overline{W}\cG.$ Denoting the right Maurer-Cartan $1$-forms on $G$ by $\theta^r\in\Omega^1(G;\g)$, the tautological 1-form and the canonical symplectic form on $T^*G=\g^*\times G$ can be written as
$$\theta=\langle \pr_1, \pr_2^*\theta^r\rangle\in\Omega^1(\g^*\times G)\quad \omega=-d\theta\in\Omega^2(\g^*\times G).$$
\begin{prop}
    For a Lie group $G$, if $\cG$ is the strict Lie 2-group associated with the crossed module $(G,\gg^*,1,Ad^*)$, the Lie $2$-group $\mathbb{K}_\bullet=\overline{W}\cG$ carries a $3$-shifted symplectic structure given by  $$p_2^*\omega\in\Omega^2(\mathbb{K}_3),$$ where 
    $p_{2}:\mathbb{K}_3\to \g^*\times G,$ $p_{2}(\xi_1, \xi_2,g_2; \xi,g_1; g_0)=(\xi_2,g_1).$ Moreover, the image of $p_2^*\omega$ under the van Est map defined in Theorem \ref{Theo:CommutesWithDiffs} coincides with the canonical symplectic form on $T^*[3]\gg[1]$, i.e.
    $$\Phi(p_2^*\omega)=\omega^{\inf}.$$
\end{prop}  
\begin{proof}
   Using \eqref{Eq:BFaceMaps}, we see that the face and degeneracy maps of $\mathbb{K}_\bullet$ are given by
   \begin{align*}\label{Eq:ExFaceMaps}
\bar{\sigma}_i(g) & =\begin{cases}
(0,1;g)& i=0 \\
(0,g;1)& i=1 ,\end{cases}\quad\bar{\partial}_i(\Vec{A}) =\begin{cases}
g_0& i=0 \\
g_1g_0& i=1 \\
g_1& i=2,\end{cases}\quad\bar{\sigma}_i(\Vec{A}) =\begin{cases}
(0,0,1;\xi,g_1;g_0)& i=0 \\
(0,\xi,g_1;1,0;g_0)& i=1 \\
(\xi,0,g_1;0,g_0;1)& i=2,\end{cases} \\
\bar{\partial}_i(\Vec{B}) & =\begin{cases}
(\xi,g_1;g_0)& i=0 \\
(\xi+Ad_{g_1}^*\xi_1,g_2g_1;g_0)& i=1 \\
(\xi_1+\xi_2,g_2;g_1g_0)& i=2 \\
(\xi_1,g_2;g_1) & i=3,\end{cases}\quad 
\bar{\partial}_i(\Vec{C}) =\begin{cases}
(\xi_1,\xi_2,g_2;\xi,g_1;g_0)& i=0 \\
(\xi_1+Ad_{g_2}^*\xi_2',\xi_2+Ad_{g_2}^*\xi_3',g_3g_2;\xi,g_1;g_0)& i=1 \\
(\xi_1'+\xi_2',\xi_3',g_3;\xi+Ad_{g_1}^*\xi_1,g_2g_1;g_0)& i=2 \\
(\xi_1',\xi_2'+\xi_3',g_3;\xi_1+\xi_2,g_2;g_1g_0) & i=3 \\
(\xi_1',\xi_2',g_3;\xi_1,g_2;g_1) & i=4,\end{cases}
\end{align*}   
for $g\in\mathbb{K}_1$, $\Vec{A}=(\xi,g_1; g_0)\in\mathbb{K}_2$, $\Vec{B}=(\xi_1, \xi_2,g_2; \xi,g_1; g_0)\in\mathbb{K}_3$, and $\Vec{C}=(\xi_1', \xi_2', \xi_3',g_3; \xi_1, \xi_2,g_2; \xi,g_1; g_0)\in\mathbb{K}_4$. As it is explained in \cite{Cueca_Zhu:2021}, we need to show that $p_2^*\omega$ is normalized, closed, and non-degenerate.  Clearly, $p_2^*\theta$ is normalized; therefore, so is $p_2^*\omega.$ Now, observe that
    \begin{equation*}
        (\delta p_2^*\theta)_{\Vec{C}}=\langle\xi_2,\theta^r_{g_1}\rangle-\langle \xi_2+Ad^*_{g_2}\xi_3',\theta^r_{g_1}\rangle+\langle\xi_3',\theta^r_{g_2g_1}\rangle-\langle\xi_2'+\xi_3',\theta^r_{g_2}\rangle+\langle\xi_2',\theta^r_{g_2}\rangle=0,
    \end{equation*}
    because Maurer-Cartan $1$-forms are multiplicative in the sense of \cite[Lemma 3.2]{wei:symp}. Thus, $(\delta-d)p_2^*\omega=0$. We are left to prove that $p_2\omega$ is non-degenerate. For $(\eta,v)\in T_{(0,1)}(\gg^*\times G)\cong\gg^*\oplus\gg$, we get
    \begin{equation*}
        \begin{split}
            \lambda(\eta,v)=&p_2^*\omega(T(\bar\sigma_2\circ\bar{\sigma}_1)(v),T\bar{\sigma}_0(\eta))-p_2^*\omega(T(\bar\sigma_2\circ\bar{\sigma}_0)(v),T\bar{\sigma}_1(\eta))+p_2^*\omega(T(\bar\sigma_1\circ\bar{\sigma}_0)(v),T\bar{\sigma}_2(\eta))\\
            =&-(p_2^*\omega)_{|(0,0,1;0,1;1)}((0,0,0;0,v;0),(0,\eta,0;0,0;0))=-\omega_{|(0,1)}(v,\eta)=\langle\eta,v\rangle.
        \end{split}
    \end{equation*}
    Lastly, we compute the image of the $3$-shifted symplectic form $p_2^*\omega$ under the van Est map. Let $(\xi;z)\in\gg^*\otimes\gg$, then 
\begin{align*}
    \Phi ^{0110}(p_2^*\theta)(\xi;z)& =(R_\xi J_z-J_zR_\xi)(p_2^*\theta) \\
        & =\Big{(}\dd{\tau}(J_zp_2^*\theta)(\tau\xi,1;1)\Big{)}-\Big{(}(R_\xi p_2^*\theta)_1(z)\Big{)} \\
        & =\dd{\tau}\Big((p_2^*\theta)_{(0,0,1;\tau\xi,1;1)}(0,0,0;0,0;z)-(p_2^*\theta)_{(\tau\xi,0,1;0,1;1)}(0,0,0;0,0;z)\\
        &\qquad \qquad\quad +(p_2^*\theta)_{(0,\tau\xi,1;0,1;1)}(0,0,0;0,z;0)\Big) \\
        & =\dd{\tau}\langle 0,\theta^r_1(0)\rangle-\langle0,\theta^r_1(0)\rangle+\langle\tau\xi,\theta^r_1(z)\rangle =\langle\xi,z\rangle
\end{align*}
    and all the other components of the van Est map are zero, thus proving the result.
\end{proof}

Observe that the $2$-form $p_2^*\omega$ does not live on the nerve of the groupoid $\g^*\rtimes G\rightrightarrows G$; therefore, it is a genuine feature of the Lie 2-group $\mathbb{K}_\bullet$. On the other hand, the symplectic NQ-manifold $(T^*[3]\g^*[1], \omega^{\inf}, \cL_{\delta_{CE}})$ is the target of the $4$-dimensional TQFT known as BF-theory, see e.g. \cite{cat:bf}. One expects that its global counterpart $(\mathbb{K}_\bullet, p_2^*\omega)$ would play a r\^ole in physics.

\subsection{The loop $2$-group and the image of the Segal form}\label{ss-Segal}

In this section, let $G$ be a Lie group whose Lie algebra $\gg$ carries a symmetric and non-degenerate pairing satisfying 
$$\langle[a,b],c\rangle+\langle b,[a,c]\rangle=0\quad\forall a,b,c\in\gg.$$
Let $L_{1}G$ be the based loop group\footnote{We avoid the more traditional notation $\Omega G$ as it may conflict with our notation above} of $G$, that is loops starting at the identity of $G$ and of Sobolev class $r$. Let $P_1G$ be the space of paths on $G$ of a fixed Sobolev class $r$ starting at the identity. In \cite[Thm 3.9]{Cueca_Zhu:2021}, it is shown that the loop group carries a strict Lie 2-group structure
    \begin{eqnarray}\label{Eq:LoopGp}
    \xymatrix{
    \mathbb{G}_\bullet= L_1G \ar[r]\ar@<1.0ex>[r]\ar@<-1.0ex>[r] & P_1G \ar@<0.5ex>[r]\ar@<-0.5ex>[r] & \ast
    }
    \end{eqnarray}
    with Lie $2$-algebra given by $\mathbbm{g}=(L_0\g\to P_0\g, \ell_1, [\cdot,\cdot])$ where $ P_0\g$ and $L_0\g$ denote respectively the spaces of paths and loops starting at $0$ in $\g$ and $\ell_1$ is the natural inclusion. Moreover, in the same Theorem, it is proven that the Segal $2$-form $\omega_S\in\Omega^2(L_1G)$
    \begin{eqnarray*}
        (\omega_S)_\lambda(a,b)=\int_{S^1}\Big{\langle}\dd{\tau}dL_{\lambda(\tau)}^{-1}a(\tau),dL_{\lambda(\tau)}^{-1}b(\tau)\Big{\rangle}d\tau, & \textnormal{for }\lambda\in L_1G, a,b\in T_\lambda L_1G,
    \end{eqnarray*}
    is a 2-shifted symplectic form on $\mathbbm{G}_\bullet$. In \cite[Thm 3.15]{Cueca_Zhu:2021}, it is shown that there is a Morita equivalence $ev_\bullet: \mathbbm{G}\to G$ that relates the symplectic form $\omega_S$ with the $2$-shifted symplectic form on $G$ given by 
    $$\omega_2=\frac{1}{2}\langle \pr_1^*\theta^l, \pr_2^*\theta^r\rangle\in\Omega^2(G^2)\quad \omega_1=-\frac{1}{12}\langle [\theta^l, \theta^l], \theta^l\rangle\in\Omega^3(G),$$
    where $\theta^l, \theta^r$ are the left and right Maurer-Cartan 1-forms on $G$, respectively.
    
    An immediate consequence of \cite[Thm. 5.1]{baez:from} is that the Lie $2$-algebra $\mathbbm{g}$ is quasi-isomorphic to $\g$, where the quasi-isomorphism is given by 
    \begin{equation}\label{eq:ev}
        ev_1:P_0\g\to \g, \quad ev_1(\alpha)=\alpha(1)\quad \text{for}\quad \alpha\in P_0\g.
    \end{equation}
    
In the sequel, we evaluate the van Est map $\Phi$ constructed in Theorem \ref{Theo:CommutesWithDiffs} at $\omega_S$, and express its relation with  $\langle\cdot,\cdot\rangle$. 
\begin{prop}\label{Prop:Segal}
 Let $(\mathbb{G}_\bullet, \omega_S)$ be the $2$-shifted symplectic Lie 2-group \eqref{Eq:LoopGp}, let $ev_1$ be the map \eqref{eq:ev}, and let $\Phi$ the van Est map. Then  
    $$\Phi(\omega_S)=ev_1^*\langle\cdot,\cdot\rangle.$$
\end{prop}
\begin{proof}
    For $\alpha,\beta\in P_0\gg$, writing $\overline{1}\in L_1G$ for the constant loop,
\begin{align}\label{Eq:Segal1}
    (\Phi^{0020}\omega_S)(\alpha,\beta) & =J_\alpha J_\beta\omega_S+J_\beta J_\alpha\omega_S =(\omega_S)_{\overline{1}}(T\sigma_1(\alpha),T\sigma_0(\beta))+(\omega_S)_{\overline{1}}(T\sigma_1(\beta),T\sigma_0(\alpha)).
\end{align}
By definition, (see  \cite[Eq.~(3.6)]{Cueca_Zhu:2021}), for $\gamma\in P_1G$,
\begin{eqnarray*}
    (\sigma_0\gamma)(\tau)=\begin{cases}
        \gamma(3\tau) & \tau\in[0,\frac{1}{3}] \\
        \gamma(1) & \tau\in[\frac{1}{3},\frac{2}{3}] \\
        \gamma(3-3\tau) & \tau\in[\frac{2}{3},1],
    \end{cases} & (\sigma_1\gamma)(\tau)=\begin{cases}
        1 & \tau\in[0,\frac{1}{3}] \\
        \gamma(3\tau-1) & \tau\in[\frac{1}{3},\frac{2}{3}] \\
        \gamma(3-3\tau) & \tau\in[\frac{2}{3},1];
    \end{cases}
\end{eqnarray*}
therefore, continuing the computation of Eq.~\eqref{Eq:Segal1},
\begin{align*}
    (\omega_S)_{\overline{1}}(T\sigma_1(\alpha),T\sigma_0(\beta)) & =\int_{1/3}^{2/3}\Big{\langle}\frac{d}{dt}\alpha(3\tau-1),\beta(1)\Big{\rangle}d\tau+\int_{2/3}^{1}\Big{\langle}\frac{d}{dt}\alpha(3-3\tau),\beta(3-3\tau)\Big{\rangle}d\tau \\
        & =\int_{0}^{1}\frac{d}{dt}\langle\alpha(\tau),\beta(1)\rangle d\tau-\int_{0}^{1}\Big{\langle}\frac{d}{dt}\alpha(\tau),\beta(\tau)\Big{\rangle}d\tau\\
        &=\langle\alpha(1),\beta(1)\rangle-\int_{0}^{1}\Big{\langle}\frac{d}{dt}\alpha(\tau),\beta(\tau)\Big{\rangle}d\tau .
\end{align*}
Here, we used a simple substitution and the observation that $\frac{d}{d\tau}\langle\alpha,b\rangle=\langle\frac{d}{d\tau}\alpha,b\rangle$ for all $b\in\gg$. Thus, integrating by parts, 
\begin{align*}
    (\Phi^{0020}\omega_S)(\alpha,\beta) & =\langle\alpha(1),\beta(1)\rangle-\int_{0}^{1}\Big{\langle}\frac{d}{dt}\alpha(\tau),\beta(\tau)\Big{\rangle}d\tau+\langle\beta(1),\alpha(1)\rangle-\int_{0}^{1}\Big{\langle}\frac{d}{dt}\beta(\tau),\alpha(\tau)\Big{\rangle}d\tau \\
        & =2\langle\alpha(1),\beta(1)\rangle-\int_{0}^{1}\frac{d}{dt}\langle\alpha(\tau),\beta(\tau)\rangle d\tau \\
        &=2\langle\alpha(1),\beta(1)\rangle-\langle\alpha(1),\beta(1)\rangle=\langle\alpha(1),\beta(1)\rangle,
\end{align*}
as claimed. 
\end{proof}

Following ideas outlined in \cite{sev:some}, it was shown in \cite{Cueca_Zhu:2021} that the $2$-shifted symplectic Lie $2$-group $(\mathbb{G}_\bullet, \omega_S)$ appears as the reduced phase space of the Hamiltonian Chern-Simons theory on a disk with fixed boundary conditions. This procedure parallels the construction of the symplectic groupoid of a Poisson manifold obtained in \cite{cat:sym}. Proposition~\ref{Prop:Segal} shows that
$(\mathbb{G}_\bullet, \omega_S)$ is an integration of the degree $2$ symplectic $NQ$-manifold $(\underline{\g}[1], Q=\delta_{CE}, \omega=\langle\cdot,\cdot\rangle)$, further exhibiting a concrete procedure to differentiate. In other words, our van Est map defines a differentiation procedure that produces a Courant algebroid over a point. We summarize this result and its relationship to \cite[Thm 3.15]{Cueca_Zhu:2021} using the following schematic diagram
\begin{equation*}
    \xymatrix{\big(\mathbbm{g}, ev_1^*\langle\cdot,\cdot\rangle\big)\ar[r]^{ev_1}&(\g, \langle\cdot,\cdot\rangle)\\
    (\mathbb{G}_\bullet, \omega_S)\ar[u]^{\Phi}\ar[r]_{ev_\bullet}&(G,\omega_2+\omega_1).\ar[u]_{VE}}
\end{equation*}

\appendix
\section{The homotopies}\label{sec:ap}

In this appendix, we give details on the homotopy that appears in the body of the article, as well as an alternative homotopy. The alternative homotopy has the drawback that it has a harder formula and is harder to work with; however, it has the advantage that it verifies a further compatibility condition with the algebra structure as required by \cite{gls:pert}.

\subsection{The main homotopy}

Let $h_0^*:C^k(\A_{p+1};S^l\H_{p+1}^*)\longrightarrow C^k(\A_{p};S^l\H_{p}^*)$ be given by pulling back along the maps in Eq.'s~\eqref{Eq:ContrHomotopiesF} and~\eqref{Eq:ContrHomotopiesH}. More precisely, if $\omega\in C^k(\A_{p+1};S^l\H_{p+1}^*)$, define 
\begin{align*}
    (h_0^*\omega)_{\vec{\gamma}}(\mathbb{X};\mathbb{Y}) & :=\omega_{\bar{\sigma}_0(\vec{\gamma})}(h^{\A}_0(X_1, \vec{\gamma}),...,h^{\A}_0(X_{k},\vec{\gamma});h^{\H}_0(Y_1, \vec{\gamma}),...,h^{\H}_0(Y_l, \vec{\gamma})),
\end{align*}
for $\vec{\gamma}\in W_pG$, $\mathbb{X}=((X_1,\vec{\gamma}),...,(X_{k}, \vec{\gamma}))$ and $\mathbb{Y}=((Y_1,\vec{\gamma})...,(Y_l,\vec{\gamma}))$, where $(X_m, \vec{\gamma})\in\A_{p}$ for each $1\leq m\leq k$ and $(Y_n,\vec{\gamma})\in\H_{p}$ for each $1\leq n\leq l$. The following result states that, thus defined, $h_0^*$ is a contracting homotopy.

\begin{lem}\label{lemma:h0}
     $h_0^*$ is a homotopy between the identity and the zero map with respect to the differential $\partial$, that is, 
$$h_0^*\circ\partial+\partial\circ h_0^*=\textnormal{Id}.$$
\end{lem}
\begin{proof}
Recall the simplicial structure defined in Proposition \ref{prop:A}. Then, computing, we get on the one hand that
\begin{align*}
    h^\A_0(\partial_0^\cA(X,\vec{\gamma}))=&h^\cA_0(Td_1X\star 0_{\pr_p\varepsilon_p(\vec{\gamma})},\ \partial_0\vec{\gamma})
    =(\hat{\sigma}_0\big(Td_1X\star 0_{\pr_p\varepsilon_p(\vec{\gamma})}\star 0^{-1}_{s_0\pr_{p-1}\partial_0\vec{\gamma}}\big),\ \bar{\sigma}_0\partial_0\vec{\gamma})\\
    =& \big(T(s_0\circ d_1)(X)\star 0^{-1}_{s_0\pr_p\vec{\gamma}},\ \bar{\sigma}_0\partial_0\vec{\gamma}\big),
\end{align*}
while, whenever $i>0$, we get
\begin{align*}
    h^\A_0(\partial_i^\cA(X,\vec{\gamma}))=&h^\cA_0(Td_{i+1}X, \partial_i\vec{\gamma})
    =(\hat{\sigma}_0\big(Td_{i+1}X\star 0^{-1}_{s_0\pr_{p-1}\partial_i\vec{\gamma}}\big),\ \bar{\sigma}_0\partial_i\vec{\gamma})\\
    =& \big(T(s_0\circ d_{i+1})(X)\star 0^{-1}_{s_0d_{i+1}s_0\pr_p\vec{\gamma}}, \ \bar{\sigma}_0\partial_i\vec{\gamma}\big).
\end{align*}
On the other hand, for $i>0$,
\begin{align*}
    \partial^\cA_i(h^{\A}_0(X,\vec{\gamma}))=&\partial^\cA_i\big(\hat{\sigma}_0(X\star 0_{s_0\pr_p(\vec{\gamma})}^{-1}), \bar{\sigma}_0(\vec{\gamma})\big)\\
    =&\Big((Td_{i+1})\circ \hat{\sigma}_0(X\star 0_{s_0\pr_p(\vec{\gamma})}^{-1}), \partial_i\bar{\sigma}_0\vec{\gamma}\Big)= h^\A_0(\partial_{i-1}^\cA(X,\vec{\gamma})),
\end{align*}
and for $i=0$,
\begin{align*}
    \partial^\cA_0(h^{\A}_0(X,\vec{\gamma}))=&\partial^\cA_0\big(\hat{\sigma}_0(X\star 0_{s_0\pr_p(\vec{\gamma})}^{-1}), \bar{\sigma}_0(\vec{\gamma})\big)\\
    =&\Big((Td_1)\circ \hat{\sigma}_0(X\star 0_{s_0\pr_p(\vec{\gamma})}^{-1})\star 0_{\pr_{p+1}\varepsilon_{p+1}\bar\sigma_0\vec{\gamma}},\ \partial_0\bar{\sigma}_0\vec{\gamma}\Big)
    =(X,\vec{\gamma}).
\end{align*}
It is easy to see that the same identities hold for $h_0^\cH$; therefore, we get the equality
\begin{align}
    h^*_0\circ\partial+\partial\circ h^*_0=h_0^*\circ\left(\sum_{i=0}^{p}(-1)^i\partial^*_i\right)+\left(\sum_{i=0}^{p+1}(-1)^i\partial^*_i\right)\circ h_0^*=\id. & \qedhere
\end{align}
\end{proof}

The following result is necessary to conclude that the operator $h$ in~\eqref{Eq:h} is indeed a homotopy between the identity and $\pi_\bullet^*\circ\iota_\bullet^*$. 

\begin{prop}\label{Prop:h'}
    $h':=h-h_0^*$, where $h$ is the homotopy operator~\eqref{Eq:h}, is a homotopy between zero and $\pi_\bullet^*\circ\iota_\bullet^*$ with respect to the differential $\partial$, that is, 
    $$\partial\circ h'+h'\circ\partial =-\pi_{p+1}^*\circ\iota_{p+1}^*.$$
\end{prop}
\begin{proof}
        Let $\omega\in \Gamma\left(\bigwedge^k\A^*_p\otimes\sym^l\H_p^*\right)$, $\vec{\gamma}\in W_{p+1}\G$, $\mathbb{X}=((X_1, \vec{\gamma}),...,(X_k, \vec{\gamma}))$ and $\mathbb{Y}=((Y_1,\vec{\gamma})...,(Y_l,\vec{\gamma}))$, where $(X_a,\vec{\gamma})\in\A_{p+1}$ for each $1\leq a\leq k$ and $(Y_b,\vec{\gamma})\in\H_{p+1}$ for each $1\leq b\leq l$.  Computing, 
    \begin{align}\label{Eq:dh'}
        (\partial h'\omega)_{\vec{\gamma}}(\mathbb{X};\mathbb{Y}) & =\sum_{m=0}^{p+1}(-1)^m(h'\omega)_{\partial_m\vec{\gamma}}(\partial^\cA_m\mathbb{X};\partial^\cH_m\mathbb{Y}) \nonumber \\
            & =\sum_{m=1}^{p+2}\sum_{n=1}^{p}(-1)^{m+n}\omega_{\bar{\sigma}_0^{n+1}(\bar{\partial}_{p+2-n}\circ...\circ\bar{\partial}_{p+1}\circ\bar{\partial}_m\vec{\gamma})}(h^\A_n(\partial^\cA_{m-1}\mathbb{X});h^\H_n(\partial^\cH_{m-1}\mathbb{Y}))
    \end{align}
    and
    \begin{align}\label{Eq:h'd}
        (h'\partial\omega)_{\vec{\gamma}}(\mathbb{X};\mathbb{Y}) & =\sum_{n=1}^{p+1}(-1)^n(\partial\omega)_{\sigma_0^{n+1}(\partial_{p+2-n}\circ...\circ\partial_{p+1}\vec{\gamma})}(h^\A_n(\mathbb{X});h^\H_n(\mathbb{Y})) \nonumber \\
             & =\sum_{n=1}^{p+1}\sum_{m=1}^{p+3}(-1)^{m+n}\omega_{\bar{\partial}_m\bar{\sigma}_0^{n+1}(\bar{\partial}_{p+3-n}\circ...\circ\bar{\partial}_{p+2}\vec{\gamma})}(\partial^\cA_{m-1}h^\A_n(\mathbb{X});\partial^\cH_{m-1}h^\H_n(\mathbb{Y})).
    \end{align}
    By displaying the terms of the Eq.'s~\eqref{Eq:dh'} and~\eqref{Eq:h'd} in the Tables 1. and 2., we explain how mostly all terms cancel one another, leaving behind only 
    \begin{align}\label{Eq:lastTerm}
        (-1)^{(p+3)+(p+1)}(\omega)_{\bar{\partial}_{p+3}\bar{\sigma}_0^{p+2}(\bar{\partial}_{2}\circ...\circ\bar{\partial}_{p+2}\vec{\gamma})}(\partial^\cA_{p+2}h^\A_{p+1}(\mathbb{X});\partial_{p+2}^\cH h^\H_{p+1}(\mathbb{Y})) & =-(\pi_{p+1}^*\circ\iota_{p+1}^*\omega)_{\vec{\gamma}}(\mathbb{X};\mathbb{Y}).
    \end{align}
    In the tables below, the Roman numbers indicate elimination schemes, and the subindices indicate which pairs cancel one another. Each elimination scheme follows from the simplicial identities as follows:
    \begin{enumerate}
        \item Elimination Scheme I --- In the table for Eq.~\eqref{Eq:dh'}, for fixed $1\leq n\leq p$, $m=p+2-n+\alpha$ cancels out with $m=p+3-n+2\alpha$ for each $\alpha\geq 0$. Notice that, whenever $n$ is odd $\alpha\leq (n-1)/2$; whereas, whenever $n$ is even $\alpha\leq n/2-1$. We denote these pairs by I$_{n,\alpha}$.
        \item Elimination Scheme II --- In the table for Eq.~\eqref{Eq:h'd}, for fixed $1\leq n\leq p+1$, $m=1+2\alpha$ cancels out with $m=2+2\alpha$ for each $0\leq\alpha\leq (n-1)/2$ whenever $n$ is odd and for each $0\leq\alpha\leq n/2-1$ whenever $n$ is even. We denote these pairs by II$_{n,\alpha}$.
        \item Elimination Scheme III --- In the table for Eq.~\eqref{Eq:h'd}, for any fixed even $1<n\leq p+1$, the term with coordinates $(n+1,n)$ cancels out with the term with coordinates $(p+3,n-1)$. We denote these pairs by III$_{p+3-n,p+4-n,...,p+1,p+2}$.
        \item Elimination Scheme IV --- For a fixed even $1<n\leq p$, the terms in the lowest row of the $n$-th columns of the tables for Eq.~\eqref{Eq:dh'} and Eq.~\eqref{Eq:h'd} cancel each other out. We denote these pairs by IV$_{p+2-n,p+3-n,...,p+1,p+2}$.
        \item Elimination Scheme V --- For any fixed $1\leq n\leq p$, the remaining terms in the $n$-th column of the table for Eq.~\eqref{Eq:dh'} cancel each other in order with the remaining terms in the $n$-th column of the table for Eq.~\eqref{Eq:h'd}. We denote these pairs by V$_{m,p+3-n,p+4-n,...,p+1,p+2}$ with $1\leq m\leq p+1-n$.
    \end{enumerate}
\begin{center}
     \textbf{Table 1.} Table for Eq.~\ref{Eq:dh'}
     \begin{tabular}{l|c|c|c|c|c|c|c|}
      $\mathbf{m\backslash n}$ & \textbf{1} & \textbf{2} & \textbf{3} & $\mathbf{\cdots}$ & $\mathbf{p-2}$ & $\mathbf{p-1}$ & $\mathbf{p}$ \\
      \hline
      \textbf{1} & V$_{\scriptscriptstyle 1,p+2}$ & V$_{\scriptscriptstyle 1,p+1,p+2}$ & V$_{\scriptscriptstyle 1,p,p+1,p+2}$ & & V$_{\scriptscriptstyle 1,5,6,...,p+1,p+2}$ & V$_{\scriptscriptstyle 1,4,5,...,p+1,p+2}$ & V$_{\scriptscriptstyle 1,3,4,...,p+1,p+2}$ \\
      \hline
      \textbf{2} & V$_{\scriptscriptstyle 2,p+2}$ & V$_{\scriptscriptstyle 2,p+1,p+2}$ & V$_{\scriptscriptstyle 2,p,p+1,p+2}$ & &  V$_{\scriptscriptstyle 2,5,6,...,p+1,p+2}$  & V$_{\scriptscriptstyle 2,4,5,...,p+1,p+2}$ & I$_{\scriptscriptstyle p,0}$ \\
      \hline
      \textbf{3} & V$_{\scriptscriptstyle 3,p+2}$ & V$_{\scriptscriptstyle 3,p+1,p+2}$ & V$_{\scriptscriptstyle 3,p,p+1,p+2}$ & &  V$_{\scriptscriptstyle 3,5,6,...,p+1,p+2}$  & I$_{\scriptscriptstyle p-1,0}$ & I$_{\scriptscriptstyle p,0}$ \\
      \hline
      \textbf{4} & V$_{\scriptscriptstyle 4,p+2}$ & V$_{\scriptscriptstyle 4,p+1,p+2}$ & V$_{\scriptscriptstyle 4,p,p+1,p+2}$ & &  I$_{\scriptscriptstyle p-2,0}$  & I$_{\scriptscriptstyle p-1,0}$ & I$_{\scriptscriptstyle p,1}$ \\
      \hline
      $\mathbf{\vdots}$ \\
      \hline 
      $\mathbf{p-1}$ & V$_{\scriptscriptstyle p-1,p+2}$ & V$_{\scriptscriptstyle p-1,p+1,p+2}$ & I$_{3,0}$ & &  I$_{p-2,p/2-3}$  & I$_{p-1,p/2-2}$ & I$_{p,p/2-2}$ \\
      \hline
      $\mathbf{p}$ & V$_{p,p+2}$ & I$_{2,0}$ & I$_{3,0}$ & &  I$_{p-2,p/2-2}$  & I$_{p-1,p/2-2}$ & I$_{p,p/2-1}$ \\
      \hline
      $\mathbf{p+1}$ & I$_{1,0}$ & I$_{2,0}$ & I$_{3,1}$ & &  I$_{p-2,p/2-2}$  & I$_{p-1,p/2-1}$ & I$_{p,p/2-1}$ \\
      \hline
      $\mathbf{p+2}$ & I$_{1,0}$ & IV$_{\scriptscriptstyle p,p+1,p+2}$ & I$_{3,1}$ & & IV$_{\scriptscriptstyle 4,5,,...,p+1,p+2}$ & I$_{\scriptscriptstyle p-1,p/2-1}$ & IV$_{\scriptscriptstyle 2,3,...,p+1,p+2}$ \\
      \hline
    \end{tabular} 
  \end{center}     
\begin{center}
    \textbf{Table 2.} Table for Eq.~\ref{Eq:h'd}
    \begin{tabular}{l|c|c|c|c|c|c|c|}
      $\mathbf{m\backslash n}$ & \textbf{1} & \textbf{2} & \textbf{3} & $\mathbf{\cdots}$ & $\mathbf{p-1}$ & $\mathbf{p}$ & $\mathbf{p+1}$ \\
      \hline
      \textbf{1} & II$_{1,0}$ & II$_{2,0}$ & II$_{3,0}$ & & II$_{p-1,0}$ & II$_{p,0}$ & II$_{p+1,0}$ \\
      \hline
      \textbf{2} & II$_{1,0}$ & II$_{2,0}$ & II$_{3,0}$ & &  II$_{p-1,0}$  & II$_{p,0}$ & II$_{p+1,0}$ \\
      \hline
      \textbf{3} & V$_{\scriptscriptstyle 1,p+2}$ & III$_{\scriptscriptstyle p+1,p+2}$ & II$_{3,1}$ & &  II$_{p-1,1}$  & II$_{p,1}$ & II$_{p+1,1}$ \\
      \hline
      \textbf{4} & V$_{\scriptscriptstyle 2,p+2}$ & V$_{\scriptscriptstyle 1,p+1,p+2}$ & II$_{\scriptscriptstyle 3,1}$ & &  II$_{\scriptscriptstyle p-1,1}$  & II$_{\scriptstyle p,1}$ & II$_{\scriptstyle p+1,1}$ \\
      \hline
      $\mathbf{\vdots}$ \\
      \hline 
      $\mathbf{p}$ & V$_{\scriptscriptstyle p-2,p+2}$ & V$_{\scriptscriptstyle p-3,p+1,p+2}$ & V$_{\scriptscriptstyle p-4,p,p+1,p+2}$ & &  II$_{\scriptscriptstyle p-2,p/2-1}$  & II$_{\scriptscriptstyle p-1,p/2-1}$ & II$_{\scriptscriptstyle p+1,p/2-1}$ \\
      \hline
      $\mathbf{p+1}$ & V$_{\scriptscriptstyle p-1,p+2}$ & V$_{\scriptscriptstyle p-2,p+1,p+2}$ & V$_{\scriptscriptstyle p-3,p,p+1,p+2}$ & &  V$_{\scriptscriptstyle 1,4,5,...,p+1,p+2}$ & III$_{\scriptscriptstyle 3,4,,...,p+1,p+2}$ & II$_{\scriptscriptstyle p+1,p/2}$ \\
      \hline
      $\mathbf{p+2}$ & V$_{\scriptscriptstyle p,p+2}$ & V$_{\scriptscriptstyle p-1,p+1,p+2}$ & V$_{\scriptscriptstyle p-2,p,p+1,p+2}$ & &  V$_{\scriptscriptstyle 2,4,5,...,p+1,p+2}$ & V$_{\scriptscriptstyle 1,3,4,...,p+1,p+2}$ & II$_{\scriptscriptstyle p+1,p/2}$ \\
      \hline
      $\mathbf{p+3}$ & III$_{\scriptscriptstyle p+1,p+2}$ & IV$_{\scriptscriptstyle p,p+1,p+2}$ & III$_{\scriptscriptstyle p-1,p,p+1,p+2}$ & & III$_{\scriptscriptstyle 3,4,,...,p+1,p+2}$ & IV$_{\scriptscriptstyle 2,3,...,p+1,p+2}$ & \eqref{Eq:lastTerm} \\
      \hline
    \end{tabular} 
\end{center}


    These elimination schemes are deduced using the following two consequences of the simplicial identities. First, for each $1\leq n\leq p$, 
    \begin{align}\label{Eq:reordFace}
        \bar{\partial}_{p+2-n}\circ\bar{\partial}_{p+3-n}\circ ...\circ\bar{\partial}_{p}\circ\bar{\partial}_{p+1}\circ\bar{\partial}_{m} & =\begin{cases}
            \bar{\partial}_{m}\circ\bar{\partial}_{p+3-n}\circ...\circ\bar{\partial}_{p+2} & \textnormal{if }1\leq m<p+2-n \\
            \bar{\partial}_{p+2-n}\circ ...\circ\bar{\partial}_{p+1}\circ\bar{\partial}_{p+2} & \textnormal{if }p+2-n\leq m\leq p+2.
        \end{cases}
    \end{align}
    Next, for each $1\leq n\leq p+1$,
    \begin{align}\label{Eq:reordDeg}
        \bar{\partial}_{m}\circ\bar{\sigma}_0^{n+1} & =\begin{cases}
            \bar{\sigma}_0^{n} & \textnormal{if }1\leq m\leq n+1 \\
            \bar{\sigma}_0^{n+1}\circ\bar{\partial}_{m-(n+1)} & \textnormal{if }n+2\leq m\leq p+3.
        \end{cases}
    \end{align}
    Note that, Eq.'s~\eqref{Eq:reordFace} and~\eqref{Eq:reordDeg} immediately imply that each pair of terms in each of the elimination schemes belongs to the same fiber; moreover, note that in each said pair, the signs are necessarily opposed to one another. Thus, the result holds at the base level and we are left to prove that the arguments coincide. 
    
    For the $k$ many $\A$ coordinates, because Eq.'s~\eqref{Eq:reordFace} and~\eqref{Eq:reordDeg} hold for every simplicial manifold, they hold in particular for the nerve of the strict Lie 2-algebra. According to Eq.~\eqref{Eq:BFaceMaps}, the only map that one needs to be careful with is $\partial_1$ in the Elimination Scheme II, but if $1\leq n\leq p+1$, for any $1\leq a\leq k$, then $$\bar{\sigma}_0^n(\bar{\partial}_{p+3-n}\circ\bar{\partial}_{p+4-n}\circ ...\circ\bar{\partial}_{p+1}\circ\bar{\partial}_{p+2}(\vec{\gamma}))=(\mathbbm{1}_{p+3};\mathbbm{1}_{p+2};\cdots)$$. As for the $l$ many $\H$ coordinates, the result follows easily as each $Y_b\star 0_{\sigma_0(\gamma^{p+1})}^{-1}$ can be rewritten as $\hat{\sigma}_1^{p+1}(y_b)$ for some $y_b\in\hh$. Lastly, Eq.~\eqref{Eq:lastTerm} holds as $\bar{\sigma}_0^{p+2}(\bar{\partial}_{1}\circ...\circ\bar{\partial}_{p+2}\vec{\gamma})=\bar{\sigma}_0^{p+2}(\ast)=\mathbbm{1}\in W_{p+2}G.$
\end{proof}

\subsection{The alternative homotopy}

In the sequel, we lay out an alternative homotopy operator. Instead of using $\bar{\sigma}_0$ on the base, we define $(\eta_0)_\bullet:W_\bullet G\longrightarrow W_{\bullet+1} G$ inductively as
\begin{align*}\label{Eq:Eta0}
        (\eta_0)_{0}(g) & =(s_0g;g^{-1}), \qquad (\eta_0)_{1}(\gamma;g) =\big((\gamma,\overline{\gamma});\overline{\gamma}^{-1}\vJoin s_0(d_0(\gamma)\vJoin g);(d_0(\gamma)\vJoin g)^{-1}), \\
        (\eta_0)_p(\vec{\gamma}) & =\left((\gamma^p,\overline{\gamma^p_1\Join...\Join\gamma^p_p});\big(\gamma^{p-1},\overline{\gamma^{p-1}_1\Join...\Join\gamma^{p-1}_p}\Join\textnormal{pr}_2(\eta_0)_{1}(\gamma^p_1;d_1\gamma^{p-1}_1)\big);\bar{\partial}_0(\eta_0)_{p-1}(\partial_0\vec{\gamma})\right), \nonumber
\end{align*}
where $g\in G_0$, $\vec{\gamma}=\gamma\in G_1$, $(\gamma^p;...;\gamma^0)\in W_p\G$ with $\gamma^k=(\gamma^k_1,...,\gamma^k_k)\in G_k$ and the bar denotes the inverse in the groupoid. 
\begin{lem}\label{lemma:AltHomotopy0}
      $\eta_0^*:C^\infty(W_{p+1} G)\longrightarrow C^\infty(W_{p} G)$ given by $\eta_0^*f=(-1)^{p}f\circ(\eta_0)_{p-1}$ defines a homotopy between the identity and the zero map with respect to $\partial$, that is, 
$$\eta_0^*\circ\partial+\partial\circ\eta_0^*=\textnormal{Id}.$$
\end{lem}
\begin{proof}
    The result follows from noticing that
    $$\partial_{p+1}\circ(\eta_0)_{p}=\textnormal{Id}\quad\text{and}\quad (\eta_0)_{p-1}\circ\partial_n=\partial_n\circ(\eta_0)_{p}$$
    for all $0\leq n\leq p$. This follows from the identities we list in the sequence. For $\gamma^k\in G_k$ as above, let $\gamma^k_{[]}:=\gamma^k_1\Join...\Join\gamma^k_k\in G_1$.
    \begin{itemize}
        \item For $n=0$, the identity follows from 
        \begin{align*}
        \scriptstyle
           \overline{\gamma^p}_{[]}\vJoin\Big{(}\overline{\gamma^{p-1}}_{[]}\Join\big{(}(\overline{\gamma^p_1})^{-1}\vJoin s_0d_1(\gamma^p_2\vJoin\gamma^{p-1}_1)\big{)}\Big{)}& =  \Big{(}\overline{(d_0\gamma^p)}_{[]}\Join\overline{\gamma^{p}_1}\Big{)}\vJoin\Big{(}\overline{\gamma^{p-1}}_{[]}\Join\big{(}(\overline{\gamma^p_1})^{-1}\vJoin s_0d_1(\gamma^p_2\vJoin\gamma^{p-1}_1)\big{)}\Big{)} \\
           & =\Big{(}\overline{(d_0\gamma^p)}_{[]}\vJoin\overline{\gamma^{p-1}}_{[]}\Big{)}\Join\Big{(}\overline{\gamma^{p}_1}\vJoin(\overline{\gamma^p_1})^{-1}\vJoin s_0d_1(\gamma^p_2\vJoin\gamma^{p-1}_1)\Big{)}\\
           &= \overline{(d_0\gamma^p\vJoin\gamma^{p-1})}_{[]}.
        \end{align*}
        \item For $n=1$, the identity follows from $$\overline{(d_0\gamma^{p-1}\vJoin\gamma^{p-2})}_{[]}\Join\big{(}(\overline{\gamma^p_1\Join\gamma^p_2})^{-1}\vJoin s_0d_1(\gamma^p_3\vJoin\gamma^{p-1}_2\vJoin\gamma^{p-2}_1)\big{)}$$
        being $$\Big{(}\overline{\gamma^{p-1}}_{[]}\Join\big{(}(\overline{\gamma^p_1})^{-1}\vJoin s_0d_1(\gamma^p_2\vJoin\gamma^{p-1}_1)\big{)}\Big{)}\vJoin\Big{(}\overline{\gamma^{p-2}}_{[]}\Join\big{(}(\overline{\gamma^p_2\vJoin\gamma^{p-1}_1})^{-1}\vJoin s_0d_1(\gamma^p_3\vJoin\gamma^{p-1}_2\vJoin\gamma^{p-2}_1)\big{)}\Big{)},$$ which is true  by the exchange law.
        \item For $n=2$, first consider 
        \begin{align*}
            \partial_2\circ(\eta_0)_{2}(\vec{\gamma}) & =\partial_2\big((\gamma^2,\overline{\gamma^2}_{[]});(\gamma^1,\overline{\gamma^1}\Join((\overline{\gamma^2_1})^{-1}\vJoin s_0d_1(\gamma^2_2\vJoin\gamma^1)));d_0\circ(\eta_0)_{1}\circ\partial_0(\vec{\gamma})\big) \\
                & =\big{(}d_2(\gamma^2,\overline{\gamma^2}_{[]});d_1(\gamma^1,\overline{\gamma^1}\Join((\overline{\gamma^2_1})^{-1}\vJoin s_0d_1(\gamma^2_2\vJoin\gamma^1)));d_1d_0\circ(\eta_0)_{1}\circ\partial_0(\vec{\gamma})\big{)} \\
                & =\big{(}(\gamma^2_1,\overline{\gamma^2_1});(\overline{\gamma^2_1})^{-1}\vJoin s_0d_1(\gamma^2_2\vJoin\gamma^1);d_1(\gamma^2_2\vJoin\gamma^1)^{-1}\big{)}\\
                &=(\eta_0)_{1}(d_2\gamma^2;d_1\gamma^1)= (\eta_0)_{1}(\partial_2\vec{\gamma}),
        \end{align*}
        and then
        \begin{align*}
            \partial_2\circ(\eta_0)_{1}(\vec{\gamma}) & =\partial_2\big((\gamma,\overline{\gamma});\overline{\gamma}^{-1}\vJoin s_0(d_0(\gamma)\vJoin g);(d_0(\gamma)\vJoin g)^{-1}) \\
                & =\big{(}d_2(\gamma,\overline{\gamma});d_1(\overline{\gamma}^{-1}\vJoin s_0(d_0(\gamma)\vJoin g))\big{)}=(\gamma;g).
        \end{align*}
    \end{itemize}  
     Due to the simplicial identities the rest follows by induction as  $$\bar{\partial}_0\circ(\eta_0)_{p-2}\circ\partial_0\circ\partial_n=\partial_{n-2}\circ\bar{\partial}_0\circ(\eta_0)_{p-1}\circ\partial_0$$ is equivalent to $$\eta_0^{(p-2)}\circ\partial_{n-1}=\partial_{n-1}\circ\eta_0^{(p-1)}$$ for $p>2$ and $2\leq n\leq p+1$; whereas, for $p> 1$, $\partial_p\circ\bar{\partial}_{0}\circ(\eta_0)_{p}\partial_1=\textnormal{Id}$ is equivalent to $\partial_{p+1}\circ(\eta_0)_{p}=\textnormal{Id}$.     
\end{proof}
The homotopy operator $\eta_0^*$ can be extended to a homotopy between the identity and $\pi_\bullet^*\circ\iota_\bullet^*$ in two steps. First, for $1\leq i\leq p$, define $(\eta_i)_p: W_p G\to W_{p+1}G$ by 
\begin{eqnarray}\label{Eq:eta i}
    (\eta_i)_{p}:=(\eta_0)_{p}\circ ...\circ(\eta_0)_{p-i}\circ\partial_{p+1-i}\circ\partial_{p+2-i}\circ...\circ\partial_{p-1}\circ\partial_{p}.
\end{eqnarray}
Next, as we did with $h$~\eqref{Eq:h}, define $\eta:=\sum_{i=0}^{p}(-1)^i\eta_i^*$, keeping in mind there is an extra sign coming from the pullback as in Lemma~\ref{lemma:AltHomotopy0}. Now, to extend the homotopy operator to the double-complex $\cC^{p,q}$ in~\eqref{Eq:AuxDoubleCxForms}, we first set  
    \begin{eqnarray*}
        \eta^{\A}_0:\A_p \longrightarrow  \A_{p+1}, &  \eta^{\A}_0(X,\vec{\gamma}):=\big((X,\overline{X}_{[]}), (\eta_0)_p(\vec{\gamma})\big),  \\ \label{Eq:AltHomotopiesH}
       \eta^{\H}_0:\H_p \longrightarrow \H_{p+1}, &  \eta^{\H}_0(Y,\vec{\gamma}):=\big(Ts_1(Y), (\eta_0)_p(\vec{\gamma})\big).
    \end{eqnarray*}
Here, if $X=(X_0,...,X_p)\in T_{s_0\pr_p\vec{\gamma}}G_{p+1}$, $X_{[]}:=X_0\Join ...\Join X_p$.
Lastly, $\eta_i^\cA$ and $\eta_i^\cH$ are defined, mutatis mutandis, by Eq.~\eqref{Eq:eta i} and $$\eta:\cC^{p+1,q}\longrightarrow \cC^{p,q}$$
is defined by the alternating sum $\sum_{i=0}^{p}(-1)^i\eta_i^*$, where the pullbacks are taken as in Eq.~\eqref{Eq:hFinal}.

\begin{prop}\label{Prop:AltHomotopy}
    Thus defined, $\eta:\cC^{p+1,q}\longrightarrow \cC^{p,q}$ is a homotopy between the identity and $\pi_\bullet^*\circ\iota_\bullet^*$ with respect to the differential $\partial$. Furthermore, $\eta$ is a $(\pi_\bullet^*\circ\iota_\bullet^*)$-derivation with respect to the natural product in the Weil algebra, that is,
    $$\eta(a\smile b)=\eta (a)\smile(\pi_\bullet^*\circ\iota_\bullet^*)(b)+(-1)^{\vert a\vert}a\smile\eta (b),$$
    and verifies the `side conditions' $\eta\circ\eta = 0$, $\iota_\bullet^*\circ\eta = 0$.    
\end{prop}
This result follows analogously to the groupoid case (cf. \cite[Prop.~3.3]{Meinrenken_Li-Bland:2015}) and, as a consequence of the refined Perturbation lemma of Gugenheim-Lambe-Stasheff \cite{gls:pert}, $\iota_\bullet^*\circ(1+\delta\circ\eta)^{-1}$ is a morphism of differential algebras. Because the formula it yields is necessarily unique, the following result follows.
\begin{thm}
    The van Est map $\Phi$ is an algebra morphism, and under the connectedness hypotheses of Theorem~\ref{Theo:IndIsosCochains}, it descends to algebra isomorphisms in cohomology.
\end{thm}

\bibliographystyle{plain}
\bibliography{refs}
\end{document}